\documentclass{amsart}
\usepackage{pstricks,pst-node,amssymb,epsfig}
\usepackage{pstcol}             

\title{Total positivity, Grassmannians, and networks}

\author{Alexander Postnikov}

\address{Department of Mathematics,
         M.I.T., Cambridge, MA 02139}

\email{apost@math.mit.edu}
\urladdr{http://math.mit.edu/\mytilde apost/}


\thanks{This work is supported in part by NSF CAREER Award DMS-0504629.}

\date{September 26, 2006}

\psset{unit=1pt,arrowsize=4pt,linewidth=.7pt}
\psset{linecolor=blue}
\newgray{grayish}{.90}
\newrgbcolor{emgreen}{0 .5 0}
\def\vblack(#1,#2)#3{\cnode*[linecolor=black](#1,#2){3}{#3}}
\def\vwhite(#1,#2)#3{\cnode[linecolor=black,fillcolor=white,fillstyle=solid](#1,#2){3}{#3}}

\theoremstyle{plain}
\newtheorem{theorem}{Theorem}[section]
\newtheorem{proposition}[theorem]{Proposition}
\newtheorem{corollary}[theorem]{Corollary}
\newtheorem{conjecture}[theorem]{Conjecture}
\newtheorem{lemma}[theorem]{Lemma}

\theoremstyle{definition}
\newtheorem{definition}[theorem]{Definition}
\newtheorem{example}[theorem]{Example}

\theoremstyle{remark}
\newtheorem{remark}[theorem]{Remark}

\numberwithin{equation}{section}
\numberwithin{figure}{section}

\def\R{\mathbb{R}}
\def\Z{\mathbb{Z}}
\def\RP{\mathbb{RP}}

\def\M{\mathcal{M}}

\def\P{\mathcal{P}}
\def\I{\mathcal{I}}
\def\Le{\hbox{\rotatedown{$\Gamma$}}}
\def\miscell{Misce\reflectbox{\sc ll}aneous}

\def\GL{\mathrm{GL}}
\def\Mat{\mathrm{Mat}}

\def\mytilde{\kern-.015in\hbox{\lower.03in\hbox{\~{}}}\kern-.01in}
\def\Tail{\mathrm{Tail}}

\def\Id{\mathit{Id}}
\def\Mes{\mathit{Meas}}
\def\Image{\mathrm{Image}}
\def\MR{\mathrm{MR}}

\def\tnn{\mathrm{tnn}}
\def\tp{\mathrm{tp}}
\def\wind{\mathit{wind}}
\def\form{\mathrm{form}}
\def\Net{\mathit{Net}}
\def\Cyc{\mathit{Cyc}}
\def\inv{\mathit{inv}}
\def\xing{\mathit{xing}}

\def\col{\mathit{col}}
\def\OO{\mathcal{O}}

\begin{document}

\begin{abstract}
The aim of this paper is to discuss a relationship between
total positivity and planar directed networks.
We show that the inverse boundary problem for these networks is naturally
linked with the study of the totally nonnegative Grassmannian.
We investigate its cell decomposition, where
the cells are the totally nonnegative parts of the matroid strata.
The boundary measurements of networks give parametrizations
of the cells.  We present several different combinatorial descriptions 
of the cells, study the partial order on the cells, and describe how
they are glued to each other.  
\end{abstract}

\maketitle

\tableofcontents

\section{Introduction}

A {\it totally positive\/} matrix is a matrix with all positive minors.  We
extend this classical notion, as follows.  Define the 
{\it totally nonnegative Grassmannian\/} $Gr_{kn}^\tnn$ as the set of elements 
in the Grassmannian $Gr_{kn}(\R)$  with all
nonnegative Pl\"ucker coordinates.  The classical set of totally positive
matrices is embedded into $Gr_{kn}^\tnn$ as an open subset.
The intersections $S_\M^\tnn$ of the {\it matroid strata\/} with $Gr_{kn}^\tnn$,
which we call the {\it nonnegative Grassmann cells,} give an interesting subdivision of
the totally nonnegative Grassmannian.  These intersections are actually cells 
(unlike the matroid strata that might have a nontrivial geometric structure)
and they form a CW-complex.  Conjecturally, this is a
regular CW-complex and the closures of the cells are homeomorphic to balls.
Fomin-Zelevinsky's {\it double Bruhat cells\/} (for type $A$) are included into
$Gr_{kn}^\tnn$ as certain special nonnegative Grassmann cells $S_\M^\tnn$.
Note that the subdivision of $Gr_{kn}^\tnn$ into the cells $S_\M^\tnn$
is a finer subdivision than the Schubert decomposition.

Our ``elementary'' approach agrees with Lusztig's general theory
total positivity~\cite{Lusz1, Lusz2, Lusz3} and with the cellular decomposition
of the nonnegative part of $G/P$ conjectured by Lusztig and proved by
Rietsch~\cite{Riet1, Riet2}.

Another main object of the paper is a planar directed {\it network\/} 
drawn inside a disk with the boundary vertices $b_1,\dots,b_n$ (and some number of 
internal vertices) and with positive weights $x_e$ on the edges.
We assume that $k$ boundary vertices $b_i$ are sources and the remaining 
$(n-k)$ boundary vertices $b_j$ are sinks.  We allow the sources and sinks to be 
interlaced with each other in any fashion.  For an acyclic network, we define the 
{\it boundary measurements\/} as
$M_{ij} = \sum_{P\,:\,b_i\to b_j} \prod_{e\in P} x_e$,
where the sum is over all directed paths $P$ in the network from 
a source $b_i$ to a sink $b_j$, and the product is over all edges $e$ of $P$.
If a network has directed cycles, we introduce the sign $(-1)^{\wind(P)}$
into the weight of $P$, where $\wind(P)$ is the winding index that
counts the number of full $360^\circ$ turns the path $P$ makes.
We show that the power series for $M_{ij}$ (which might be infinite
if $G$ has directed cycles) always sums to a {\it subtraction-free\/} 
rational expression. 

We discuss the {\it inverse boundary problem\/} for such planar directed
networks.  In other words, we are interested in the information about networks
that can be recovered from all boundary measurements $M_{ij}$.  We characterize
all possible collections of the measurements, describe all transformations of
networks that preserve the measurements, and show how to reconstruct a
network from the measurements (up to these transformations).
Our work on this problem is parallel to results of
Curtis-Ingerman-Morrow \cite{CIM, Inger, CM} on the inverse problem for (undirected)
transistor networks.

The inverse boundary problem for directed networks has deep connections with
total positivity.  The collection of all boundary measurements $M_{ij}$ of a
network can be encoded as a certain element of the Grassmannian $Gr_{kn}$.
This gives the boundary measurement map $\Mes:\{networks\}\to Gr_{kn}$.  We
show that the image of the map $\Mes$ is exactly the totally nonnegative
Grassmannian $Gr_{kn}^\tnn$.  Moreover, the image of the set of networks with
fixed combinatorial structure given by a graph $G$ (with arbitrary positive
weights on the edges of $G$) is a certain nonnegative Grassmann cell
$S_\M^\tnn$.  This gives the map from graphs $G$ to the set of Grassmann cells.
If the graph $G$ is reduced (that is minimal in a certain sense) then the map
$\Mes$ induces a rational subtraction-free parametrization of the corresponding
cell $S_\M^\tnn$.

For each cell $S_\M^\tnn$ we describe one particular graph $G$ given by a
$\Le$-diagram.  These $\Le$-diagrams are fillings of Young diagrams of shape
$\lambda$ with $0$'s and $1$'s that satisfy certain $\Le$-property.  The shape
$\lambda$ corresponds to the Schubert cell $\Omega_\lambda$ that contains
$S_\M^\tnn$.  The $\Le$-diagrams have interesting combinatorial properties.

There are several types of transformations of networks that preserve the
boundary measurements.  First of all, there are quite obvious rescaling of the
edge weights $x_e$ at each internal vertex, which we call the {\it gauge
transformations.} Then there are transformations that allows up to switch
directions of edges in the network.  We can easily transform any network into a
special form (called a perfect network) and then color the vertices into two
colors according to some rule.   We prove that the boundary measurement map
$\Mes$ is invariant under switching directions of edges that preserve colors of
vertices.  Thus the boundary measurement map can now be defined for {\it
undirected\/} planar networks with vertices colored in two colors.  We call them
{\it plabic networks\/}  (abbreviation for ``planar bicolored'').  Finally, there
are several {\it moves\/} (that is local structure transformations) of plabic
networks that preserve the boundary measurement map.  We prove that any two
networks with the same boundary measurements can be obtained from each other by
a sequence of these transformations.

Essentially, the only nontrivial transformation of networks is a certain {\it
square move.}  This move can be related to cluster transformation from
Fomin-Zelevinsky theory of {\it cluster algebras\/} \cite{FZ2, FZ3, FZ4, BFZ2}.
It is a variant of the {\it octahedron recurrence\/} in a disguised form.

We show how to transform each plabic graph into a {\it reduced\/} graph.  
We define {\it trips\/} in such graph as directed paths in these
(undirected) graphs that connect boundary vertices $b_i$ and obey certain ``rules of
the road.''   The trips give the {\it decorated trip permutation\/} of the boundary
vertices.  We show that any two reduced plabic graphs as related by the moves
and correspond to the same Grassmann cell $S_\M^\tnn$ if and only if they have
the same decorated trip permutation.  Thus the cells $S_\M^\tnn$ are in
one-to-one correspondence with decorated permutations.

Plabic graphs can be thought of as generalized wiring diagrams, which are
graphical representations of reduced decompositions in the symmetric group.  The
moves of plabic graphs are analogues of the {\it Coxeter moves\/} of wiring
diagrams.  Plabic graphs also generalize Fomin-Zelevinsky's {\it double wiring
diagrams\/} \cite{FZ1}.

We also define {\it alternating strand diagrams,} which are in bijection with
plabic graphs.  These diagrams consist of $n$ directed strands that connect $n$
points on a circle and intersect with each other inside the circle in an
alternating fashion.  
Scott~\cite{Sco1, Sco2} used our alternating strand diagrams to 
study Leclerc-Zelevinsky's \cite{LZ} quasi-commuting families of quantum minors
and cluster algebra on the Grassmannian.

We discuss the partial order on the cells $S_\M^\tnn$ by containment of their
closures and describe it in terms of decorated permutations.  We call this order
the {\it circular Bruhat order\/} because it reminds the usual (strong) Bruhat
order on the symmetric group. Actually, the usual Bruhat order is a certain
interval in the circular Bruhat order.  

We use our network parametrizations of the cells to describe how they are
glued to each other.  The gluing of a cell $S_\M^\tnn$ to the lower
dimensional cells inside of its closure $\overline{S_\M^\tnn}$ is described by
sending some of the edge weights $x_e$ to $0$.  Thus, for the cell
$S_\M^\tnn$ associated with a graph $G$, the lower dimensional cells in its
closure are associated with subgraphs $H\subseteq G$
obtained from $G$ by removing some edges.  In a sense, this is an analogue of
the statement that, for a Weyl group element with a reduced decomposition
$w=s_{i_1}\cdots s_{i_l}$, all elements below $w$ in the Bruhat order are
obtained by taking subwords in the reduced decomposition.

For each plabic graph $G$ associated with a cell $S_\M^\tnn$, we describe 
a different parametrization of the cell by a certain subset of 
the Pl\"ucker coordinates.
This parametrization is related to the boundary measurement parametrization by
the {\it chamber anzatz\/} and a certain {\it twist map\/} $S_\M^\tnn/T \to
S_\M^\tnn/T$, where $T$ is the ``positive torus'' $T=\R_{>0}^n$ acting on
$Gr_{kn}^\tnn$.  This construction is analogous to a similar construction of
Berenstein-Fomin-Zelevinsky~\cite{BFZ, FZ1} for double Bruhat cells.  In our
setup, instead of chambers in (double) wiring diagrams, we work with regions of
plabic graphs.

As an application, we obtain a description of Berenstein-Zelevinsky's {\it
string cones\/} and polytopes~\cite{BZ1, BZ2} (of type $A$) as {\it
tropicalizations\/} of the boundary measurements $M_{ij}$.  Integer lattice
points in these polytopes count the {\it Littlewood-Richardson coefficients}.
This explains the combinatorial description of the string cones from our earlier
work~\cite{GlPo} and the rule for the Littlewood-Richardson coefficients.

Our construction produces several different combinatorial objects associated with
the cells $S_\M^\tnn$.  We give explicit  bijections between all
these objects.  Here is the (incomplete) list of various objects: totally
nonnegative Grassmann cells, totally nonnegative matroids, $\Le$-diagrams,
decorated permutations, circular chains, move-equivalence classes of (reduced)
plabic graph, move-equivalence classes of alternating strand diagrams.

We also construct bijections between $\Le$-diagrams 
and other combinatorial objects such as permutations of
a certain kind, rook placements on skew Young diagrams, etc.  Williams~\cite{W1},
Steingrimsson-Williams~\cite{S-W}, and Corteel-Williams~\cite{CW} obtained
several enumerative results on $\Le$-diagrams and related objects,
and studied their combinatorial properties.

\medskip
Throughout the paper we use the following notation
$[n]:=\{1,\dots,n\}$ and $[k,l]:=\{k,k+1,\dots,l\}$.
The word ``network'' means a graph (directed or undirected)
together with some weights assigned to edges or faces of the graph.

\medskip

Many results of this paper were obtained in 2001. 
Some results were announced by Williams in \cite[Sect.~2--3]{W1}
and \cite[Appendix]{W2}.

\medskip
\noindent
\textsc{Acknowledgments:}
I would like to thank 
(in alphabetical order)
Sergey Fomin,
Alexander Goncharov,
Alberto Gr\"unbaum,
Xuhua He,
David Ingerman,
Allen Knutson,
George Lusztig,
James Propp,
Konni Rietsch,
Joshua Scott,
Michael Shapiro,
Richard Stanley,
Bernd Sturmfels,
Dylan Thurston,
Lauren Williams,
and Andrei Zelevinsky
for helpful conversations.

\section{Grassmannian}
\label{sec:grassmannian}

In this section we review some classical facts about Grassmannians,
their stratifications, and matroids. 
For more details, see~\cite{Fult}.
\medskip

\subsection{Schubert cells}
\label{ssec:schubert_cells}

For $n\geq k\geq 0$, let the {\it Grassmannian\/} $Gr_{kn}$ be the 
manifold of $k$-dimensional subspaces $V\subset \R^n$.  It can be presented
as the quotient $Gr_{kn} = GL_k\backslash \Mat^*_{kn}$, where 
$\Mat^*_{kn}$ is the space of real $k\times n$-matrices of rank $k$.
Here we assume that the subspace $V$ associated with a $k\times n$-matrix 
$A$ is spanned by the row vectors of $A$.

Recall that a {\it partition\/} $\lambda=(\lambda_1,\dots,\lambda_k)$ 
is a weakly decreasing sequence of nonnegative integers.  It is
graphically represented by its Young diagram which is the collection of boxes
with indexes $(i,j)$ such that $1\leq i\leq k$, $1\leq j\leq \lambda_i$
arranged on the plane in the same fashion as one would arrange matrix elements.

There is a cellular decomposition of the Grassmannian $Gr_{kn}$ into 
a disjoint union of Schubert cells $\Omega_\lambda$ indexed by partitions 
$\lambda \subseteq (n-k)^k$ 
whose Young diagrams fit 
inside the $k\times(n-k)$-rectangle $(n-k)^k$, 
that is $n-k\geq \lambda_1\geq \cdots \geq \lambda_k\geq 0$.

\begin{figure}[ht]
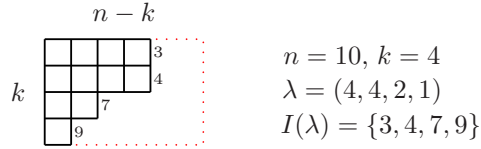

\pspicture(0,0)(60,56)
\psline[linecolor=black](0,40)(40,40)
\psline[linecolor=black](0,30)(40,30)
\psline[linecolor=black](0,20)(40,20)
\psline[linecolor=black](0,10)(20,10)
\psline[linecolor=black](0,0)(10,0)
\psline[linecolor=black](0,40)(0,0)
\psline[linecolor=black](10,40)(10,0)
\psline[linecolor=black](20,40)(20,10)
\psline[linecolor=black](30,40)(30,20)
\psline[linecolor=black](40,40)(40,20)
\psline[linecolor=red,linestyle=dotted](42,40)(60,40)(60,0)(12,0)
\rput(5,5){}
\rput(5,15){}
\rput(5,25){}
\rput(5,35){}
\rput(15,15){}
\rput(15,25){}
\rput(15,35){}
\rput(25,25){}
\rput(25,35){}
\rput(35,25){}
\rput(35,35){}
\rput(45,35){}
\rput(43,35){\tiny $3$}
\rput(43,25){\tiny $4$}
\rput(23,15){\tiny $7$}
\rput(13,5){\tiny $9$}
\rput(-10,20){$k$}
\rput(30,50){$n-k$}
\rput(120,33){$n=10$, $k=4$}
\rput(120,20){$\lambda=(4,4,2,1)$}
\rput(127,7){$I(\lambda)=\{3,4,7,9\}$}
\endpspicture
\caption{A Young diagram $\lambda$ 
and the corresponding subset $I(\lambda)$}
\label{fig:lambda_I}
\end{figure}

The partitions $\lambda \subseteq (n-k)^k$ are in one-to-one correspondence
with $k$-element subsets $I\subset [n]$.
The boundary of the Young diagram of such partition $\lambda$ 
forms a lattice path from the the upper-right corner to the lower-left corner
of the rectangle $(n-k)^k$. Let us label the $n$ steps in this path
by the numbers $1,\dots,n$ consecutively, and define $I = I(\lambda)$ 
as set of labels of $k$ vertical steps in the path.
The inverse map $I=\{i_1<\dots<i_k\}\mapsto \lambda$ is given by 
$\lambda_j = |[i_j,n]\setminus I|$, for $j=1,\dots,k$.
As an example, Figure~\ref{fig:lambda_I} shows 
a Young diagram of shape $\lambda=(4,4,2,1)\subseteq 6^4$
that corresponds to the subset $I(\lambda) = \{3,4,7,9\}\subseteq [10]$.

For $\lambda\subseteq (n-k)^k$,
the {\it Schubert cell\/} $\Omega_\lambda$ in $Gr_{kn}$ is defined 
as the set of $k$-dimensional subspaces $V\subset \R^n$ with
prescribed dimensions of intersections with the elements of
the opposite coordinate flag: 
$$
\Omega_\lambda :=\{V\in Gr_{kn}\mid \dim(V\cap \left<e_i,\dots,e_n\right>) = 
|I(\lambda)\cap [i,n]|,\text{ for } i=1,\dots,n\},
$$ 
where 
$\left<e_i,\dots,e_n\right>$ is the linear span of coordinate vectors.

The decomposition of $Gr_{rk}$ into Schubert cells can be described
using {\it Gaussian elimination.}  We can think of points in
$Gr_{kn}$ as nondegenerate $k\times n$-matrices modulo row
operations.  Recall that according to Gaussian elimination 
any nondegenerate $k\times n$-matrix can be transformed by 
row operations to the canonical matrix in {\it echelon form,}
that is a matrix $A$ such that $a_{1,i_1}=\cdots=a_{k,i_k}=1$
for some $I=\{i_1<\cdots <i_k\}\subset [n]$,
and all entries of $A$ to the left of these $1$'s and in same columns 
as the $1$'s are zero.  
In other words, matrices in echelon form are representatives
of the left cosets in $GL_k\backslash \Mat^*_{kn} = Gr_{kn}$.
Let us also say that such echelon matrix $A$ is in {\it $I$-echelon form\/} 
if we want to specify that the $1$'s are located in the column set $I$.
For example, a matrix in $\{3,4,7,9\}$-echelon form, for $n=10$ and $k=4$,
looks like
$$
A=\begin{pmatrix}
0 & 0 & 1 & 0 & * & * & 0 & * & 0 & * \\
0 & 0 & 0 & 1 & * & * & 0 & * & 0 & * \\
0 & 0 & 0 & 0 & 0 & 0 & 1 & * & 0 & * \\
0 & 0 & 0 & 0 & 0 & 0 & 0 & 0 & 1 & * 
\end{pmatrix}
$$
where ``$*$'' stands for any element of $\R$.

The Schubert cell $\Omega_\lambda$ is exactly the set of elements 
in the Grassmannian $Gr_{kn}$ that are represented by
matrices $A$ in $I$-echelon form, where $I=I(\lambda)$.
If we remove the columns with indices $i\in I$ from $A$,
i.e., the columns with the $1$'s,
and reflect the result with respect the vertical axis, 
the pattern formed by the $*$'s is exactly
the Young diagram of shape $\lambda$.  
So an $I$-echelon matrix has exactly $|\lambda|:=\lambda_1+\cdots+\lambda_k$ 
such $*$'s, which can be any elements in $\R$.
This shows that the Schubert cell $\Omega_\lambda$
homeomorphic to $\R^{|\lambda|}$.
Thus the Grassmannian $Gr_{kn}$ has the disjoint decomposition
$$
Gr_{kn} = \bigcup_{\lambda\subseteq (n-k)^k} \Omega_\lambda \simeq
\bigcup_{\lambda\subseteq (n-k)^k} \R^{|\lambda|}.
$$

For example, the $*$'s in the $\{3,4,7,9\}$-echelon form above
correspond to boxes of the Young diagram of shape $\lambda = (4,4,2,1)$.
Thus the Schubert cell $\Omega_{(4,4,2,1)}$, whose elements
are represented by matrices in 
$\{3,4,7,9\}$-echelon form, is isomorphic to $\R^{|\lambda|}=\R^{11}$.

\subsection{Pl\"ucker coordinates}

For a $k\times n$-matrix $A$ and a $k$-element subset $I\subset [n]$,
let $A_I$ denote the $k\times k$-submatrix of $A$ in the column set $I$,
and let $\Delta_I(A):=\det(A_I)$ denote the corresponding {\it maximal minor\/}
of $A$. If we multiply $A$ by $B\in GL_k$ on the left, all minors
$\Delta_I(A)$ are rescaled by the same factor $\det(B)$.
If $A=(a_{ij})$ is in $I$-echelon form then $A_I = \Id_k$ and 
$a_{ij} = \pm \Delta_{(I\setminus\{i\})\cup\{j\}}(A)$.
Thus the $\Delta_I$ form projective coordinates 
on the Grassmannian $Gr_{kn}$, called the {\it Pl\"ucker coordinates,}
and the map $A\mapsto (\Delta_I)$ induces the {\it Pl\"ucker embedding\/} 
$Gr_{kn}\hookrightarrow \RP^{\binom{n}{k} -1}$
of the Grassmannian into the projective space.
The image  of the Grassmannian $Gr_{kn}$ under the Pl\"ucker embedding
is the algebraic subvariety in $\RP^{\binom{n}{k} -1}$ given by 
the {\it Grassmann-Pl\"ucker relations:} 
\begin{equation*}
\Delta_{(i_1,\dots,i_k)}\cdot \Delta_{(j_1,\dots,j_k)} 
= 
\sum_{s=1}^k  \,
\Delta_{(j_s,i_2,\dots,i_k)}\cdot 
\Delta_{(j_1,\dots,j_{s-1},i_1,j_{s+1},\dots,j_k)},
\label{eq:Gr-Pl}
\end{equation*}
for any $i_1,\dots,i_k,j_1,\dots,j_k\in [n]$.
Here we assume that $\Delta_{(i_1,\dots,i_k)}$
(labelled by an ordered sequence rather than a subset) 
equals to $\Delta_{\{i_1,\dots,i_k\}}$ if $i_1<\cdots <i_k$
and $\Delta_{(i_1,\dots,i_k)} = (-1)^{\mathrm{sign}(w)} 
\Delta_{(i_{w(1)},\dots,i_{w(k)})}$ for any permutation $w\in S_k$.

\subsection{Matroid strata}
\label{ssec:Matroid_strata}

An element in the Grassmannian $Gr_{kn}$ can also be
understood as a collection of $n$ vectors $v_1,\dots,v_n\in\R^k$ spanning the
space $\R^k$, modulo the simultaneous action of $GL_k$ on the vectors.  The
vectors $v_i$ are the columns of a $k\times n$-matrix $A$ that represents the
element of the Grassmannian.

Recall that a {\it matroid\/} of {\it rank\/} $k$ on the set $[n]$ 
is a nonempty collection $\M\subseteq \binom{[n]}{k}$ of 
$k$-element subsets in $[n]$, called {\it bases} of $\M$, that satisfies the {\it exchange axiom:}

\begin{center}
For any $I,J\in \M$ and $i\in I$ there exists $j\in J$ such that 
$(I\setminus\{i\})\cup\{j\}\in M$.
\end{center}

An element $V\in Gr_{kn}$ of the Grassmannian represented by a 
$k\times n$-matrix $A$ gives the matroid $\M_V$ whose bases are 
the $k$-subsets $I\subset[n]$
such that $\Delta_I(A)\ne 0$, or equivalently, $I$ is a base of $\M_V$
whenever $\{v_i\mid i\in I\}$ is a basis of $\R^k$.
This collection of bases satisfies the exchange axiom because, if
the left-hand side in a Grassmann-Pl\"ucker relation is nonzero,
then at least one term in the right-hand side is nonzero.

The Grassmannian $Gr_{kn}$ has a subdivision into 
{\it matroid strata,} also known as {\it Gelfand-Serganova strata,}
$S_\M$ labelled by some matroids $\M$:
$$
S_\M:=\{V\in Gr_{kn}\mid \M_V = \M\} 
$$
In other words, the elements of the stratum $S_\M$ are represented 
by matrices $A$ such that $\Delta_I(A) \ne 0$ if and only if 
$I\in \M$.
The matroids $\M$ with nonempty strata $S_\M$ are called 
{\it realizable\/} over $\R$.  The geometrical structure of 
the matroid strata $S_\M$ can be highly
nontrivial.  Mn\"ev~\cite{Mnev} showed that they 
can be as complicated as essentially any algebraic variety.

Note that, for an element $V$  of the Schubert cell $\Omega_\lambda$,
the subset $I(\lambda)$ is exactly the lexicographically minimal base
of the matroid $\M_V$.  This fact it transparent when $V$ is represented
by a matrix in $I$-echelon form.  In other words, the Schubert cells 
can also be described as
$$
\Omega_\lambda = \{V\in Gr_{kn}\mid I(\lambda) 
\text{ is the lexicographically minimal base of } \M_V\}. 
$$
This implies that the decomposition of $Gr_{kn}$ into matroid
strata $S_\M$ is a finer subdivision than the decomposition into 
Schubert cells $\Omega_\lambda$.

The Schubert decomposition depends on a choice of ordering of 
the coordinates in $\R^n$.  The symmetric group $S_n$ acts on $\R^n$ by 
permutations of the coordinates.  For a permutation $w\in S_n$, 
let $\Omega_\lambda^w :=w(\Omega_\lambda)$ be the permuted Schubert
cell.  In other words, the cell $\Omega_\lambda^w$ is the set of 
elements $V\in Gr_{kn}$ such that $I(\lambda)$ is the lexicographically
minimal base of $\M_V$ with respect to the total order
$w(1)<w(2)<\cdots < w(n)$ of the set $[n]$.

\begin{remark}
\label{rec:common_ref_factorial}
The decomposition of the Grassmannian $Gr_{kn}$ into matroid strata
$S_\M$ is the common refinement of the $n!$ 
permuted Schubert decompositions $Gr_{kn}=
\bigcup_{\lambda\subseteq (n-k)^k} \Omega_\lambda^w$,
for $w\in S_n$; see~\cite{GGMS}. 
Indeed, if we know the lexicographically minimal base in $\M_V$
with respect any total order on $[n]$ then we can determine
all bases of $\M_V$ and thus determine the matroid stratum containing 
the element $V\in Gr_{kn}$.
\end{remark}

It will be convenient for us use local affine coordinates on the Grassmannian.
Let us pick a $k$-subset $I\subset[n]$.
Let $A$ be a $k\times n$-matrix that represents an element  in $Gr_{kn}$
such that $\Delta_I(A)\ne 0$, that is $I$ is a base of the corresponding
matroid. 
Then $A'=(A_I)^{-1} A$ is the unique representative the left coset of 
$GL_k\cdot A$ such that $A'_I$ is the identity matrix.  
Then matrix elements of $A'$ located in columns indexed $j\not\in I$
give local affine coordinates on $Gr_{kn}$.  In other words,
we have the rational isomorphism
$$
Gr_{kn}\setminus \{\Delta_I = 0\} \simeq \R^{k(n-k)}.
$$
In the case when $I$ is the lexicographically minimal base,
the matrix $A'$ is exactly the representative in echelon form.

\section{Totally nonnegative Grassmannian}
\label{set:tot_nonneg_Grass}

A matrix is called {\it totally positive\/} (resp., {\it totally
nonnegative\/}) if all its minors of all sizes are strictly positive (resp.,
nonnegative).  In this section we discuss analogues of these classical 
notions for the Grassmannian.  
\medskip

\begin{definition} 
Let us define the {\it totally nonnegative Grassmannian\/}
$Gr_{kn}^\tnn \subset Gr_{kn}$ as the quotient 
$Gr_{kn}^\tnn  = \GL_k^{+}\backslash \Mat^\tnn_{kn}$, where 
$\Mat^\tnn_{kn}$ is the set of real $k\times n$-matrices $A$ 
of rank $k$ with nonnegative {\it maximal\/} minors $\Delta_I(A)\geq 0$
and $\GL_k^+$ is the group of $k\times k$-matrices with
positive determinant.
The {\it totally positive Grassmannian\/}
$Gr_{kn}^\tp \subset Gr_{kn}^\tnn$ is the 
subset of $Gr_{kn}$ whose elements can be represented by
$k\times n$-matrices with strictly positive maximal minors $\Delta_I(A)>0$.
\end{definition}  

For example, the totally positive Grassmannian contains all
$k\times n$-matrices $A=(x_j^i)$ with $x_1<\cdots<x_n$, because any maximal
minor $\Delta_I(A)$ of such matrix is a positive Vandermonde determinant.
Clearly, $Gr_{kn}^\tp$ is an open subset in $Gr_{kn}$ and
$Gr_{kn}^\tnn$ is a closed subset in $Gr_{kn}$ of dimension
$k(n-k)=\dim Gr_{kn}$. 

\begin{definition}
Let us define {\it totally nonnegative Grassmann cells $S_\M^\tnn$} 
in $Gr_{kn}^\tnn$ as the intersections 
$S_\M^\tnn = S_\M\cap Gr_{kn}^\tnn$ of the matroid strata $S_\M$ 
with the totally nonnegative Grassmannian, i.e.,
$$
S_\M^\tnn = \{GL_k^+\cdot A \in Gr_{kn}^\tnn \mid \Delta_I(A)> 0
\textrm{ for } I\in\M, \textrm{ and } \Delta_I(A) = 0
\textrm{ for } I\not\in\M\}.
$$

Let us say that a matroid $\M$ is {\it totally nonnegative\/} 
if the cell $S_\M^\tnn$ is nonempty.
\end{definition}

Note that the totally positive Grassmannian $Gr_{kn}^\tp$ is just 
the top dimensional cell $S_\M^\tnn \subset Gr_{kn}^\tnn$, 
that is the cell corresponding to the complete matroid $\M=\binom{[n]}{k}$.

\begin{remark}
\label{rem:cyclic_symmetry}
Clearly, the notion of total positivity is not invariant under
permutations of the coordinates in $\R^n$,  
and the class of totally nonnegative matroids is not preserved 
under permutations of the elements.
This notion does however have some symmetries.
For a $k\times n$-matrix $A=(v_1,\dots,v_k)$ with the column vectors 
$v_i\in \R^k$,
let $A' = (v_2,\dots,v_n,(-1)^{k-1} v_1)$ be the matrix obtained from $A$ by
the cyclic shift of the columns and then multiplying the last column by
$(-1)^{k-1}$.  Note that $\Delta_I(A) = \Delta_{I'}(A')$ where
$I'$ is the cyclic shift of the subset $I$.
Thus $A$ is totally nonnegative (resp., totally positive)
if and only if $A'$ is totally nonnegative (resp., totally positive).
This gives an action of the cyclic group $\Z/n\Z$ on 
the sets $Gr_{kn}^\tnn$ and $Gr_{kn}^\mathrm{tp}$.
This also implies that cyclic shifts of elements in $[n]$
preserve the class of totally nonnegative matroids on $[n]$.
\end{remark}

\begin{example}  For $n=4$ and $k=2$, 
there are only three rank $2$ matroids on $[4]$ which are not 
totally nonnegative:
$\M = \{\{1,2\},\{2,3\},\{3,4\},\{1,4\}\}$,
$\M \cup\{\{1,3\}\}$,
$\M \cup\{\{2,4\}\}$.
This set of matroids is closed under cyclic shifts of $[4]$.
\end{example}

Interestingly, the totally nonnegative Grassmann cells $S_\M^\tnn$
have a much simpler geometric structure than the matroid strata $S_\M$.

\begin{theorem}  
Each totally nonnegative Grassmann cell $S_\M^\tnn$ is homeomorphic to an 
open ball of appropriate dimension.  
The decomposition of the totally nonnegative Grassmannian $Gr_{kn}^\tnn$
into the union of the cells $S_\M^\tnn$ is a CW-complex.
\end{theorem}

In Section~\ref{sec:Le-diagrams} we will explicitly construct a rational
parametrization for each cell $S_\M^\tnn$, i.e., an isomorphism between the
space $\R_{>0}^d$ and $S_\M^\tnn$; see Theorem~\ref{th:g_D}.  In
Section~\ref{sec:gluing} we will describe how the cells are glued to each
other.

This next conjecture follows a similar conjecture  
by Fomin-Zelevinsky~\cite{FZ1} on double Bruhat cells.

\begin{conjecture}
The CW-complex formed by the cells $S_\M^\tnn$ is regular.  The closure
of each cell is homeomorphic to a closed ball 
of appropriate dimension.
\end{conjecture}

According to Remark~\ref{rec:common_ref_factorial},
the matroid stratification of the Grassmannian is the common 
refinement of $n!$ Schubert decompositions.  For the totally
nonnegative part of the Grassmannian it is enough to take
just $n$ Schubert decompositions. 

\begin{theorem}
\label{th:cycles_are_enough}
The decomposition of $Gr_{kn}^\tnn$ into the cells $S_\M^\tnn$
is the common refinement of the $n$ Schubert decompositions
$Gr_{kn}^\tnn = 
\bigcup_{\lambda\subseteq (n-k)^k} (\Omega_\lambda^w \cap Gr_{kn}^\tnn)$,
where $w$ run over cyclic shifts $w:i\mapsto i+k\pmod n$, for $k\in [n]$.
\end{theorem}

This theorem will follow from Theorem~\ref{th:SM=circ_necklace}
in Section~\ref{sec:circ_bruhat}.

Lusztig~\cite{Lusz1, Lusz2, Lusz3} developed general theory of total positivity
for a reductive group $G$ using canonical bases.  He defined the totally
nonnegative part $(G/P)_{\geq 0}$ of any
generalized partial flag manifold $G/P$ and conjectured that it is made up
of cells.  This conjecture was proved by Rietsch~\cite{Riet1, Riet2}.
Marsh-Rietsch~\cite{M-R} gave a simpler proof and constructed 
parametrization of the totally nonnegative cells in $(G/B)_{\geq 0}$.
This general approach to totally positivity agrees with our ``elementary'' 
approach. 

\begin{theorem}  
\label{th:Rietsch=equal}
In case of the Grassmannian $Gr_{kn}$, Rietsch's cell decomposition
coincides with the decomposition of $Gr_{kn}^\tnn$ into the cells $S_\M^\tnn$.
\end{theorem}

I thank Xuhua He and Konni Rietsch for the following explanation.  According
to~\cite[Proposition 12.1]{M-R}, Rietsch cells are given by conditions
$\Delta_I>0$ and $\Delta_J=0$ for {\it some\/} minors.  Actually, the
paper~\cite{M-R} concerns with the case of $G/B$,  but the case of $G/P$ (which
includes the Grassmannian) can be obtained by applying the projection map
$G/B\to G/P$, as it was explained in~\cite{Riet1}.  It follows that our cell
decomposition of $Gr_{kn}^\tnn$ into the cells $S_\M^\tnn$ is a {\it
refinement\/} of Rietsch's cell decomposition.  In
Section~\ref{sec:Bruhat_intervals} we will construct a combinatorial bijection
between objects that label our cells and objects that label Rietsch's cells,
which will prove Theorem~\ref{th:Rietsch=equal}.

Let us show how total positivity on the Grassmannian is related to 
the classical notion of total positivity of matrices.
For a $k\times n$-matrix $A$ such that the square submatrix 
$A_{[k]}$ in the first $k$ columns is the identity matrix $A_{[k]}=\Id_k$,
define $\phi(A) = B$, where $B=(b_{ij})$ is the $k\times (n-k)$-matrix 
with entries $b_{ij} = (-1)^{k-j} a_{i+k,j}$:
$$
\phi:
\begin{pmatrix}
1  & \cdots & 0 & 0 & a_{1,k+1} & \cdots & a_{1n} \\
\vdots & \vdots & \ddots & \vdots & \vdots  & \ddots & \vdots  \\
0  & \cdots & 1 & 0 & a_{k-1,k+1} & \cdots & a_{k-1,n} \\
0  & \cdots & 0 & 1 & a_{k,k+1}  & \cdots & a_{kn} \\
\end{pmatrix}
\mapsto
\begin{pmatrix}
\pm a_{1,k+1} & \cdots & \pm a_{1n} \\
\vdots & \ddots & \vdots \\
- a_{k-1,k+1} & \cdots & - a_{k-1,n} \\
a_{k,k+1} & \cdots & a_{kn}
\end{pmatrix}
.
$$

Let $\Delta_{I,J}(B)$ denote the minor of matrix $B$ 
(not necessarily maximal) in the row set $I$ and the column set $J$.
By convention, we assume that $\Delta_{\emptyset,\emptyset}(B) = 1$.

\begin{lemma}
Suppose that $B= \phi(A)$.  
There is a correspondence between the maximal minors of $A$ and 
all minors of $B$ such that 
each maximal minor of $A$ equals to the corresponding 
minor of $B$.  Explicitly,
$\Delta_{I,J}(B) = \Delta_{([k]\setminus I)\cup \tilde J}(A)$,
where $\tilde J$ is obtained by increasing all elements in $J$ by $k$.
\end{lemma}

\begin{proof}
Exercise for the reader.
\end{proof}

Note that matrices $A$ with $A_{[k]}=\Id_k$ are representatives
(in echelon form) of elements of the top Schubert cell 
$\Omega_{(n-k)^k}\subset Gr_{kn}$, i.e., the set of elements in the
Grassmannian with nonzero first Pl\"ucker coordinate 
$\Delta_{[k]}\ne 0$.  Thus $\phi$ gives the isomorphism
$\phi:\Omega_{(n-k)^k}\to \Mat_{k,n-k}$.

\begin{proposition}
The map $\phi$ induces the isomorphism between $\Omega_{(n-k)^k}\cap
Gr_{kn}^\tnn$ and the set of classical totally nonnegative $k\times
(n-k)$-matrices, i.e., matrices with nonnegative minors of all sizes.  This map
induces the isomorphism between each totally nonnegative cell $S_\M^\tnn
\subset \Omega_{(n-k)^k}\cap Gr_{kn}^\tnn$ and a set of $k\times
(n-k)$-matrices given by prescribing some minors to be positive and the
remaining minors to be zero.  In particular, it induces the isomorphism between
the totally positive part $Gr_{kn}^\tp$ of the Grassmannian and the classical
set of all totally positive $k\times (n-k)$-matrices.  
\end{proposition}

\begin{remark}
\label{rem:FZ_double_Bruhat}
Fomin-Zelevinsky~\cite{FZ1} investigated the decomposition
of the totally nonnegative part of $GL_k$ into cells,
called the {\it double Bruhat cells}.  These cells are parametrized
by pairs of permutations in $S_k$.  The partial order by containment
of closures of the cells is isomorphic to the direct product of two
copies of the Bruhat order on $S_k$.
The map $\phi$ induces the isomorphism between the totally nonnegative
part of the Grassmannian $Gr_{k,2k}$ such that $\Delta_{[k]}\ne 0$
and $\Delta_{[k+1,2k]}\ne 0$ and the totally nonnegative part of
$GL_k$.  Moreover, it gives isomorphisms between the double Bruhat 
cells in $GL_k$ and {\it some\/} totally nonnegative cells
$S_\M^\tnn\subset Gr_{k,2k}^\tnn$;  namely, the cells
such that $[k], [k+1,2k]\in\M$. 
\end{remark}

In this paper we will extend Fomin-Zelevinsky's \cite{FZ1} results 
on (type $A$) double Bruhat cells to all totally nonnegative cells 
$S_\M^\tnn$ in the Grassmannian.
We will see that these cells have a rich combinatorial structure
and lead to new combinatorial objects.

\section{Planar networks}

\begin{definition}
\label{def:planar_networks}
A {\it planar directed graph\/} $G$ is a directed graph drawn inside a disk
(and considered modulo homotopy).
We allow $G$ to have loops and multiple edges.
We will assume that $G$ has $n$ {\it boundary vertices\/} on the boundary of
the disk labelled $b_1,\dots,b_n$ clockwise.  The remaining vertices, called
the {\it internal vertices,} are located strictly inside the disk.  We will
always assume that each boundary vertex $b_i$ is either a {\it source\/} 
or a {\it sink.}   Even if $b_i$ is an {\it isolated\/} boundary vertex,
i.e., a vertex not incident to any edges, we will assign $b_i$ to be
a source or a sink.
A {\it planar directed network\/} $N=(G,x)$ is a planar directed graph
$G$ as above together with {\it strictly positive\/} real weights
$x_e>0$ assigned to all edges $e$ of $G$.

For such network $N$, the {\it source set\/} $I\subset[n]$ and the {\it
sink set\/} $\bar I:=[n]\setminus I$ of $N$ are the sets such that $b_i$, $i\in
I$, are the sources of $N$ (among the boundary vertices) and the $b_j$, $j\in
\bar I$, are the boundary sinks. 
\end{definition}

If the network $N$ is acyclic, that is it does not have closed directed paths,
then, for any $i\in I$ and $j\in \bar I$, we define
the {\it boundary measurement\/} $M_{ij}$ as the finite sum
$$
M_{ij} := \sum_{P\,:\,b_i \to b_j} \prod_{e\in P} x_e ,
$$
where the sum is over all directed paths $P$ in $N$ from the boundary source
$b_i$ to the boundary sink $b_j$, and the product is over all edges $e$ in $P$.

If the network is not acyclic, we have to be more careful because the above sum
might be infinite.   We will need the following definition.

\begin{figure}[ht]
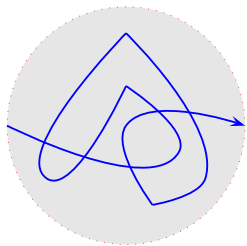

\pspicture(0,0)(100,100)
\pscircle[linecolor=red,linestyle=dotted,fillstyle=solid,fillcolor=grayish,linewidth=0](50,50){45}
\pscurve{-}(5,50)(70,40)(50,65)
\pscurve{-}(50,65)(20,30)(50,85)
\pscurve{-}(50,85)(80,30)(60,20)
\pscurve{->}(60,20)(50,50)(95,50)
\endpspicture
\caption{A path $P$ with the winding index $\wind(P)=-1$}
\label{fig:2}
\end{figure}

For a path $P$ from a boundary vertex $b_i$ to a boundary vertex $b_j$, we 
define its winding index, as follows.
We may assume that all edges of the network are given by smooth curves;
thus the path $P$ is given by a continuous piecewise-smooth curve.
We can slightly modify the path and smoothen it around each junction,
so that it is given by a smooth curve $f:[0,1]\to\R^2$,
and furthermore make the initial tangent vector $f'(0)$ to
have the same direction as the final tangent vector $f'(1)$.
We can now define the {\it winding index\/} $\wind(P)\in \Z$ of the path $P$ as
the signed number of full $360^\circ$ turns the tangent vector $f'(t)$ makes 
as we go from $b_i$ to $b_j$ (counting counterclockwise turns 
as positive);
see example in Figure~\ref{fig:2}.

Let us also give a recursive combinatorial definition of the
winding index for a path $P$ with vertices $v_1,v_2,\dots,v_l$ 
(where $v_1$ and $v_l$ are boundary vertices).  If the path $P$ 
has no self-intersections, i.e., all vertices $v_i$ in $P$ are different, 
then $\wind(P)=0$.  Suppose that $P$ has at least one self-intersection.  
Let us find a {\it cycle\/} in $P$, that is a segment $C$ with vertices 
$v_i$, $i\in[a,b]$, such that the vertices 
$v_a,v_{a+1},\dots,v_{b-1}$ are distinct and $v_a = v_b$.
Let $\tilde P$ be the path with vertices $v_1,\dots,v_a,v_{b+1},\dots,v_l$,
i.e., the path $P$ with erased cycle $C$.  We now have 
$$
\wind(P) = 
\left\{
\begin{array}{cl}
\wind(\tilde P) + 1 & \text{if $C$ is a counterclockwise cycle;}\\
\wind(\tilde P) - 1 & \text{if $C$ is a clockwise cycle.}
\end{array}
\right.
$$
We can now determine the winding index of $P$ by repeatedly erasing 
cycles in $P$ until we get a path without self-intersections.

\begin{remark}
Note that in general the number of cycles in $P$ is not well-defined.
Indeed the number of cycles that we need to erase until we get a path
without self-intersections may depend on the order in which we erase
the cycles.
However the number $\wind(P)= \#\{\textrm{counterclockwise cycles}\}
-\#\{\textrm{clockwise cycles}\}$ is a well-defined
invariant of a path in a planar graph.
For example, for the path shown on Figure~\ref{fig:2},
we can erase a counterclockwise cycle and then two
clockwise cycles.  On the other hand, for the same path, we can also erase
just one big clockwise cycle to get a path without self-intersections.
\end{remark} 

Let us now return to boundary measurements.
Let $N$ be a planar directed network as above, which is now allowed
to have cycles.  Let us assume for a moment that the weights $x_e$ 
of edges in $N$ are formal variables.
For a source $b_i$, $i\in I$, and a sink $b_j$, $j\in \bar I$,
we define the {\it formal boundary measurement\/} $M_{ij}^\form$ as 
the formal power series
\begin{equation}
\label{eq:Mij}
M_{ij}^\form:=
\sum_{P\,:\,b_i\to b_j} (-1)^{\wind(P)}\prod_{e\in P} x_e\,,
\end{equation}
where the sum is over all directed paths $P$ in $N$ from 
$b_i$ to $b_j$ and the product is over all edges $e$ in $P$.

Recall that a {\it subtraction-free rational expression\/} is an
expression with positive integer coefficients
that can be written with the operations of addition, multiplication, 
and division (but subtraction is strictly forbidden), or
equivalently, it is a quotient of two polynomial expressions
with positive coefficients.  For example, 
$\frac{x+y/x}{z^2+25y/(x+t)} = \frac{(x^2 + y)(x+t)}{xz^2(x+t) + 25xy}$ 
is subtraction-free.

\begin{lemma}
\label{lem:Mij}
The formal power series $M_{ij}^\form$ 
sums to a subtraction-free rational expression in the variables $x_e$.
\label{lem:Mij-subtraction-free}
\end{lemma}
This statement follows from
Proposition~\ref{prop:minors-exact} in the next section.

\begin{definition}
We can now define the {\it boundary measurements\/} $M_{ij}$ 
as the specializations of the formal boundary measurements $M_{ij}^\form$, written as 
subtraction-free expressions, when we assign the $x_e$ to be the 
positive real weights of edges $e$ in the network $N$.
Since the $M_{ij}^\form$ are
subtraction-free and the weights $x_e$ of edges in $N$ are strictly positive, we will 
never get a zero in denominator.  Thus the boundary measurements $M_{ij}$
are well-defined nonnegative real numbers for an arbitrary network.
\end{definition}

\begin{example}
For the network
\begin{center}
\pspicture(-50,-20)(50,25)
\cnode*[linecolor=black](-50,0){2}{B1}
\cnode*[linecolor=black](50,0){2}{B2}
\cnode*[linecolor=black](-20,0){2}{V1}
\cnode*[linecolor=black](20,0){2}{V2}
\nccurve[angleA=0,angleB=180]{->}{B1}{V1}
\aput*{:0}{$x$}
\nccurve[angleA=60,angleB=120]{->}{V1}{V2}
\aput*{:0}{$y$}
\nccurve[angleA=-60,angleB=-120]{<-}{V1}{V2}
\bput*{:0}{$z$}
\nccurve[angleA=0,angleB=180]{->}{V2}{B2}
\aput*{:0}{$t$}
\rput(-60,0){$b_1$}
\rput(60,0){$b_2$}
\endpspicture
\end{center}
we have $M_{12}^\form = xyt - x y z y t + x y z y z y t - \cdots 
= xyt/(1 + y z)$, which is a subtraction-free rational expression.
If all weights of edges are $x=y=z=t=1$, then 
the boundary measurement is $M_{12} = 1/(1+1) = 1/2$.
\end{example}

\smallskip
\noindent{\bf Inverse Boundary Problem.} \
{\it What information about a planar directed network can be recovered 
from the collection of boundary measurements $M_{ij}$?
How to recover this information?
Describe all possible collections of boundary measurements.
Describe transformation of networks that preserve 
the boundary measurements.
}
\smallskip

Let us describe the {\it gauge transformations\/} of the weights $x_e$.
Pick a collection of positive real numbers $t_v>0$, for each internal vertex 
$v$ in $N$; and also assume that $t_{b_i} = 1$ for each boundary vertex 
$b_i$.  Let $N'$ be the network with the same directed graph 
as the network $N$ and with the weights 
\begin{equation}
\label{eq:gauge_transformation}
x_e' = x_e \,t_{u} {t_{v}}^{-1},
\end{equation} 
for each directed edge $e=(u,v)$.
In other words, for each internal vertex $v$ we multiply by $t_v$
the weights of all edges outgoing from $v$,  divide by $t_v$
the weights of all edges incoming to $v$.  
Then the network $N'$ has the same boundary measurements as the network $N$.
Indeed, for a directed path $P$ between two boundary vertices and for an internal vertex $v$, 
we have to divide the weight $\prod_{e\in P} x_e$ of $P$ by $t_v$ every time when $P$ enters  $v$
and multiply it by $t_v$ every time when $P$ leaves $v$.

We will see that there are also some local structure transformations 
of networks that preserve the boundary measurements.

Let us now describe the set of all possible collections of boundary measurements.
For a network $N$ with $k$ boundary sources $b_i$, $i\in I$,
and $n-k$ boundary sinks $b_j$, $j\in \bar I$,
it will be convenient to encode the $k(n-k)$ boundary measurements $M_{ij}$,
$i\in I$, $j\in \bar I$, as a certain point in the Grassmannian $Gr_{kn}$.
Recall that $\Delta_J(A)$ is the maximal minor of a matrix $A$ 
in the column set~$J$.  The collection of all $\Delta_J$, for 
$k$-subsets $J\subset[n]$, form projective Pl\"ucker coordinates
on $Gr_{kn}$.

\begin{definition}
\label{def:bounary_measurement_map}
Let $\Net_{kn}$ be the set of planar directed
networks with $k$ boundary sources and $n-k$ boundary sinks.
Define the {\it boundary measurement map}
$$
\Mes: \Net_{kn}\to Gr_{kn},
$$
as follows.  For a network $N\in\Net_{kn}$ with the source set $I$
and with the boundary measurements $M_{ij}$, 
the point $\Mes(N)\in Gr_{kn}$ is given 
in terms of its Pl\"ucker coordinates $\{\Delta_J\}$ by 
the conditions that $\Delta_I\ne 0$ and 
$$
M_{ij} = \Delta_{(I\setminus\{i\})\cup \{j\}}/\Delta_I
\text{ for any } i\in I \text{ and } j\in \bar I.
$$
More explicitly, if $I=\{i_1<\cdots<i_k\}$, then the point 
$\Mes(N)\in Gr_{kn}$ is represented by the 
{\it boundary measurement matrix\/} $A(N)=(a_{ij})\in \Mat_{kn}$ such that 
\begin{enumerate}
\item
The submatrix $A(N)_I$ in the column set $I$ is the identity matrix $\Id_k$.
\item
The remaining entries of $A(N)$ are $a_{rj} = (-1)^s M_{i_r,j}$,
for $r\in [k]$ and $j\in \bar I$,
where $s$ is the number 
of elements of $I$ strictly between $i_r$ and $j$.
\end{enumerate}
\end{definition}

Note that the choice of signs of entries in $A(N)$ ensures that 
$\Delta_{(I\setminus\{i\})\cup\{j\}}(A(N)) = M_{ij}$,
for $i\in I$ and $j\in \bar I$.  Clearly, we have $\Delta_I(A(N)) = 1$.


\begin{example}
\label{ex:network_13_24}
For a network $N$ with four boundary vertices, with the source set $I=\{1,3\}$, and the sink set $\bar I=\{2,4\}$, 
we have
\begin{center}
\pspicture(-55,-35)(200,35)
\rput(-55,0){$N=$}
\pscircle[linecolor=grayish,fillstyle=solid,fillcolor=grayish,linewidth=0](0,0){30}
\psline{->}(10,10)(21.2,21.2)
\psline{<-}(10,-10)(21.2,-21.2)
\psline{<-}(-10,10)(-21.2,21.2)
\psline{->}(-10,-10)(-21.2,-21.2)
\rput(-27,27){$b_1$}
\rput(27,27){$b_2$}
\rput(27,-27){$b_3$}
\rput(-27,-27){$b_4$}

\rput(140,0){$A(N)=
\begin{pmatrix}
1 & M_{12} & 0 & -M_{14}\\
0 & M_{32} & 1 & M_{34}
\end{pmatrix}
.
$}
\endpspicture
\end{center}
In this case, we have 
$M_{12} = \frac{\Delta_{23}}{\Delta_{13}}$, 
$M_{14}=\frac{\Delta_{24}}{\Delta_{13}}$,
$M_{32}=\frac{\Delta_{12}}{\Delta_{13}}$, 
$M_{34}=\frac{\Delta_{14}}{\Delta_{13}}$.
\end{example}

The following two results establish a relationship between networks and total
positivity on the Grassmannian.  

\begin{theorem}
\label{th:Net=Gr}
The image of the boundary measurement
map $\Mes$ is exactly the totally nonnegative Grassmannian:
$$
\Mes(\Net_{kn}) = Gr_{kn}^\tnn.
$$
\end{theorem}

This theorem will follow from Corollary~\ref{cor:N_subs_Gr}
and Theorem~\ref{th:g_D}.


\begin{definition}
\label{def:subtraction_free_parametrization}
Let us say that a {\it subtraction-free rational parametrization\/}
of a cell $S_\M^\tnn\subset Gr_{kn}^\tnn$ 
is an isomorphism $f:\R_{>0}^d\to S_\M^\tnn$ such that 
\begin{enumerate}
\item The quotient of any two Pl\"ucker coordinates $\Delta_J/\Delta_I$,
$I,J\in \M$, of the point $f(x_1,\dots,x_d)\in S_\M^\tnn$
can be written as a subtraction-free rational expression in the 
usual coordinates $x_i$ on $\R_{>0}^d$.
\item For the inverse map $f^{-1}$, the $x_i$ can 
be written as subtraction-free rational expressions in terms
of the Pl\"ucker coordinates $\Delta_J$.
\end{enumerate}
Moreover, 
we say that such subtraction-free parametrization is {\it $I$-polynomial\/} 
for given $I\in \M$, if the quotients $\Delta_J/\Delta_I$, for $J\in \M$, 
are given by polynomials in the $x_i$ with nonnegative integer coefficients.
\end{definition}

Let $G$ be a planar directed graph with the set of edges $E(G)$.
Clearly, we can identify the set of all networks
on the given graph $G$ with the set $\R_{>0}^{E(G)}$ of positive real-valued
functions on $E(G)$.  The boundary measurement map induces the map
\begin{equation}
\label{eq:M_G_map}
\Mes_G:\R_{>0}^{E(G)}/\{\text{gauge transformations}\}
\to Gr_{kn}.
\end{equation}

\begin{theorem}
\label{th:Image_G_M_is_cell}
For a planar directed graph $G$, the image of the map $\Mes_G$ is a certain totally 
nonnegative Grassmann cell $S_\M^\tnn$.
\end{theorem}

For a cell $S_\M^\tnn$,  let $I_i\subset [n]$ 
as the lexicographically minimal base of the matroid $\M$
with respect to the linear order $i<i+1<\cdots <n < 1<\cdots <i-1$ on $[n]$.
In particular, $I_1=I(\lambda)$, whenever $S_\M^\tnn\subset \Omega_\lambda$.

\begin{theorem}
\label{th:net_to_Gr_parametrization}
For any cell $S_\M^\tnn$, one can find a graph $G$ such that the map $\Mes_G$ 
is a subtraction-free rational parametrization of this cell.
Moreover, for $i=1,\dots,n$,
there is an acyclic planar directed graph $G$ with the source set $I_i$ such that
$\Mes_G$ is an $I_i$-polynomial parametrization of the cell $S_\M^\tnn$.
\end{theorem}

This theorem will follow from Theorem~\ref{th:g_D}.


In the next section we will prove that $\Mes(\Net_{kn})\subseteq Gr_{kn}^\tnn$.
In other words, we will prove that all maximal minors of the boundary 
measurement matrix $A(N)$ are nonnegative.

\section{Loop-erased walks}
\label{sec:loop_erased}

In this section we generalize the well-known {\it Lindstr\"om lemma\/} 
to suit our purposes.  
Since our graphs may not be acyclic, we will follow Fomin's approach 
from the work on {\it loop-erased walks\/}~\cite{F}, which
extends the Lindstr\"om lemma to non-acyclic graphs.
The sources and sinks in our graphs may be interlaced with each other,
which gives an additional complication. 

\medskip

\begin{figure}[ht]
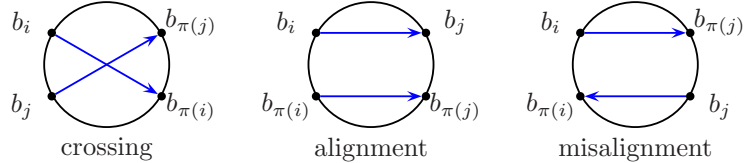

\pspicture(-50,-40)(50,40)
\rput(0,-40){crossing}
\pscircle[linecolor=black](0,0){30}
\cnode*[linewidth=0,linecolor=black](25.98,15){2}{B1}
\cnode*[linewidth=0,linecolor=black](25.98,-15){2}{B2}
\cnode*[linewidth=0,linecolor=black](-25.98,-15){2}{B4}
\cnode*[linewidth=0,linecolor=black](-25.98,15){2}{B5}
\ncline{->}{B5}{B2}
\ncline{->}{B4}{B1}
\rput(-40,20){$b_i$}
\rput(40,20){$b_{\pi(j)}$}
\rput(-40,-20){$b_j$}
\rput(40,-20){$b_{\pi(i)}$}
\endpspicture
\qquad
\pspicture(-50,-40)(50,40)
\pscircle[linecolor=black](0,0){30}
\rput(0,-40){alignment}
\cnode*[linewidth=0,linecolor=black](25.98,15){2}{B1}
\cnode*[linewidth=0,linecolor=black](25.98,-15){2}{B2}
\cnode*[linewidth=0,linecolor=black](-25.98,-15){2}{B4}
\cnode*[linewidth=0,linecolor=black](-25.98,15){2}{B5}
\ncline{->}{B5}{B1}
\ncline{->}{B4}{B2}
\rput(-40,20){$b_i$}
\rput(40,20){$b_j$}
\rput(-40,-20){$b_{\pi(i)}$}
\rput(40,-20){$b_{\pi(j)}$}
\endpspicture
\qquad
\pspicture(-50,-40)(50,40)
\pscircle[linecolor=black](0,0){30}
\rput(0,-40){misalignment}
\cnode*[linewidth=0,linecolor=black](25.98,15){2}{B1}
\cnode*[linewidth=0,linecolor=black](25.98,-15){2}{B2}
\cnode*[linewidth=0,linecolor=black](-25.98,-15){2}{B4}
\cnode*[linewidth=0,linecolor=black](-25.98,15){2}{B5}
\ncline{->}{B5}{B1}
\ncline{<-}{B4}{B2}
\rput(-40,20){$b_i$}
\rput(40,20){$b_{\pi(j)}$}
\rput(-40,-20){$b_{\pi(i)}$}
\rput(40,-20){$b_j$}
\endpspicture
\caption{Crossings, alignments, and misalignments}
\label{fig:crossing_alignment_misalignment}
\end{figure}

For two $k$-subsets $I,J\subset [n]$, let $K=I\setminus J$ and 
$L=J\setminus I$.  Then $|K| = |L|$.
For a bijection $\pi:K\to L$, we say that a pair of indices $(i,j)$,
where $i<j$ and $i,j\in K$, is a {\it crossing,} an {\it alignment,} or a 
{\it misalignment\/} of $\pi$, if the two chords 
$[b_i,b_{\pi(i)}]$ and $[b_{j},b_{\pi(j)}]$ 
are arranged with respect to each other as shown on
Figure~\ref{fig:crossing_alignment_misalignment}.
Define the {\it crossing number\/} $\xing(\pi)$ of $\pi$ 
as the number of crossings of $\pi$.

\psset{unit=0.8pt}
\psset{unit=1pt}

\begin{lemma}
\label{lem:plucker_coord_identity}
Let $I,J$ be two $k$-element subsets in $[n]$, $K=I\setminus J$
and $L=J\setminus I$.  Also let $r=|K|=|L|$.  Then the following
identity holds for the Pl\"ucker coordinates in $Gr_{kn}$
$$
\Delta_J \cdot \Delta_I^{r-1} = \sum_{\pi:K\to L}
(-1)^{\xing(\pi)} \prod_{i\in K} \Delta_{(I\setminus\{i\})\cup\{\pi(i)\}}\,,
$$
where the sum is over all $r!$ bijections $\pi:K\to L$.
\end{lemma}

Note that, in the case $r=2$, this identity is equivalent to 
a 3-term Grassmann-Pl\"ucker relation.

\begin{proof}
Let us prove this identity for maximal minors $\Delta_J(A)$ of a matrix $A$. 
We first show that the identity is invariant under permutations of
columns of the matrix.  Let us switch two adjacent rows of $A$ with
indices $a$ and $a+1$ and correspondingly modify the subsets $I$ and $J$.
Then a minor $\Delta_{M}$ switches its sign if $a,a+1\in M$;
and otherwise the minor and does not change.
Also note that the crossing number $\xing(\pi)$ changes by $\pm 1$ if
$a,a+1\in K\cup L$ and $\pi^{\pm 1}(a) \ne a+1$;
and otherwise $\xing(\pi)$ does not change.
Considering several cases according to which of the subsets
$K$, $L$, $I\cap J$, or $[n]\setminus (I\cup J)$ contain the elements $a$ 
and $a+1$, we verify in all cases that both sides of the identity 
make the same switch of sign.
Using this invariance under permutations of columns 
we can reduce the problem to the case when $I=[k]$ and $J=[k-r]\cup [k+1,k+r]$.
Let us assume that $\Delta_I\ne 0$.  Multiplying the matrix $A$
by $(A_I)^{-1}$ on the left, we reduce the problem to the case when $A_I$ is 
the identity matrix.  In this case, $\Delta_I(A) = 1$ and $\Delta_J(A)$
is the determinant of $r\times r$-matrix $B=(b_{ij})$, where
$b_{ij} = a_{i+k-r,\,j+k}$, for $i,j\in [r]$.  Note that
in this case
$\Delta_{(I\setminus\{i+k-r\})\cup\{j+k\}} = (-1)^{r-j}b_{ij}$ and 
the crossing number equals $\xing(\pi)=\binom{r}{2}-\inv(\pi)$,
where $\inv(\pi)$ is number of inversions in $\pi$.  Thus we have 
the same determinant $\det(B)$ in the right-hand side.
The case when $\Delta_I(A)=0$ follows by the continuity since we
have already proved the identity for the dense set of matrices $A$ with
$\Delta_I(A)\ne 0$.
\end{proof}

The following proposition gives an immanant expression
for the maximal minors of the boundary measurement matrix 
$A(N)$; see Definition~\ref{def:bounary_measurement_map}.

\begin{proposition}
\label{prop:Delta_pi_cn}
Let $N$ be a network with $n$ boundary vertices, including
$k$ sources $b_i$, $i\in I$, and the boundary measurements
$M_{ij}$, $i\in I$, $j\in\bar I$.
Then the maximal minors of the boundary measurement matrix $A(N)$ 
are equal to
$$
\Delta_J(A(N)) = \sum_{\pi:K\to L} (-1)^{\xing(\pi)} 
\prod_{i\in K} M_{i,\,\pi(i)}\,,
$$
for any $k$-subset $J\subset[n]$, where
the sum is over all bijections $\pi:K\to L$
from $K = I\setminus J$ to $L=J\setminus I$.
\end{proposition}

For example, for a network as in Example~\ref{ex:network_13_24},
we have $\Delta_{24}(N) = M_{12}M_{34} + M_{14}M_{32}$ because
both bijections $\pi:\{1,3\}\to \{2,4\}$ have just one misalignment
and no crossings, i.e., $\xing(\pi)=0$.

\begin{proof}
Let us express both sides of the needed identity in terms of the 
Pl\"ucker coordinates of the point $\Mes(N)\in Gr_{kn}$; 
see Definition~\ref{def:bounary_measurement_map}.
Then the identity can be reformulated as
$\Delta_J/\Delta_I = \sum_{\pi:K\to L}(-1)^{\xing(\pi)}
\prod_{i\in K} (\Delta_{(I\setminus \{i\})\cup\{\pi(i)}/\Delta_I)$,
which is equivalent to Lemma~\ref{lem:plucker_coord_identity}.
\end{proof}

For a network $N$, let us define the {\it formal boundary measurement matrix\/}
$A(N)^\form$ by in exactly the same way as the matrix $A(N)$
(see Definition~\ref{def:bounary_measurement_map}), 
but with the formal boundary measurements $M_{ij}^\form$ instead 
of the $M_{ij}$.   The elements of the matrix $A(G)^\form$ are formal 
power series in the edge variables $x_e$.  
We will show that each maximal 
minor of  $A(G)^\form$ sums up to a subtraction-free rational expression 
in the $x_e$.

Note that Proposition~\ref{prop:Delta_pi_cn} gives an expression for each
maximal minor of  $A(G)^\form$ as a certain alternating sum of products of the
$M_{ij}^\form$, which corresponds to the generating function for collections of
paths in the network $N$ that connect the boundary vertices $b_i$, $i\in K$
with the boundary vertices $b_j$, $j\in L$.  Let us show how to cancel
the negative terms in this generating function.

For the two sets $K=\{k_1<\cdots<k_r\}$ and $L$ as above, we define an 
{\it admissible collection\/} $\P=(\tilde P_1,\dots,\tilde P_r)$ of 
directed paths in $N$ by the following conditions:
\begin{enumerate}
\item The paths $\tilde P_i$ have no self-intersections.
\item There is a bijection $\pi:K\to L$ such that 
the path $\tilde P_i$ starts at the boundary vertex 
$b_{k_i}$ and ends at the boundary vertex $b_{\pi(k_i)}$, 
for $i=1,\dots,r$.
\item The bijection $\pi$ has no crossings, that is $\xing(\pi)=0$.
\item The paths $\tilde P_i$ and $\tilde P_j$ have no
intersections (i.e., common vertices) whenever 
$(k_i,k_j)$ is an alignment\footnote{However 
       the paths $\tilde P_i$ and $\tilde P_j$ are
       allowed to intersect when $(k_i,k_j)$ is a misalignment of $\pi$.}
of the bijection $\pi$;
see Figure~\ref{fig:crossing_alignment_misalignment}.
\end{enumerate}
Note that there are finitely many admissible collections of paths
because there are finitely many paths without self-intersections.
Suppose that the $i$th path an admissible collection $\P$ has vertices
$\tilde P_i = (\tilde v_{i1},\tilde v_{i2},\dots,\tilde v_{i,m_i})$.
Let $\Cyc_{ij}(\P)$ be the set of all cycles (without self-intersections)
that start and end at the $j$th vertex $\tilde v_{ij}$ of $\tilde P_i$ 
and do not pass through
any of the earlier vertices $\tilde v_{i1},\dots,\tilde v_{i,j-1}$ 
in $\tilde P_i$ and any vertex of 
a path $\tilde P_{i'}$ such that $i'<i$ and $(k_{i'}, k_i)$ 
is an alignment of $\pi$.
Again note that all sets $\Cyc_{ij}(\P)$ are finite.

\begin{proposition}
\label{prop:subtraction_free_minor}
For any $k$-element subset $J\subset[n]$
and the sets $K$ and $L$ as above, 
the maximal minor $\Delta_J(A)$ of the matrix $A:=A(N)^\form$ is given 
by the following subtraction-free rational expression:
\begin{equation}
\Delta_J(A) = \sum_{\P} x_\P \, 
\prod_{i,j} \left(1+\sum_{C\in\Cyc_{ij}(\P)} x_C\right)^{-1},
\label{eq:Delta-loop-erased}
\end{equation}
where the sum is over all admissible collections $\P$ of paths,
$x_\P$  denotes the product of the variables $x_e$ over all edges $e$ 
in all paths in $\P$, and $x_C$ is the product of the $x_e$ over 
the edges of a cycle $C$.
\label{prop:minors-exact} 
\end{proposition}

\begin{proof}
Let $P=(v_1,\dots,v_m)$ be a directed path in $G$.
If $P$ has at least one self-intersection, then find 
the first self-intersection, that is 
the minimal index $i$ such that 
$v_i=v_j$ for some $j<i$.  Let $P'=(v_1,\dots,v_j,v_{i+1},\dots,v_m)$
be the path obtained from $P$ by erasing the 
cycle
$C=(v_j,\dots,v_i)$.  If $P'$ still has a self-intersection
then again erase the first cycle in $P'$ to get another path $P''$, etc.
Finally, we obtain the path $\tilde P$ without self-intersections,
called the {\it loop-erased part\/}\footnote{The term
  {\it loop-erased\/} was borrowed from~\cite{F}, 
  where cycles are called {\it loops.}}
of $P$.

For any path $\tilde P = (\tilde v_1,\dots,\tilde v_s)$ without
self-intersections, all paths $P$ that have the loop-erased part $\widetilde P$
can be obtained from $\tilde P$ by the following inverse procedure.  Let
$\Cyc_j(\tilde P)$ be the set of cycles that start and end at the vertex
$\tilde v_j$ and avoid the vertices $\tilde v_1,\dots, \tilde v_{j-1}$.  Then
the paths $P$ are obtained from their loop-erased part $\tilde P$ by first
inserting any number of cycles from $\Cyc_1(\widetilde P)$ (in any order), then
inserting any number of cycles from $\Cyc_2(\widetilde P)$, etc. 

According to Proposition~\ref{prop:Delta_pi_cn} and the 
definition~\eqref{eq:Mij} of the formal boundary measurements 
$M_{ij}^\form$, the minor
$\Delta_J(A)$ equals to the sum
$$
\Delta_J(A) = \sum_{\pi:K\to L} (-1)^{\xing(\pi)}
\sum_{(P_1,\dots,P_r)} \prod_{i=1}^r (-1)^{\wind(P_i)} x_{P_i}
$$
where the first sum is over all bijections $\pi:K\to L$,
the second sum is over all collections of paths $(P_1,\dots,P_r)$
such that $P_i$ starts at the boundary vertex $b_{k_i}$
and end at the boundary vertex $b_{\pi(k_i)}$,
and $x_{P}:=\prod_{e\in P} x_e$.

Let $\tilde P_1,\dots,\tilde P_r$ be the loop-erased parts of 
the paths $P_1,\dots,P_r$. 
Suppose that there is a pair of indices $i<j$ such that 
such that the loop-erased part
$\tilde P_i$ intersects with $P_j$ (i.e., these paths have a common vertex)
and $(k_i,k_j)$ is either a crossing or an alignment for the bijection $\pi$.
Let us find the lexicographically minimal such pair $(i,j)$ and find the first
intersection vertex $v$ in the paths $\tilde P_i$ and $P_j$.  Let
$(P_1,\dots,P_{i-1},P_i',P_{i+1},\dots,P_{j-1},P_j',P_{j+1},\dots,P_s)$
be the family of paths, where the paths $P_i'$ and $P_j'$ are obtained
from the path $P_i$ and $P_j$ by switchings their tails at the
common vertex $v$.  This family of paths corresponds to the 
bijection $\pi':K\to L$ obtained from $\pi$ by switching $\pi(i)$ with 
$\pi(j)$.  
Note that we have $(-1)^{\xing(\pi')} = - (-1)^{\xing(\pi)}$,
because in $\pi'$ we have replaced a crossing from $\pi$ by an alignment 
(and maybe killed several other {\it pairs\/} of crossings of $\pi$),
or vise versa.
Also note that, if we apply this transformation of paths twice, we get
the original family of paths\footnote{Is is essential here that we use
the loop-erased part $\tilde P_i$ and the whole path $P_j$.  
If we find the first intersection
of $\tilde P_i$ and $\tilde P_j$ and switch the tails, the operation may not
be involutive.}  $(P_1,\dots,P_r)$.
This implies that the contributions to $\Delta_J(A)$ of all families,
for which one can find a pair $(i,j)$ as above, cancel each other.

The surviving terms in $\Delta_J(A)$ correspond to families 
of paths $(P_1,\dots,P_r)$
with bijections $\pi$ such that $\xing(\pi)=0$ and, 
for any alignment $(i<j)$ in $\pi$, the paths 
$\tilde P_i$ and $P_j$ have no common points.
That exactly means that the collection of the loop-erased parts 
$\P=(\tilde P_1,\dots,\tilde P_r)$ is an admissible
collection of paths and that all erased cycles at the $j$th vertex of
$\tilde P_i$ belong to the set $\Cyc_{ij}(\P)$.

Finally, notice that $(-1)^{\wind(P_i)}$ equals the number of 
erased cycles in $P_i$.  Thus the contribution of all terms
with a given admissible family $\P=(\tilde P_1,\dots,\tilde P_r)$ is
$$
x_\P\, \prod_{i,j} \left(1 - \sum_{C\in\Cyc_{ij}(\P)}x_C + 
\sum_{C_1,C_2\in\Cyc_{ij}(\P)} x_{C_1}x_{C_2} - \cdots\right),
$$
which is equal to
$x_\P \,\prod_{i,j} \left(1+\sum_{C\in\Cyc_{ij}(\P)} x_C\right)^{-1}$,
as needed.
\end{proof}

Proposition~\ref{prop:minors-exact} implies
 Lemma~\ref{lem:Mij}, because $M_{ij}^\form = 
\Delta_{(I\setminus\{i\})\cup\{j\}}(A(N)^\form)$.
Since any minor of $A(N)^\form$ is a subtraction-free rational
expression in the $x_e$, we also deduce that any minor of the boundary measurement
matrix $A(N)$ is a nonnegative real number.  Thus the boundary measurement
map $\Mes$ sends any network into a point in the totally nonnegative Grassmannian. 

\begin{corollary} 
\label{cor:N_subs_Gr}
We have $\Mes(\Net_{kn})\subseteq Gr_{kn}^\tnn$.
\end{corollary}

\section{\protect\Le-diagrams}
\label{sec:Le-diagrams}

In this section we define \Le-diagrams which are on one-to-one correspondence with 
totally nonnegative Grassmann cells.
\medskip

\begin{figure}[ht]
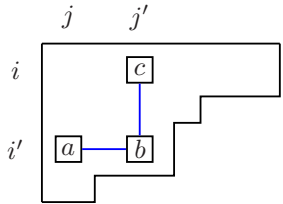

\pspicture(0,0)(90,75)
\psline[linecolor=black]{-}(0,60)(90,60)
\psline[linecolor=black]{-}(0,60)(0,0)
\psline[linecolor=black]{-}(0,0)(20,0)(20,10)(50,10)(50,30)(60,30)(60,40)(90,40)(90,60)
\rput(10,20){\rnode{A}{\rput(5.1,1.5){$a$}\psframebox[linecolor=black]{\phantom{\vrule height3pt width3pt}}}} 
\rput(37,20){\rnode{B}{\rput(5.4,1.5){$b$}\psframebox[linecolor=black]{\phantom{\vrule height3pt width3pt}}}} 
\rput(37,50){\rnode{C}{\rput(5.1,1.5){$c$}\psframebox[linecolor=black]{\phantom{\vrule height3pt width3pt}}}} 
\ncline[linecolor=blue]{-}{A}{B}
\ncline[linecolor=blue]{-}{B}{C}
\rput(-10,20){$i'$}
\rput(-10,50){$i$}
\rput(10,70){$j$}
\rput(37,70){$j'$}
\endpspicture
\caption{\protect\Le-property: if $a,c\ne 0$ then $b\ne 0$}
\label{fig:Le_property}
\end{figure}

\begin{definition}
For a partition $\lambda$, let us define a {\it \Le-diagram\/} $D$ of shape $\lambda$ as
a filling of boxes of the Young diagram of shape $\lambda$ 
with $0$'s and $1$'s such that, for any three boxes
indexed $(i',j)$, $(i',j')$, $(i,j')$, 
where $i<i'$ and $j<j'$, filled with $a$, $b$, $c$, correspondingly,
if $a,c\ne 0$ then $b\ne 0$; see Figure~\ref{fig:Le_property}.
Note that these three boxes should form a ``\Le'' 
shaped pattern.\footnote{
      The letter ``\Le'' should be pronounced as [le], because it is 
      the mirror image of ``L'' [el].
      We follow English notation for drawing Young diagrams on the plane.
      A reader who prefers another notation may opt to use one
      of the following alternative terms instead of $\Le$-diagrams:
      $\mathrm{L}$-diagram, $\Gamma$-diagram,  $\daleth$-diagram,
      $V$-diagram, $\Lambda$-diagram, $<$-diagram, or $>$-diagram.}
For a \Le-diagram $D$, let $|D|$ be the number of $1$'s it contains.
Let $\Le_{kn}$ be the set of $\Le$-diagrams whose shape $\lambda$
fits inside the $k\times (n-k)$-rectangle.
\end{definition}


\begin{figure}[ht]
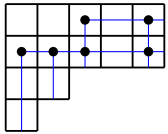

\pspicture(0,-5)(100,85)
\psline[linecolor=black](0,80)(100,80)
\psline[linecolor=black](0,60)(100,60)
\psline[linecolor=black](0,40)(100,40)
\psline[linecolor=black](0,20)(40,20)
\psline[linecolor=black](0,0)(20,0)
\psline[linecolor=black](0,80)(0,0)
\psline[linecolor=black](20,80)(20,0)
\psline[linecolor=black](40,80)(40,20)
\psline[linecolor=black](60,80)(60,40)
\psline[linecolor=black](80,80)(80,40)
\psline[linecolor=black](100,80)(100,40)

\psline[linecolor=blue, linewidth=.4pt](10,0)(10,50)(100,50)
\psline[linecolor=blue, linewidth=.4pt](30,20)(30,50)
\psline[linecolor=blue, linewidth=.4pt](50,40)(50,70)(100,70)
\psline[linecolor=blue, linewidth=.4pt](90,40)(90,70)
\pscircle*[linecolor=black](10,50){3}
\pscircle*[linecolor=black](30,50){3}
\pscircle*[linecolor=black](50,50){3}
\pscircle*[linecolor=black](50,70){3}
\pscircle*[linecolor=black](90,50){3}
\pscircle*[linecolor=black](90,70){3}
\endpspicture
\caption{A \protect\Le-diagram $D$ of shape $\lambda = (5,5,2,1)$ with $|D|=6$}
\label{fig:Le_diagram}
\end{figure}

Figure~\ref{fig:Le_diagram} shows an example of \Le-diagram.  Here dots in
boxes of the Young diagram indicate that they are filled with $1$'s, and empty
boxes are assumed to be filled with $0$'s.  Let us draw the {\it hook\/} for
each box with a dot, i.e., two lines going to the right and down
from the dotted box.  The \Le-property means that
every box of the Young diagram located at an intersection of two lines 
should contain a dot.

\psset{unit=0.6pt}
\psset{unit=1pt}

For a Young diagram filled with $0$'s and $1$'s, let us say that that a $0$ is
{\it blocked\/} if there is a $1$ somewhere above it in the same column.  
For example, the diagram shown on Figure~\ref{fig:Le_diagram} has three blocked 
$0$'s: two in the first column and one in the second column.
The \Le-property can be reformulated in terms of blocked $0$'s, as follows. 
For each blocked $0$, all entries to the left and in the same row as this $0$
are also $0$'s.

\begin{remark}
\label{rem:blocked_zeros}
One can use this observation to recursively
construct \Le-diagrams.  Suppose that we have a \Le-diagram $D$ 
of shape $\lambda$ whose last column contains $d$ boxes 
and $b$ blocked $0$'s.  Let $\tilde D$ be
the \Le-diagram of shape $\tilde \lambda$ obtained from $D$ 
by removing the last column and the $b$ rows
(filled with all $0$'s) that contain these blocked $0$'s in the last column.
The shape $\tilde \lambda$ of this diagram  is obtained from $\lambda$
by removing $b$ rows of maximal length and removing the last column.
Then $\tilde D$ can be any \Le-diagram of shape $\tilde \lambda$.
Thus an arbitrary \Le-diagram $D$ as above with prescribed $0$'s and $1$'s 
in the last column is constructed by picking an arbitrary \Le-diagram 
$\tilde D$ as above and inserting rows filled with all $0$'s in the positions 
corresponding to the blocked zeros and then inserting the last column.
\end{remark}


\begin{definition}
A {\it $\Gamma$-graph\/} is a planar directed graph $G$
satisfying the conditions:
\begin{enumerate}
\item The graph $G$ is drawn inside a 
closed boundary curve in $\R^2$.
\item $G$ contains only vertical edges oriented downward
and horizontal edges oriented to the left.
\item For any internal vertex $v$, the graph $G$
contains the line going down from $v$ until it hits the boundary 
(at some boundary sink)
and the line going to the right from $v$ until it hits the boundary
(at some boundary source).
\item All pairwise intersections of such lines should also 
be vertices of $G$.
\item The graph may also contain some number of isolated boundary 
vertices, which are assigned to be sinks or sources.
\end{enumerate}
In other words, a $\Gamma$-graph $G$ is obtained by drawing several
$\Gamma$-shaped hooks inside the boundary curve.
A {\it $\Gamma$-network\/} is a network with a $\Gamma$-graph.
\end{definition}

Note that for an arbitrary $\Gamma$-network $N$ there is a unique 
gauge transformation of edge weights~\eqref{eq:gauge_transformation} 
that transforms the weights of all vertical edges into $1$.

\begin{figure}[ht]
\includegraphics[height=1.2in]{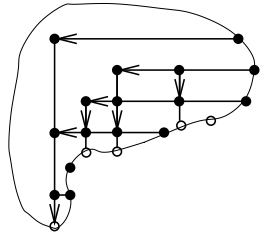}
\caption{A $\Gamma$-graph}
\label{fig:Gamma_network}
\end{figure}

Figure~\ref{fig:Gamma_network} shows an example of $\Gamma$-graph.
We displayed boundary sources by black vertices and boundary sinks
by white vertices.

For a \Le-diagram $D$ of shape $\lambda$, define {\it \Le-tableaux\/} $T$ as
nonnegative real-valued functions $T$ on boxes $(i,j)$ of the Young
diagram of shape $\lambda$ such that and $T(i,j)>0$ if and only if the box
$(i,j)$ of the diagram $D$ is filled with a $1$.  

There is a simple correspondence between \Le-tableaux of shape $\lambda$ that
fit inside the rectangle $(n-k)^k$ and $\Gamma$-networks with $k$ boundary
sources and $n-k$ boundary sinks modulo gauge transformations.
Let $T$ be a \Le-tableau 
of shape $\lambda\subseteq(n-k)^k$.
The boundary of the Young diagram of $\lambda$ gives the lattice path of length $n$ 
from the upper right corner to the lower left 
corner of the rectangle $(n-k)^k$.
Let us place a vertex in the middle of 
each step in the lattice path and mark these vertices by $b_1,\dots,b_n$
as we go downwards and to the left.
The vertices $b_i$, $i\in I$, corresponding to the vertical steps
in the lattice path will be the sources of the network and the remaining 
vertices $b_j$, $j\in \bar I$, corresponding to horizontal steps 
will be the sinks.   Notice that the source set $I$ is exactly the set 
$I(\lambda)$ as defined in Section~\ref{ssec:schubert_cells}.
Then connect the upper right corner with the lower left corner 
of the rectangle by another path so that together with the lattice
path they form a closed curve containing the Young diagram in its 
interior.   For each box $(i,j)$ of the Young diagram such that 
$T(i,j)\ne 0$, draw an internal vertex in the middle of this box and draw
the line that goes downwards from this vertex until it hits a 
boundary sink and another line that goes to the right from this vertex
until it hits a boundary source.
As we have already mentioned in Section~\ref{sec:Le-diagrams},
the \Le-property means that any intersection of such lines
should also be a vertex; cf.~Figure~\ref{fig:Le_diagram}.
Orient all edges of the obtained graph to the left and downwards.
Finally, for each internal vertex $v$ drawn in the middle of
the box $(i,j)$ assign the weight $x_e = T(i,j)>0$ to the 
horizontal edge $e$ that enters $v$ (from the right).
Also assign weights $x_e=1$ to all vertical edges of the network.
Let us denote the obtained network $N_T$.
It is not hard to see that any $\Gamma$-network (with the weights of all 
vertical edges equal to 1) comes from a \Le-tableau in this fashion.
We leave it as an exercise for the reader to rigorously prove this claim.

\begin{figure}[ht]
\includegraphics[height=1.2in]{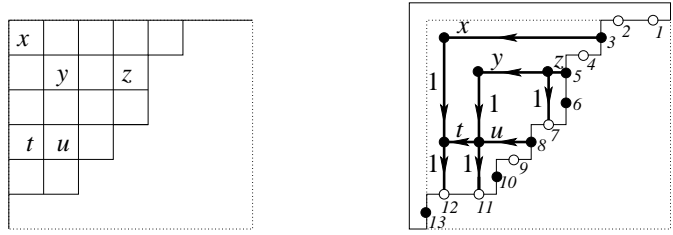}
\caption{A \protect\Le-tableau $T$ and the corresponding $\Gamma$-network
$N_T$}
\label{fig:le_tableau_gamma}
\end{figure}

\begin{example}
Figure~\ref{fig:le_tableau_gamma} gives an example of a \Le-tableau and the
corresponding $\Gamma$-network.  In the \Le-tableau only nonzero entries are
displayed.   In the $\Gamma$-network we marked the boundary vertices
$b_1,\dots,b_n$ just by the numbers $1,\dots,n$.   
The dotted lines indicate the boundary of the rectangle $(n-k)^k$; 
they are not edges of the $\Gamma$-network.  
In this example, $n=13$ and the source set is $I=\{3,5,6,8,10,13\}$.
\end{example}

For a \Le-diagram $D\in \Le_{kn}$, let $\R^D_{>0} \simeq \R^{|D|}_{>0}$ be the set of 
\Le-tableaux $T$ associated with $D$.  The map
$T\mapsto N_T$ gives the isomorphism 
$$
\R^D_{>0} \simeq \R^{E(G)}_{>0}/\{\text{gauge transformations}\}
$$
between the set of \Le-tableaux $T$ with fixed \Le-diagram $D$
and the set of $\Gamma$-networks (modulo gauge transformations)
with the fixed graph $G$ corresponding to the \Le-diagram $D$ as above.
The boundary measurement map $\Mes$ 
(see Definition~\ref{def:bounary_measurement_map})
induces the map
$$
\Mes_D:\R^D_{>0}\to Gr_{kn},
\quad \Mes_D:T\mapsto \Mes(N_T).
$$


Recall Definition~\ref{def:subtraction_free_parametrization} of a
subtraction-free parametrization.

\begin{theorem}
\label{th:g_D}
For each \Le-diagram $D\in \Le_{kn}$, the map $\Mes_D$ is a subtraction-free
parametrization a certain certain totally nonnegative Grassmann cell $S_\M^\tnn
= \Mes_D(\R_{>0}^\tnn) \subset Gr_{kn}^\tnn$.  This gives a bijection between
\Le-diagrams $D\in \Le_{kn}$ and all cells $S_\M^\tnn$ in $Gr_{kn}^\tnn$.  The
\Le-diagram $D$ has shape $\lambda$ if and only if $S_\M^\tnn\subset
\Omega_\lambda$.  The dimension of $S_\M^\tnn$ equals to $|D|$.  Moreover, for
$D$ of shape $\lambda$, the map $\Mes_D$ is $I$-polynomial, where
$I=I(\lambda)$.  
\end{theorem}

Theorem~\ref{th:g_D}, together with Corollary~\ref{cor:N_subs_Gr}, implies
Theorem~\ref{th:Net=Gr}.  Moreover, it implies Theorem~\ref{th:net_to_Gr_parametrization}.
Indeed, the parametrization $\Mes_D$ is $I_1$-polynomial,
where $I_1=I(\lambda)$.
The claim for other bases $I_i$ follows from the cyclic symmetry; 
see Remark~\ref{rem:cyclic_symmetry}.
In other words, take the \Le-diagram corresponding to a cell $S_\M^\tnn$ as above, but assuming 
that the boundary vertices are ordered as $b_i<\cdots <b_n<b_1<\cdots <b_{i-1}$.  It gives
an $I_i$-polynomial parametrization of $S_\M^\tnn$.

\section{Inverting the boundary measurement map}

In this section we prove Theorem~\ref{th:g_D} by 
constructing the bijective map 
$Gr_{kn}^\tnn  \to \{\text{\Le-tableaux}\}$, 
which is inverse to the boundary measurement map.  

\medskip

%


The construction is based on the following four lemmas.
Let $A=(a_{ij})$ be a $k\times n$-matrix in $I$-echelon form;
see Section~\ref{ssec:schubert_cells}.
Let $A_{d+1}$ be the first column-vector of $A$ that can be expressed
as a linear combination of the previous column-vectors: 
$A_{d+1} = (-1)^{d-1} x_1\,A_1+\cdots+x_{d-2}\,A_{d-2}-x_{d-1}\,A_{d-1} 
+ x_{d}\,A_{d}$.  
In other words,
$d$ is the maximal integer such that $[d]\subseteq I$.
(Here we exclude the trivial case $k=n$ when $A$ should 
be the identity matrix.
But we allow $d=0$ when the first column of $A$ is zero.)
Then $A_i$ is the $i$th coordinate vector, for $i=1,\dots,d$,
and $A_{d+1}= ((-1)^{d-1} x_1,...,x_{d-2},-x_{d-1},x_{d},0,\dots,0)^T$, that is 
the matrix $A$ has the following form:
\begin{equation}
\label{eq:A_xxx}
A=\begin{pmatrix}
1 &  \cdots & 0 & 0 & (-1)^{d-1} x_1 & * & \cdots & * \\
\vdots  & \ddots & \vdots & \vdots & \vdots & \vdots &  & \vdots \\
0 &  \cdots & 1& 0 & -x_{d-1} & * & \cdots & * \\
0 &  \cdots & 0 & 1 & x_d & * & \cdots & * \\
0 & \cdots & 0 & 0 & 0 & * & \cdots & * \\
\vdots & & \vdots & \vdots & \vdots & \vdots & & \vdots \\
0 & \cdots & 0 & 0 & 0 & * & \cdots & * \
\end{pmatrix}
.
\end{equation}

Recall that $\Mat^\tnn_{kn}$ is the set of $k\times n$-matrices
of rank $k$ with nonnegative maximal minors $\Delta_J\geq 0$.

\begin{lemma}  
We have $\Delta_{(I\setminus\{i\})\cup\{d+1\}}(A) = x_i$. 
Thus the condition $A\in \Mat^\tnn_{kn}$ implies that
$x_i\geq 0$ for $i=1,\dots,d$.
\label{lem:Le-1}
\end{lemma}

\begin{proof}  The matrix $A_{(I\setminus\{i\})\cup\{d+1\}}$
is obtained from the identity matrix $A_{I}$ by skipping 
its $i$th column and inserting the column 
$v_{d+1}=(\pm x_1,\dots,-x_{d-1},x_d,0,\dots,0)^T$ in $d$th position.
\end{proof}

\begin{lemma}
Let $r\in[d]$ be an index such that $x_r=0$ and there exists 
$i<r$ such that $x_i\ne 0$.
Then the condition $A\in\Mat^\tnn_{kn}$ implies
that $a_{rj}=0$ for all $j>d$.
In other words, the $r$-row of $A$ has only one nonzero entry
$a_{rr}=1$.
\label{lem:Le-2}
\end{lemma}

We will call an entry $x_r=0$ of the vector $(x_1,\dots,x_d)$
satisfying the condition in this lemma a {\it blocked zero.}

\begin{proof}
For $j\in I$ we have $a_{rj} = 0$ because the $j$-th column of $A$
has only one nonzero entry $a_{sj}=1$ for some $s>d$.
Suppose that $j\not\in I$.  Then 
$\Delta_{(I\setminus\{r\})\cup \{j\}} = (-1)^t \,b_{rj} \geq 0$,
where $t: = |I\cap [r+1,j-1]|$. 
On the other hand, $\Delta_{(I\setminus\{i,r\})\cup\{d+1,j\}}= 
(-1)^{t+1}\, x_i \,a_{rj} \geq 0$.
Since by Lemma~\ref{lem:Le-1} $x_i > 0$,
we deduce that $(-1)^{t+1}\,a_{rj}\geq 0$.
This implies that $a_{rj} =0$.
\end{proof}

\begin{lemma}
\label{lem:Le-3}
Assume that the $r$-row of $A$ has only one nonzero entry
$a_{rr}=1$ for some $r\in[d]$.
Let $B=(b_{ij})$ be the $(k-1)\times (n-1)$-matrix obtained from $A$ 
by removing the $r$-row and the $r$-column and inverting signs 
of the entries $a_{ij}$ for $i=1,\dots,r-1$ and $j\geq d+1$.
Then $A\in \Mat^\tnn_{kn}$ if and only if $B\in\Mat^\tnn_{k-1,n-1}$.

Moreover, the maximal minors of the matrices $A$ and $B$ are equal to 
each other.  More explicitly, $\Delta_J(A) = 0$ if $r\not\in J$, and
$\Delta_J(A) = \Delta_{\tilde J\setminus\{r\}}(B)$ if $r\in J$,
where $\tilde J$ means that we decrease  elements  $>r$ in 
$J$ by $1$.  
\end{lemma}

\begin{proof}
The equality of the minors is straightforward;
it implies the first claim.
\end{proof}

\begin{lemma}
\label{lem:Le-4}
Assume that there are no blocked zeros, that is 
$x_1=\cdots=x_s=0$ and $x_{s+1},x_{s+2},\dots,x_d > 0$,  
for some $s\in[0,d]$.
Let $C=(c_{ij})$ be the $k\times(n-1)$-matrix whose first $d$ columns
are the first coordinate vectors (as in the matrix $A$) 
and the remaining entries are
$$
c_{i,\,j-1}=\left\{
\begin{array}{cl}
a_{ij}&\text{if } i\in[s]\cup[d+1,k],\\[.05in]
\frac{a_{ij}}{x_i} + 
\frac{a_{i+1,\,j}}{x_{i+1}} & 
      \text{if } i\in [s+1,d-1],\\[.05in]
\frac{a_{dj}}{x_d} & \text{if } i=d,
\end{array}
\right.
$$
for $j=d+2,\dots,n$.
Then $A\in \Mat^\tnn_{kn}$ if and only if $C\in\Mat^\tnn_{k,n-1}$.

Moreover, if we fix a totally nonnegative cell 
$S_\M^\tnn \subset Gr_{kn}^\tnn$ and require that (the coset of) 
the matrix $A$ belongs to $S_\M^\tnn$, then 
we can write all maximal minors $\Delta_J(C)$ as subtraction-free rational 
expressions in terms of the minors $\Delta_K(A)$, $K\in\M$.
On the other hand, 
we can write the minors $\Delta_K(A)$ as nonnegative integer
polynomials in terms of the minors $\Delta_J(C)$
and the $x_{i}$, $i\in[s+1,d]$.

This gives a bijective correspondence between
totally nonnegative cells $S_\M^\tnn\subseteq Gr_{kn}^\tnn$
that can contain a matrix $A$ of this form
and totally nonnegative cells $S_{\M'}^\tnn\subseteq Gr_{k,n-1}^\tnn$
that can contain a matrix $C$ of this form.
\end{lemma}

We will prove this lemma in Section~\ref{sec:tails}.
Note that the matrix $C$ is still in echelon form.
Let us illustrate this lemma by an example.

\begin{example}  Let 
$$
A=\begin{pmatrix}
1& 0 & - x_1 & y \\
0& 1 & x_2 & z 
\end{pmatrix}
\quad\text{and}\quad
C=\begin{pmatrix}
1& 0 &  \frac{y}{x_1} + \frac{z}{x_2} \\
0& 1 &   \frac{z}{x_2}
\end{pmatrix},
$$
where $x_1>0$ and $x_2>0$.
Then $\Delta_{12}(C)  = 1$,
$\Delta_{13}(C) = \frac{\Delta_{14}(A)}{\Delta_{13}(A)}$,
$\Delta_{23}(C) = \frac{\Delta_{34}(A)}{\Delta_{13}(A)\cdot\Delta_{23}(A)}$.
On the other hand,
$\Delta_{12}(A)=1$, $\Delta_{13}(A) = x_2$, $\Delta_{23}(A)=x_1$,
$\Delta_{14}(A)= x_2\,\Delta_{13}(C)$, $\Delta_{24}(A) = 
x_1\, (\Delta_{13}(C) + \Delta_{23}(C))$,  
$\Delta_{34}(A) = x_1\,x_2\, \Delta_{23}(C)$.
\end{example}


\begin{proof}[Proof of Theorem~\ref{th:g_D}]
Let us prove the theorem, together with the additional claim that the map
$\R_{>0}^D\to\Mat_{kn}$ given by $T\mapsto A(N_{T})$ produces
$k\times n$-matrices in echelon form.  The proof is by induction on $n$.
The cases when $k=n$ or $k=0$ are trivial, which provides the base of 
induction.  Assume that $k\in[n]$.  Assume by induction that 
the theorem is valid for all $Gr_{k'n'}^\tnn$ with $n'<n$.



Let $A$ be a $k\times n$-matrix in $I$-echelon form that represents a point
in $S_\M^\tnn\subseteq \Omega_\lambda\cap Gr_{kn}^\tnn$, where
$I=I(\lambda)$.  Let us find the integer $d$ and the real numbers
$x_1,\dots,x_d$ as above in this section; see~\eqref{eq:A_xxx}.  According to
Lemma~\ref{lem:Le-1}, we have $x_i\geq 0$, for $i=1,\dots,d$.
Note that we can uniquely determine the number $d$ and the set of indices $i$
with $x_i\ne 0$ from the matroid $\M$, because these numbers
are certain maximal minors of $A$; see Lemma~\ref{lem:Le-1}. 



Suppose that there are $b>0$ blocked zeros in the 
vector $(x_1,\dots,x_d)$.  For each blocked zero $x_r=0$,
all entries in the $r$th row of $A$ are zero, except $a_{rr}=1$.
Let $A'$ be the $(k-b)\times (n-b)$-matrix obtained 
from $A$ by skipping the $r$th row and the $r$th column, for each blocked zero $x_r=0$,
and inverting signs of some entries, as in Lemma~\ref{lem:Le-3}. 
Namely, we need to invert the sign of $a_{ij}$, $i\leq d<j$, if and only if there is
an odd number of blocked zeros $x_r=0$ with $r>i$.
According Lemma~\ref{lem:Le-3}, the maximal minors of $A$ are equal 
to the corresponding maximal minors of $A'$ (or to zero).
Thus the matroid $\M'$ associated with 
$A'$ can be uniquely constructed from the matroid $\M$;
and vise versa, if we know the matroid $\M'$ and the positions of blocked
zeros, then we can uniquely reconstruct the matroid $\M$.
In particular, the Schubert cell $\Omega_{\lambda'}$ of $A'$ corresponds to the
Young diagram of shape $\lambda'$ obtained by removing $b$ (longest) rows with
$n-k$ boxes from the Young diagram of $\lambda$.

Note that $A'$ can be any $I'$-echelon matrix, where $I'=I(\lambda')$,
that represents a cell $S_{\M'}^\tnn \subset \Omega_{\lambda'}$ such that the
vector $(x_1',\dots,x_{d-b}')$ with  $x_i' = \Delta_{(I'\setminus \{i\})\cup
\{d-b+1\}}(A')$ has no blocked zeros.  By the induction hypothesis, we already
know that the map $T'\mapsto A(N_{T'})$ is a bijection between \Le-tableaux
$T'$ of shape $\lambda'$ such that the last column of $T'$ contains no blocked
zeros and the set of matrices $A'$ as above.  Moreover, it gives a bijection
between \Le-diagrams $D'$ corresponding to such tableaux and cells
$S_{M'}^\tnn$ containing such matrices $A'$; and the map $T'\mapsto A(N_{T'})$
gives a subtraction-free parametrization for each cell $S_{\M'}^\tnn$.

Let $T'$ be the \Le-tableau such that $A(N_{T'}) = A'$ and $D'$ be its
\Le-diagram.  Let $T$ be the \Le-tableau (and $D$ be its \Le-diagram) obtained
from $T'$ (reps., from $D$) by inserting $b$ rows filled with all $0$'s in the
positions corresponding to blocked zeros in $(x_1,\dots,x_d)$.  Then we have
$A(N_T)=A$.  Indeed, the network $N_T$ is obtained from $N_{T'}$ by inserting
$b$ isolated sources in the positions corresponding to the blocked zeros.
Thus, according to Definition~\ref{def:bounary_measurement_map}, its boundary
measurement matrix $A(N_T)$ is obtained from $A'=A(N_{T'})$ by inserting, for
each blocked zero $x_r$, a row and a column in the $r$th positions with a
single nonzero entry $a_{rr}=1$, and switching the signs of some entries
$a_{ij}$.  (The signs are switched because we insert additional elements $r$
into $I$.  These switches are exactly the same as in Lemma~\ref{lem:Le-3}.)
Note that all \Le-diagrams and all \Le-tableaux of shape $\lambda$ with 
given entries in the last column are of this form;
see Remark~\ref{rem:blocked_zeros}.
This implies that $T\mapsto A(N_T)$ gives subtraction-free 
parametrizations
with needed properties for the cells $S_\M^\tnn$ with some blocked zeros.





Let us now suppose that $(x_1,\dots,x_d)$ contains no blocked zeros,
that is $x_1=\cdots=x_s=0$ and $x_{s+1},\dots,x_d> 0$.  
Let $C$ be the matrix obtained from $A$ as in Lemma~\ref{lem:Le-4}.
The matrix $C$ represents a point in $S_{\M''}^\tnn\subset
\Omega_{\lambda''}\cap Gr_{k,n-1}^\tnn$.
According to Lemma~\ref{lem:Le-4}, the matroid $\M''$ is uniquely
determined by the matroid $\M$, and vise versa $\M$ is uniquely
determined by $\M''$ and positions of zeros in $(x_1,\dots,x_d)$.
In this case, the Young diagram of $\lambda''$ is obtained from the 
Young diagram of $\lambda$ by removing the last column (with $d$ boxes).

Note that $S_{\M''}^\tnn$ can be any cell in $\Omega_{\lambda''}$ and $A''$ can
be the echelon representative of any point in such cell.  Again, by the
induction hypothesis, we already know that $T''\mapsto A(N_{T''})$ gives
a subtraction-free parametrization for cells $S_{\M''}^\tnn\subset 
\Omega_{\lambda''}$ 
and that this map induces a bijection between all \Le-diagrams $D''$ of shape $\lambda''$ and 
these cells. 
Let $T''$ be the \Le-tableau such that $A(N_{T''}) = C$, and $D''$ be its
\Le-diagram.  Let $T$ be the \Le-tableau obtained from $T''$ by inserting the last 
column filled with $(x_1,\dots,x_d)$, and $D$ be the \Le-diagram of $T$,
i.e., $D$ is obtained from $D''$ by inserting the last row filled with 
$(0,\dots,0,1,\dots,1)$ ($s$ zeros followed by $d-s$ ones).
According to Remark~\ref{rem:blocked_zeros}, any \Le-diagram $D$ and any \Le-tableau $T$
of shape $\lambda$ without blocked zeros in the last column have this form.

\begin{figure}[ht]
\input{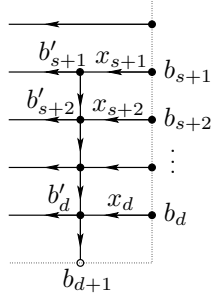}
\caption{A piece of the $\Gamma$-network $N_T$}
\label{fig:gamma_network_induction}
\end{figure}

We claim that for this \Le-tableau $T$, we have $A(N_T) = A$.
Indeed, the network $N_{T''}$ is obtained from $N_T$ by removing the $d-s$ internal
vertices $b_{s+1}',\dots,b_{d}'$ adjacent to the boundary sources
$b_{s+1},\dots, b_d$, removing the vertical edges incident to the vertices
$b_i'$, merging $d-s$ pairs of horizontal edges, erasing the weights
$x_{s+1},\dots,x_d$, removing the sink $b_{d+1}$, and shifting the labels of
the boundary vertices $b_i$, $i>d+1$, by one; see
Figure~\ref{fig:gamma_network_induction}.
The boundary measurements $M_{ij}$ of the network $N_T$ are 
obtained from the boundary measurements of $M_{ij}''$ of $N_{T''}$
by $M_{i,j+1} = x_i (M_{ij}'' + M_{i+1,j}'' + \cdots
+ M_{d,j}'')$ for $i\in[s+1,d]$, $j\geq d+1$;
$M_{i,j+1} =  M_{ij}''$ for $i\not\in [s+1,d]$, $j\geq d+1$;
and $M_{i, d+1} = x_i$ for $i\in [k]$;
see Figure~\ref{fig:gamma_network_induction}.
Thus the boundary measurement matrix $\tilde A=(\tilde a_{ij})$ of $N_T$ is obtained
from the matrix $C=A(N_{T''})$, by inserting the column
$(x_1,\dots,x_d,0,\dots,0)^T$ in the $(d+1)$st position
and changing other entries as
$$
\tilde a_{i,j+1} = 
\left\{
\begin{array}{cl}
x_i (c_{ij} - c_{i+1,j}+ \cdots \pm c_{dj})
& \text{if } i\in [s+1,d],\\
c_{ij} & \text{otherwise,}
\end{array}
\right.
$$
for $j\geq d+1$; and of course $\tilde a_{ij}=c_{ij} = \delta_{ij}$, 
for $j\leq d$.
Notice that the $C\mapsto \tilde A$ is the inverse of the transformation
$A\mapsto C$ from Lemma~\ref{lem:Le-4}.
Thus $A(N_T)$ is the original matrix $A$, as needed. 
Now Lemma~\ref{lem:Le-4} implies that 
the map $T\mapsto A(N_T)$ gives needed subtraction-free parametrizations
for the cell $S_\M^\tnn$ without blocked zeros.

Finally, note that the network $N_T$ is acyclic and its source set is
$I=I(\lambda)$.  Thus all boundary measurements $M_{ij}$ are nonnegative
polynomials in the edge weights.   By
Proposition~\ref{prop:subtraction_free_minor} all maximal minors
$\Delta_J(A(N_T))$ are given nonnegative polynomials and $\Delta_I(A(N_T)) =
1$.  This implies the last claim of Theorem~\ref{th:g_D} about
$I$-polynomiality of the map $\Mes_D$.  
\end{proof}




We can describe the map $Gr_{kn}^\tnn\to\{\textrm{$\Le$-tableaux}\}$
implicitly constructed in the above proof, 
which is inverse to the boundary measurement map,
 via the following recursive procedure.
This procedure transforms points in the totally nonnegative part of
a Schubert cell $\Omega_\lambda\cap Gr_{kn}^\tnn$
into \Le-tableaux $T$ of shape $\lambda$.  
It inserts  nonnegative real numbers into boxes of the Young diagram 
of shape $\lambda$ starting with the rightmost column, then filling
the next available rightmost column, etc.
The procedure uses a variable matrix $A$ of variable dimensions.

\medskip
\noindent
{\bf Procedure.}
Map from $\Omega_\lambda \cap Gr_{kn}^\tnn$ to \Le-tableaux of shape $\lambda$.
\begin{enumerate}
\item
Take the $k\times n$-matrix $A$ in echelon form 
representing a point in $\Omega_\lambda \cap 
Gr_{kn}^\tnn$.
\item For the matrix $A$, find the integer $d$ and 
real numbers $x_1,\dots,x_d$ as in~\eqref{eq:A_xxx}.
Insert the nonnegative real numbers $x_1,\dots,x_d$
(see Lemma~\ref{lem:Le-1}) into the empty boxes of the rightmost 
available column of $T$, skipping the boxes of $T$ 
which are already filled with $0$'s.
\item
Let $B = \{r \mid x_r = 0\textrm{ and } x_i \ne 0
\textrm{ for some } i<r\}$ be the set of blocked indices.
\item
If $B\ne \emptyset$,
then for each index $r\in B$, 
invert the sign of entries $a_{ij}$, $i<r$, $j\geq d+1$, in the matrix $A$,
remove the $r$th row and the $r$th column from $A$
(see Lemma~\ref{lem:Le-3}), and
insert zeros in all boxes of $T$ to the left of the blocked zero $x_r$.
\item 
If the obtained matrix has $0$ rows, then stop.  Otherwise go to step (2).
\item
If $B = \emptyset$, then
transform the matrix $A$ into the $(n-1)\times k$-matrix
as in Lemma~\ref{lem:Le-4}.
\item
If the obtained matrix has dimensions $k'\times n'$ with $k'=n'$
then stop.
Otherwise, go to step (2).
\end{enumerate}

%
%
%
%
%

\section{Lemma on tails}
\label{sec:tails}

In this section we prove Lemma~\ref{lem:Le-4}
essential in the proof of Theorem~\ref{th:g_D}.
\medskip

%

Let us identify a sequence $(v_1,\dots,v_n)\in(\R^k)^n$  of $k$-vectors
with the $k\times n$-matrix with the column vectors $v_i$.
We say that such sequence is {\it totally nonnegative\/} if all maximal 
$k\times k$-minors of the corresponding matrix are nonnegative
and at least one of these minors is strictly positive.

\begin{definition}
For a sequence $u=(u_1,\dots,u_m)\in(\R^k)^m$,
we define the {\it $r$-tail\/} $\Tail_r(u)$ of $u$ as the set of sequences 
$w=(w_1,\dots,w_r)\in (\R^k)^r$ such that the concatenation 
$(u,v):=(u_1,\dots,u_m,w_1,\dots,w_r)$ of $u$ and $w$
is totally nonnegative:
$$
\Tail_r(u):=\{w\in (\R^k)^r\mid 
(u,w)\textrm{ is totally nonnegative}\}.
$$
The set $\Tail_r(u)$ comes equipped with a {\it stratification,}
that is a subdivision into the disjoint union of strata
defined as follows.  We say that $w,w'\in \Tail_r(u)\subset (\R^k)^r$ 
are in the same {\it stratum\/} of $\Tail_r(u)$ if
the the corresponding $k\times (m+r)$-matrices $(u,w)$ and $(u,w')$
are in the same matroid strata, see Section~\ref{ssec:Matroid_strata},
that is the maximal minor $\Delta_I(u,w)$ 
is nonzero if and only if the maximal minor 
$\Delta_I(u,w')$ is nonzero, for any $k$-element subset $I\subseteq[m+r]$.
\end{definition}

\begin{definition}
Let $u$ and $v$ be two sequences of vectors in $\R^k$,
which are allowed to have different lengths.
We say that the sequences $u$ and $v$ are 
{\it tail-equivalent,} and write $u\sim_\mathrm{tail} v$,
if, for any $r\geq 1$, we have
\begin{enumerate}
\item $\Tail_r(u)=\Tail_r(v)$;
\item the stratifications of $\Tail_r(u)$ and $\Tail_r(v)$ are the same;
\item for any stratum $S$ of $\Tail_r(u)$, 
there is a subtraction-free rational expression 
for each maximal minor $\Delta_I(v,w)$
in terms of the nonzero maximal minors $\Delta_J(u,w)$,
for $w\in S$;
and vise versa the minors $\Delta_J(u,w)$
can be written as subtraction-free rational expressions 
in terms of the $\Delta_I(v,w)$
on each stratum of $\Tail_r(v)$.
\end{enumerate}
Note that condition (3) actually implies conditions (1) and (2).
\end{definition}

The following lemma is the main technical tool in 
our proof of Theorem~\ref{th:g_D}

\begin{lemma} 
Let $u_1,\dots,u_d$ be a linearly independent set of 
vectors in $\R^k$.
Then the following two sequences are tail-equivalent:
$$
(u_1,\dots,u_d)\sim_\mathrm{tail}
(v_1,\dots,v_{d+1}),
$$
where $v_1=u_1,\,v_2=u_1+u_2,\,v_3=u_2+u_3,\dots,
v_d=u_{d-1}+u_d,\,v_{d+1}=u_d$.
\label{le:tail-equivalent}
\end{lemma}

\begin{proof}
For three subsets $I\subseteq[d]$,
$J\subseteq[d+1]$, and $K\subseteq[r]$ with the
total cardinality $|I|+|J|+|K|=k$,
let $\Delta_{I,J,K}=\Delta_{I,J,K}(u,v,w)$ denote 
the determinant  of the $k\times k$-matrix
with columns $\{u_i\}_{i\in I}$, $\{v_j\}_{j\in J}$, $\{w_l\}_{l\in K}$, 
where the vectors $u_i$, $v_j$, $w_l$ are taken in the 
order that agrees with the total order
$$
v_1<u_1<v_2<u_2<v_3<\dots<u_{d-1}<v_d<u_d<v_{d+1}<w_1,\dots<w_r.
$$
In other words, the $\Delta_{I,J,K}$ are the maximal $k\times k$-minors 
of the $k\times(2d+1+r)$-matrix formed by these column vectors.


We can write each minor $\Delta_{\emptyset,J,K}$ as a sum the 
minors $\Delta_{I',\emptyset,K'}$ by replacing each column vector
$v_j$ in $\Delta_{\emptyset,J,K}$ with the sum $u_{j-1}+u_j$ 
and expanding the result by the linearity.  Thus the nonnegativity 
of all minors $\Delta_{I,\emptyset,K}$ implies
that all minors $\Delta_{\emptyset,J',K'}$ are also nonnegative.
Moreover, if we know which minors 
$\Delta_{I,\emptyset,K}$ are strictly positive, then
we can determine which minors 
$\Delta_{\emptyset,J',K'}$ are strictly positive.
In other words, we have $\Tail_r(u)\subseteq\Tail_r(v)$ and
the stratification of $\Tail_r(u)$ refines 
the stratification of $\Tail_r(v)$.

The proof of the opposite claim is a little bit more elaborate.  
Let us fix a stratum of $\Tail_r(v)$.  We will show that,
for $w$ in this stratum, all minors  $\Delta_{I,J,K}$ can be written
as subtraction-free rational expressions in terms of the 
nonzero minors $\Delta_{\emptyset,J',K'}$.
In particular, this would imply that all minors $\Delta_{I,\emptyset, K}$
are nonnegative and we can determine which of these minors 
are strictly positive.  Thus we will get
$\Tail_r(v)\subseteq\Tail_r(u)$ and 
deduce that the stratification of $\Tail_r(v)$ refines 
the stratification of $\Tail_r(u)$.

Let us obtain subtraction-free rational expressions for 
the $\Delta_{I,J,K}$ in terms of the $\Delta_{\emptyset, J',K'}$
by induction.  Clearly, we get such expressions when $I=\emptyset$.
Let us assume that $I$ is nonempty.
If $|J|=d$ or $d+1$, then $\Delta_{I,J,K}=0$,
because any vector $u_i$ and any $d$ vectors $v_j$ are linearly dependent.
This provides the base of the induction.  
Let us assume by induction we have already determined which of the 
minors $\Delta_{I',J',K'}$ are nonzero and wrote them
as subtraction-free rational expressions
for all triples of subsets $I'$, $J'$, $K'$ such that
$|J'|>|J|$.

Let us show that either $\Delta_{I,J,K} = 0$ or 
$\Delta_{I,J,K}=\Delta_{I',J',K'}$ for some triple $(I',J',K')$ 
with $|J'|> |J|$, 
or else there is a triple of indexes $(a,b,c)$ such that 
$1\leq a\leq b<c\leq d+1$;
$b\in I$; $a,c\not\in J$; and the collection of $|I|+|J|+2$ vectors 
$v_a,v_c$ and $u_i, v_j$, for $i\in I$, $j\in J$,
is linearly independent.

Let us give a criterion when a subsequence $S$ of vectors in 
$u_1,\dots,u_d,v_1,\dots,v_{d+1}$  is linearly independent.  
We can graphically present such collection $S$ of vectors by the graph $G_S$ 
on the set $\{0,\dots,d+1\}$ with some marked vertices where $G_S$ 
has the edge $(j-1,j)$ for each $v_j\in S$, $G_S$ has the marked 
vertex $i$ for each $u_i\in S$, and the vertices $0$ and $d+1$ 
are always marked.
Then $S$ is a linearly independent collection of vectors if and only if
each connected component of the graph $G_S$ has at most one marked vertex.


The following figure shows the graph $G_S$ for 
$d=8$ and $S = (u_5, v_1,v_4,v_5,v_6,v_8)$.  The vertex labels are
displayed below the graph and the edge labels are displayed above the graph. 
The vector $u_5$ corresponds to the marked vertex $5$
and the vectors $v_1,v_4,v_5,v_6,v_8$ correspond to edges of $G_S$
labelled $1,4,5,7,8$.  The marked vertex $5$ belongs to the connected
component of $G_S$ with the edges labelled $4,5,6$.  For this 
connected component, we have $(a,b,c) = (3,5,7)$.

\begin{center}
\includegraphics[height=0.4in]{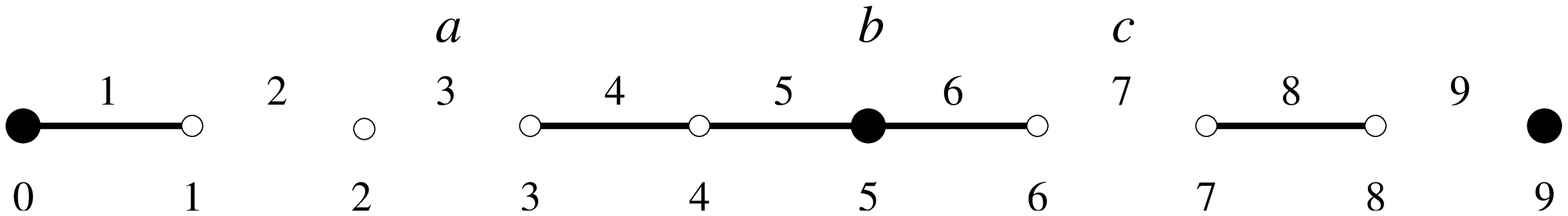}
\end{center}

Now let $S$ be the collection of vectors $u_i, v_j$, for $i\in I$ and 
$j\in J$.  If $S$ is linearly dependent, then $\Delta_{I,J,K} = 0$.  
So we can assume that $S$ is linearly independent.
Let us pick any element $b\in I$ and find the connected component
of $G_S$ that contains the marked vertex $b$.  Let $a,c$ be as above,
i.e., $a,c$ are the integers such that $a\leq b< c$, $[a+1,c-1]\subset J$,
and $a,c\not\in J$.

If $a=0$, that is $[1,b]\subseteq J$, then $\Delta_{I,J,K}=0$
because the vectors $v_1,\dots,v_b,u_b$ are linearly dependent.  
(In this case the connected component of $G_S$ contains two marked points 
$0$ and $b$.)  If $a=1$, that is $1\not\in J$ and $[2,b]\subset J$, then 
the vectors $v_2,\dots,v_b,u_b$ in $S$ can be expressed
in term of the vectors $v_1=u_1,v_2,\dots,v_b$ by a unimodular
linear transformation.   Thus we can replace the columns
$v_2,\dots,v_b,u_b$ in the minor $\Delta_{I,J,K}$ with the columns
$v_1,\dots,v_b$ without affecting the value of the minor.
So we have $\Delta_{I,J,K}=\Delta_{I\setminus\{b\},\,J\cup\{1\},\,K}$,
where we have already expressed the minor in the right-hand side by 
the induction hypothesis.
Thus we can assume that $a\geq 2$.
Analogously, we eliminate the cases $c=d+2$ and $c=d+1$ and assume that 
$c\leq d$.

Suppose that the vertex $a-1$ belongs to a connected component of $G_S$ with 
a marked vertex $a'$, that is we have $a'\in I$ and $[a'+1,a-1]\in J$ for some
$a'\leq a-1$ (or just $[1,a-1]\in J$ when $a'=0$). 
In this case we can express the sequences of vectors 
$u_{a'},v_{a'+1},\dots,v_{a-1}$ and
$v_{a+1},\dots,v_b,u_b$ (that belong to $S$) in terms of the consecutive 
sequences of vectors $u_{a'},\dots,u_{a-1}$ and $u_a,\dots,u_b$ 
by unimodular transformations, then join these two sequences into
a single sequence and unimodularly transform it again into
$v_{a'+1},\dots,v_{b}, u_b$.  The graph of the new sequence is obtained
by merging two connected components of $G_S$ with marked points $a'$
and $b$ into a single connected connected component.
The new sequence of vectors is obtained
from the old by skipping $u_{a'}$ and adding $v_a$.  This transformation
will not affect the value of the minor: $\Delta_{I,J,K} = 
\Delta_{I\setminus\{a'\},J\cup\{a\},K}$,
when $a'\geq 0$.  In the case when $a'=0$, we get 
$\Delta_{I,J,K} = \Delta_{I\setminus \{b\},J\cup\{a\},L}$. 
In all cases, the minor $\Delta_{I,J,K}$ is reduced to a minor
which has already been expressed  subtraction-freely 
by the induction hypothesis.

So we can assume that the vertex $a-1$ does not belong to a connected component
of $G_S$ with a marked vertex.  Similarly, we can eliminate the case when the
vertex $c$ belongs to a connected component with a marked point.
Now the graph obtained from $G_S$ by adding the edges $(a-1,a)$
and $(c-1,c)$ still has at most one marked vertex in each connected
component.  By the above criterion this means that the set of vectors
$\{u_i,v_j\mid i\in I,\,j\in J\}\cup \{v_a,v_c\}$ is linearly 
independent.

If $\Delta_{I,J,K}\ne 0$ then the set of vectors 
$B=\{u_i,v_j, w_l\mid i\in I,\, j\in J,\, l\in K\}$ forms a basis in $\R^k$, 
and thus the set of $k+2$ vectors $B\cup\{v_a,v_c\}$ linearly spans $\R^k$.
In this case (according to the exchange axiom) we can complete the linearly 
independent subset
$\{u_i,u_j\mid i\in I,j\in J\}\cup \{v_a,v_c\}$ of the latter set to
a basis of $\R^k$.  Thus there exist a pair of elements $x,y\in K$ such that 
$(B\setminus \{w_x,w_y\})\cup \{v_a,v_c\}$ is a basis of $\R^k$,
or equivalently, $\Delta_{I,J\cup\{a,c\},K\setminus \{x,y\}}\ne 0$.
It follows that if $\Delta_{I,J\cup\{a,c\},K\setminus \{x,y\}}= 0$
for all pairs $x,y\in K$ then $\Delta_{I,J,K}=0$.
By the induction hypothesis we have already determined which 
of the minors $\Delta_{I,J\cup\{a,c\},K\setminus \{x,y\}}$ are nonzero.
Thus we either deduce that $\Delta_{I,J,K}=0$ or find a pair $x<y$
with nonzero (and thus strictly positive)
minor $\Delta_{I,J\cup\{a,c\},K\setminus \{x,y\}}$.

We can now assume that there exist five indexes $a\leq b<c$ and $x<y$ 
such that $b\in I$; $a,c\in [d]\setminus J$; $x,y\in [r]$; 
and $\Delta_{I,J\cup\{a,c\},K\setminus \{x,y\}}>0$.
For a subword $s$ in the word $abcxy$, let $[s]$ be a shorthand 
for the triple of subsets $I',J',K'$ such that
$I'$ contains $b$ if and only if $s$ contains the letter $b$,
$J'$ contains $a$ or $c$ if and only if $s$ contains the corresponding letter,
$K'$ contains $x$ or $y$ if and only if $s$ contains the corresponding letter,
and all other entries in $I', J', K'$ are the same as $I,J,K$.
For example, $\Delta_{[bxy]}=\Delta_{I,J,K}$
and $\Delta_{[abc]} = \Delta_{I,J\cup\{a,c\},K\setminus\{x,y\}}$. 

Let us write the following 3-term Grassmann-Pl\"ucker relations:
\begin{eqnarray}
\label{eq:abc-b}
&&
\Delta_{[bxy]}\cdot
\Delta_{[abc]} +
\Delta_{[aby]}\cdot
\Delta_{[bcx]} =
\Delta_{[abx]}\cdot
\Delta_{[bcy]}\,,\\[.1in]
\label{eq:abc-a}
&&
\Delta_{[axy]}\cdot
\Delta_{[abc]} +
\Delta_{[aby]}\cdot
\Delta_{[acx]} =
\Delta_{[abx]}\cdot
\Delta_{[acy]}\,,\\[.1in]
\label{eq:abc-c}
&&
\Delta_{[cxy]}\cdot
\Delta_{[abc]} +
\Delta_{[acy]}\cdot
\Delta_{[bcx]} =
\Delta_{[acx]}\cdot
\Delta_{[bcy]}\,.
\end{eqnarray}
By the induction hypothesis, for each minor in the above
equations, except $\Delta_{[bxy]}=\Delta_{I,J,K}$, 
we have already proved its nonnegativity,
determined if it is zero or strictly positive,
and found a subtraction-free rational expression in terms
of the minors $\Delta_{\emptyset, J',K'}$.
We also know that $\Delta_{[abc]}>0$.

Suppose that $\Delta_{[acx]}= 0$. 
Replacing the column $v_{a} = u_{a-1} + u_{a}$
in this minor by the sum of two columns $u_{a-1}$ and $u_a$, we get
$\Delta_{[acx]} = \Delta_{(I\setminus\{b\})\cup\{a-1\},J\cup\{c\},
K\setminus\{y\}}
+\Delta_{I,J\cup\{c\},K\setminus\{y\}} = 0$. 
The both summands in the right-hand size are nonnegative
by the induction hypothesis and thus they should be zero.
In particular, the  second summand is $\Delta_{[bcx]}=
\Delta_{I,J\cup\{c\},L\setminus \{y\}}= 0$.
Similarly, replacing the column $v_c$
in $\Delta_{[acx]}$ by the sum of $u_{c-1}$ and $u_c$,
we get $\Delta_{[acx]} = 
\Delta_{[abx]}
+ 
\Delta_{(I\setminus\{b\})\cup\{c\},J\cup\{a\},L\setminus\{y\}} = 0$
and thus $\Delta_{[abx]} = 0$.
So the second term in the left-hand side 
and the term in the right-hand size of~\eqref{eq:abc-b} 
are both zero, which implies that 
$\Delta_{[bxy]}\cdot \Delta_{[abc]}  = 0$
and thus $\Delta_{[bxy]} = 0$.
Similarly, when $\Delta_{[acy]}= 0$, we deduce that 
$\Delta_{[aby]}=\Delta_{[bcy]}=0$ and thus again
$\Delta_{[bxy]} = 0$.
We obtain that, if $\Delta_{[acx]}= 0$ or $\Delta_{[acy]}= 0$, then 
$\Delta_{I,J,K} = \Delta_{[bxy]} = 0$.

We can now assume that $\Delta_{[acx]}>0$ and $\Delta_{[acy]}>0$.
Let us multiply both sides of equations~\eqref{eq:abc-a} 
and~\eqref{eq:abc-c} and divide the result by 
the nonzero expression $\Delta_{[acx]}\cdot \Delta_{[acy]}$. 
Then subtract the resulting equation from~\eqref{eq:abc-b}
and finally divide it by $\Delta_{[abc]}>0$.
We obtain the following subtraction-free expression for 
the minor $\Delta_{[bxy]}$:
$$
\Delta_{[bxy]} = 
\frac{\Delta_{[axy]}\cdot\Delta_{[cxy]}\cdot \Delta_{[abc]}}
{\Delta_{[acx]}\cdot \Delta_{[acy]}}
+
\frac{\Delta_{[axy]}\cdot \Delta_{[bcx]}}{\Delta_{[acx]}}
+\frac{\Delta_{[aby]}\cdot\Delta_{[cxy]}}{\Delta_{[acy]}}.
$$
Since all minors in the right-hand side have already been
expressed subtraction-freely from the minors $\Delta_{\emptyset,J',K'}$
we finally obtain a subtraction-free rational expression for 
$\Delta_{I,J,K} = \Delta_{[bxy]}$.   This finishes the proof of lemma.
\end{proof}

\begin{example}
Let us consider two matrices
$$
A=
\begin{pmatrix}
1 & 1 & 0 & 0 & z_1 & t_1 \\
0 & 1 & 1 & 0 & z_2 & t_2 \\
0 & 0 & 1 & 1 & z_3 & t_3 
\end{pmatrix}
\quad
\text{and}
\quad
B=
\begin{pmatrix}
1 & 0 & 0 & z_1 & t_1 \\
0 & 1 & 0 & z_2 & t_2 \\
0 & 0 & 1 & z_3 & t_3 
\end{pmatrix}.
$$
Suppose that all maximal minors $\Delta_I(A)$ are nonnegative.
According to Lemma~\ref{le:tail-equivalent},  all maximal minors 
$\Delta_J(B)$ should also be nonnegative and, if 
we prescribe which minors $\Delta_I(A)$ are zero and which are 
strictly positive, then we can express each minor $\Delta_J(B)$ 
as a subtraction-free rational expression in terms of the 
$\Delta_I(A)$.
The only nontrivial case is $\Delta_{245}(B) = z_3 t_1 - z_1 t_3$.

Let us show how to express the minor $\Delta_{245}(B)$.
We have $\Delta_{345}(A) = z_1 \geq 0$,
$\Delta_{145}(A) = - z_2 \geq 0$,
$\Delta_{125}(A) =  z_3 \geq 0$,
$\Delta_{346}(A) = t_1 \geq 0$,
$\Delta_{146}(A) = - t_2 \geq 0$,
$\Delta_{126}(A) =  t_3 \geq 0$.
If $\Delta_{235}(A) = z_1 - z_2 + z_3 =0$
or  $\Delta_{236}(A) = t_1 - t_2 + t_3 =0$,
then either $z_1=z_2=z_3 = 0$ or $t_1 = t_2 = t_3 = 0$.
In both cases $\Delta_{245}(B) = 0$.
If both minors $\Delta_{235}(A)$ and $\Delta_{236}(A)$ are nonzero,
then
$$
\Delta_{235}(B) = 
\frac{\Delta_{256}(A)\cdot \Delta_{356}(A)}
{\Delta_{235}(A)\cdot\Delta_{236}(A)}
+
\frac{\Delta_{256}(A)\cdot \Delta_{345}(A)}
{\Delta_{235}(A)}
+
\frac{\Delta_{126}(A)\cdot \Delta_{356}(A)}
{\Delta_{236}(A)}.
$$
Indeed, this equality can be written as
\begin{eqnarray*}
&& z_3t_1 - z_1t_3 = 
\frac{ (z_2t_3 - z_3t_2 - z_1 t_3 + z_3 t_1)\,
       (z_1 t_2 - z_2 t_1 - z_1 t_3 + z_3 t_1)}
     {(z_1 - z_2 + z_3)\,(t_1 - t_2 + t_3)} + \\
&& \qquad 
{}+ \frac{ (z_2t_3 - z_3t_2 - z_1t_3 + z_3 t_1)\,z_1}{z_1 - z_2 + z_3} +
\frac{ t_3 \,(z_1 t_2 - z_2 t_1 - z_1 t_3 + z_3 t_1)}
     {t_1 - t_2 + t_3}.
\end{eqnarray*}
We leave the verification of the last identity as an exercise for the reader.
\end{example}

\begin{lemma}
\label{lem:tail_plus_v}
For two tail-equivalent sequences of $k$-vectors 
$(v_1,\dots,v_r)$, $(v_1',\dots,v_s')$
and a $k$-vector $v$, the sequences
$(v,v_1,\dots,v_r)$ and $(v,v_1',\dots,v_s')$ are tail-equivalent.
\end{lemma}

\begin{proof}
Follows form the cyclic symmetry; 
see Remark~\ref{rem:cyclic_symmetry}.
The maximal minors of the $(r+1+m)$-matrix
$(v,v_1,\dots,v_r,w_1,\dots,w_m)$ are equal to the corresponding 
maximal minors of $(v_1,\dots,v_r,w_1,\dots,w_m,(-1)^{k-1}v)$.
\end{proof}


Let us now prove Lemma~\ref{lem:Le-4}.

\begin{proof}[Proof of Lemma~\ref{lem:Le-4}]
We have $x_1,\dots,x_s=0$ and $x_i >0$ for $i\in [s+1,d]$.
Let $v_1,\dots,v_{d+1}$ be the first column vectors of the matrix $A$ 
rescaled by the positive factors:
$v_i = A_i$ for $i\in[s]\cup\{d+1\}$,
and $v_i = x_i\, A_i$ for $i\in[s+1,d]$.
Then $v_{d+1} = v_{d} - v_{d-1} + v_{d-2} - \cdots + (-1)^{d-s-1} v_{s+1}$.

Also let $u_1,\dots,u_d$ be the linearly independent vectors given by
$u_i = v_i$ for $i\in [s]$, and $u_i = v_i - v_{i-1} + v_{i-2} -
\cdots + (-1)^{i-s-1} v_{s+1}$ for $i\in [s+1,d]$.
Equivalently, these two sequences of vectors are related to each other, 
as follows: 
$v_i = u_i$ for $i\in[s]$,
$v_{i} = u_i + u_{i-1}$ for $i\in[s+1,d]$,
and $v_{d+1} = u_{d}$.
Lemmas~\ref{le:tail-equivalent} and~\ref{lem:tail_plus_v}
imply that $(u_1,\dots,u_d)$ and $(v_1,\dots,v_{d+1})$ are 
tail-equivalent.

Let $\tilde A = (v_1,\dots,v_{d+1},A_{d+2},\dots,A_n)$ and $\tilde C =
(u_1,\dots,u_d,A_{d+2},\dots,A_n)$.  The matrix $\tilde A$ is obtained
from the matrix $A$ by multiplying the columns with indices $j\in [s+1,d]$ 
by the factors $x_j$.  Also note that the matrix $\tilde C$ 
is obtained from $C$ by first multiplying it on the left
the unipotent upper-triangular matrix $(t_{ij})$ 
(which preserves all maximal minors) and
then multiplying the rows by indices $i\in [s+1,d]$ by the factors
$x_i$.  Here $t_{ij} = (-1)^{i-j}$ if $s+1\leq i\leq j \leq d$,
and $t_{ij}=\delta_{ij}$ otherwise.
Thus any maximal minor $\Delta_I(\tilde A)$ equals $\Delta_I(A)$ times
some square-free product of the $x_j$, $j\in[s+1,d]$,
and $\Delta_J(\tilde C) =\Delta_J(C) \cdot \prod_{i\in [s+1,d]} x_j$.

Now the tail-equivalence $(u_1,\dots,u_d)\sim_\mathrm{tail}
(v_1,\dots,v_{d+1})$ implies that, assuming that $A$ belongs
to a fixed cell $S_\M^\tnn$, any maximal minor of $C$ 
can be written as a subtraction-free rational expression 
in the nonzero minors of $A$; and vise versa minors of $A$
can be expressed subtraction-freely in terms of the minors
of $C$ and the $x_i$.  Moreover, note that the minors 
of $\tilde A$ can be written just as sums of some minors of $\tilde C$;
see the beginning of proof of Lemma~\ref{le:tail-equivalent}.
Thus a minor of $A$ can be written as a sum of some minors
of $C$ multiplied by some products of the $x_i$.
\end{proof}

\section{Perfection through trivalency}
\label{sec:perfection}

In this section we show how to simplify the structure of a planar
directed network without changing the boundary measurements 
$M_{ij}$.
\medskip

First of all, if the network $N$ contains any internal sources or sinks,
then we can remove them with all their incident edges.
Indeed, a directed path connecting two boundary vertices can never
pass through an internal source or sink.

If $N$ contains an internal vertex of degree 2, then we can 
remove this vertex, glue two incident edges $e_1$ and $e_2$ 
into a single edge $e$ with weight $x_e = x_{e_1} x_{e_2}$.

If $N$ contains a boundary vertex $b_i$, say, a boundary
source of degree $\ne 1$, with incident edges $e_1,\dots,e_d$,
then we can pull these edges away from $b_i$, create 
a new internal vertex $b_i'$ connected with the source $b_i$ by an 
edge of weight 1 (directed from $b_i$ to $b_i'$), 
and reattach the edges $e_i$ to $b_i'$.  If $d=0$, that is $b_i$
was an isolated vertex in $N$, then we create an additional loop 
at the new vertex $b_i'$, so that $b_i'$ has degree $3$.

Now suppose that $N$ contains an internal vertex $v$ of degree $>3$.
Let $e_1,\dots, e_d$ be the edges incident to $v$ in the clockwise order.
If two adjacent edges $e_i$ and $e_{i+1}$ have the same orientation,
say, towards the vertex $v$, then we can again pull these two edges away from 
the vertex $v$ by creating a new vertex $v'$ and a new edge from
$v'$ to $v$ with weight $1$ and attaching the edges 
$e_i$ and $e_{i+1}$ to the vertex $v'$.  This transformation does not
change the boundary measurements of the network.

\begin{figure}[ht]
\includegraphics[height=1in]{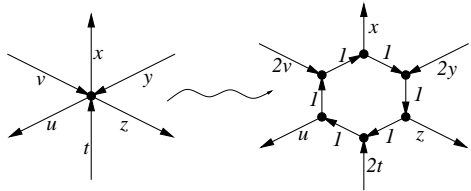}
\caption{Blowing up a vertex into a cycle with trivalent vertices}
\label{fig:vertex_inflation}
\end{figure}

However, if the orientations of the edges $e_1,\dots,e_d$ alternate,
then we cannot make the above reduction.
In this case we can do the following transformation of the network.
Remove the vertex $v$ and create a new cycle $C$ with vertices 
$v_1,\dots,v_d$ and, say, the clockwise orientation of edges.
Assign weights $1$ to all edges of the cycle.
Reattach the edges $e_1,\dots,e_d$ to the new vertices $v_1,\dots,v_d$.
Multiply by 2 the weights of the edges $e_i$ which are oriented towards 
from the cycle; see Figure~\ref{fig:vertex_inflation}.

\begin{lemma}
\label{lem:inflaction_of_vertices}
The above transformation of networks does not change the 
boundary measurements $M_{ij}$.
\end{lemma}

\begin{proof}
First, let us assign variable weights $x_i$ to the edges $(v_i,v_{i+1})$
of the new cycle.
Let $P$ be a directed path in the network $N$ that passes through a vertex $v$.
This path arrives to $v$ through some edge $e_i$ and leaves through $e_j$.
In the new network consider the infinite family
of paths that go exactly as the path $P$ with the only exception
that upon arrival to $C$ through the edge $e_i$ they can make 
several turns around $C$ before the departure through the edge $e_j$.
We may assume that $i<j$. (Otherwise cyclicly shift the labels).
The contribution of a new path from this family to the corresponding boundary 
measurement will be the same as the weight of $P$
times the additional factor $(-1)^s \,2 \,x_i\cdots x_{j-1}(x_1\cdots x_d)^s$,
where $s$ is the number of full turns around $C$ the path makes.   Indeed, 
we get the factor $(-1)^s$ because every new turn changes 
the winding index by $1$; the factor 2 comes from the rescaled weight
of the incoming edge $e_i$ in the network; and the other factors correspond
to the edges of the cycle $C$ the new path passes through.
The total contribution of these paths to the boundary measurement
is the weight of $P$ times $2\,x_i\cdots x_{j-1}\sum_{s\geq 0}
(-1)^s(x_1\cdots x_d)^s=
2\,x_i\cdots x_{j-1}/(1+x_1\cdots x_n)$.
If the path $P$ passes through the vertex $v$ several times, then
for each passage of $P$ through $v$ we will get a similar factor in the new
network.  If we now specialize the weights $x_i$ to $1$ the additional
factors become $2/(1+1) = 1$.  Thus the boundary measurement do not 
change.
\end{proof}

\begin{figure}[ht]
\includegraphics[height=.8in]{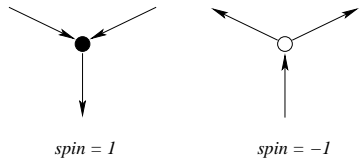}
\caption{Two types of trivalent vertices in a perfect network}
\label{fig:black_white_vertices}
\end{figure}

Thus a planar network can be transformed, without changing the boundary
measurements, into a network with only trivalent internal vertices which are
neither sources nor sinks.  In other words, the obtained network has only 2
types of internal vertices: the vertices with two incoming edges and one
outgoing edge, and the vertices with one incoming edge and two outgoing edges;
see Figure~\ref{fig:black_white_vertices}.  Such trivalent networks belong
to a more general class of networks, defined as follows.

\begin{definition}
\label{def:perfect}
A {\it perfect network\/} $N$ is a planar directed network,
in the sense of Definition~\ref{def:planar_networks}, such that
\begin{enumerate}
\item Each internal vertex $v$ either has exactly one 
incident edge directed away from $v$ (and all other edges 
directed towards $v$), or has exactly one incident edge
directed towards $v$ (all other edges directed away from $v$).
\item All boundary vertices $b_i$ in $N$ have degree 1. 
\end{enumerate}
\end{definition}

The previous discussion implies the following claim.

\begin{proposition}
Any planar directed network can be transformed
without changing the boundary measurements
into a perfect network with trivalent internal vertices.
\end{proposition}

For an internal vertex $v$ of degree $\deg(v)\ne 2$ in a perfect network $N$,
define the number $\col(v) = \col_N(v) \in \{1,-1\}$, called the {\it color\/}
of $v$,
such that, if a vertex $v$ has exactly
one outgoing edge then $\col(v)=1$, and if $v$ has exactly
one incoming edge then $\col(v)=-1$.  
(If $v$ has degree $2$, then we can assign $\col(v)$
in any way.)
We will display internal vertices with $\col(v)=1$ in black color,
and vertices with $\col(v)=-1$ in white color.



\begin{lemma}
\label{lem:knN}
For a perfect network $N$, we have $\sum
\col(v) \cdot (\deg(v)-2) = k-(n-k)$,
where $k$ is the number of sources and $n-k$ is the number 
of sinks in $N$, the sum is over internal vertices $v$.
\end{lemma}

\begin{proof}
Write two numbers at both ends of each directed edge $e=(v_1,v_2)$ in 
the network $N$:
a ``$1$'' at the target $v_2$ and a ``$-1$'' at $v_1$.  Clearly, the total sum
of all written numbers is $0$.  The sum of written numbers at
each internal vertex $v$ equals $\col(v)\cdot (\deg(v)-2)$.  
On the other hand,
at each boundary source we have a ``$-1$'' and at each boundary sink 
we have a ``$1$.''
\end{proof}

\section{Forget the orientation}
\label{sec:edge_switch}

In this section we show that it does not really matter how 
edges are directed in a perfect network.
\medskip

\begin{theorem}
\label{th:switch_of_orientation}
 Let $N = (G,x)$ and $N'=(G',x')$ be two perfect networks 
with $k$ sources and $n-k$ sinks such that:
\begin{enumerate}
\item The graphs $G$ and $G'$ are isomorphic as undirected graphs.
\item Each internal vertex $v$ of degree $\deg(v)\ne 2$ has the
same color $\col_N(v) = \col_{N'}(v)$ in the networks $N$ and $N'$.
\item If $e$ is an edge that has the same orientation in $N$ and in $N'$, 
then $x_e = x_e'$.
If $N$ contains an edge $e=(u,v)$  and $N'$ contains the same undirected 
edge in the opposite orientation $e'=(v,u)$, then $x_{e} = (x'_{e'})^{-1}$.
\end{enumerate}
Then the boundary measurement map $\Mes$ 
maps the networks $N$ and $N'$ to the same point $\Mes(N) =\Mes(N')$ in 
the Grassmannian $Gr_{kn}$.
\end{theorem}

\begin{figure}[ht]
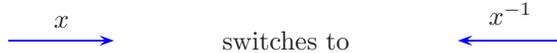

\pspicture(0,0)(210,0)
\psline{->}(0,0)(40,00)
\rput(20,7){$x$}
\rput(105,0){switches to}
\psline{<-}(170,0)(210,00)
\rput(190,10){$x^{-1}$}
\endpspicture
\caption{A switch of edge direction}
\label{fig:edge_switch}
\end{figure}

In other words, if we switch directions of some edges
in a perfect network $N$ at the same time inverting their weights 
(as shown in Figure~\ref{fig:edge_switch})
so that the colors of internal vertices are preserved, then 
the boundary measurement $\Mes(N)\in Gr_{kn}$ will not change.

\begin{figure}[ht]
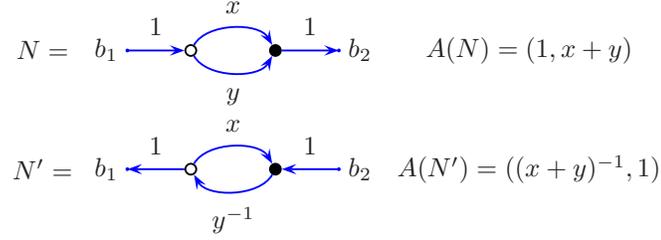

\pspicture(-100,-30)(170,15)
\rput(-90,0){$N =$}
\rput(140,0){$A(N)=(1,x+y)$}
\cnode(-50,0){0}{B1}
\cnode(50,0){0}{B2}
\vwhite(-20,0){V1}
\vblack(20,0){V2}
\nccurve[angleA=0,angleB=180]{->}{B1}{V1}
\aput*{:0}{$1$}
\nccurve[angleA=60,angleB=120]{->}{V1}{V2}
\aput*{:0}{$x$}
\nccurve[angleA=-60,angleB=-120]{->}{V1}{V2}
\bput*{:0}{$y$}
\nccurve[angleA=0,angleB=180]{->}{V2}{B2}
\aput*{:0}{$1$}
\rput(-60,0){$b_1$}
\rput(60,0){$b_2$}
\endpspicture

\pspicture(-100,-15)(170,25)
\rput(-90,0){$N' =$}
\rput(140,0){$A(N')=((x+y)^{-1},1)$}
\cnode(-50,0){0}{B1}
\cnode(50,0){0}{B2}
\vwhite(-20,0){V1}
\vblack(20,0){V2}
\nccurve[angleA=0,angleB=180]{<-}{B1}{V1}
\aput*{:0}{$1$}
\nccurve[angleA=60,angleB=120]{->}{V1}{V2}
\aput*{:0}{$x$}
\nccurve[angleA=-60,angleB=-120]{<-}{V1}{V2}
\bput*{:0}{$y^{-1}$}
\nccurve[angleA=0,angleB=180]{<-}{V2}{B2}
\aput*{:0}{$1$}
\rput(-60,0){$b_1$}
\rput(60,0){$b_2$}
\endpspicture
\caption{Two networks that map to the same point in the Grassmannian}
\label{fig:switch_equivalent_network}
\end{figure}

\begin{example}
Figure~\ref{fig:switch_equivalent_network} shows
two perfect networks obtained from each other by such switch 
of edge directions.
Their boundary measurement matrices are
$A(N) = (1,x+y)$ and $A(N') = ((x+y)^{-1},1)$.
Indeed, we have $M_{12}=x+y$ in $N$, and $M_{21}' = y^{-1}/(1+x\,y^{-1}) = 
(x+y)^{-1}$ in $N'$.  These two matrices represent the same point in 
the Grassmannian $Gr_{1,2} = \mathbb{RP}^1$.
\end{example}

\psset{unit=0.8pt}
\psset{unit=1pt}

Theorem~\ref{th:switch_of_orientation} means that the maximal minors of the
boundary measurement matrices for two networks $N$ and $N'$ obtained by such
switches of orientations of edges are related to each other, as follows.

\begin{corollary}
Let $N$ and $N'$ be two perfect networks 
satisfying conditions of Theorem~\ref{th:switch_of_orientation}.
Let $I'$ be the source set of $N'$.
Then the minor $\Delta_{I'}(A(N))$ is nonzero and,
for any $k$-subset $J\subset [n]$, we have
$$
\Delta_J(A(N')) = \frac{\Delta_J(A(N))}{\Delta_{I'}(A(N))}\,.
$$
\end{corollary}

We now proceed to proving Theorem~\ref{th:switch_of_orientation}
by first checking two special cases: switching edges in a closed
directed cycle and switching edges in a path between two boundary 
vertices.

\begin{remark}
\label{rem:perfect_network}
Perfect networks $N$ have the following important property.
If we pick a directed closed cycle $C$ in $N$, then, for any
vertex $v$ in $C$, the edges that are incident to $v$ and do not 
belong to $C$ are either all directed towards $v$ or all directed away from $v$.
Thus, for any directed path that hits $C$ at some point and
later departs from $C$, the arrival and departure points should be
different from each other. 
The same property holds for a directed path joining two boundary vertices:
another path cannot arrive and depart from it at the same point.
\end{remark}

\begin{lemma}
\label{lem:switch_a_cycle}
Let $N$ be a perfect network.  Let $N'$ be the network obtained
from $N$ by switching the directions of all edges in a closed directed
cycle $C$ and inverting their weights.
Then the network $N'$ has the same boundary measurements $M_{ij}$ as the 
network $N$.
\end{lemma}

\begin{proof}
Let $v_1,\dots,v_d$ be the vertices in the cycle $C$, and let 
$x_i= x_{(v_i,v_{i+1})}$, $i=1,\dots,d$, be the weights of edges
in the cycle.  (Here we assume that $v_{d+1} = v_1$.) 
Consider a collection $\P$ of directed paths connecting two boundary 
vertices in $N$ that are identical 
outside the cycle $C$, but every time when paths enter $C$ they can 
make any number of turns around $C$.
Suppose that paths from this collection enter $C$ through the vertex $v_i$ 
and leave through the vertex $v_j$.  
Note that $i\ne j$ because the network $N$ is perfect;
see Remark~\ref{rem:perfect_network}.
We may assume that $i<j$.  (Otherwise cyclicly shift the labels.)  
Since paths can make any number of turns around $C$, this passage through 
$C$ contributes the factor $x_{i} x_{i+1}\cdots x_{j-1}\sum_{s\geq 0}
(-1)^s(x_1\cdots x_d)^s= x_i\cdots x_{j-1}/(1+x_1\cdots x_n)$
to the sum of terms $\sum_{P\in \P} (-1)^{\wind(P)} \prod_{e\in P} x_e$
from the corresponding boundary measurement.  For each passage of paths in $\P$ 
through $C$ we get a similar factor.
Similarly, we have the collection of paths $\P'$ 
in the network $N'$ identical to paths in $\P$ outside of $C$.
Since we switch directions of the cycle $C$ and invert its weights
in the network $N'$, for each  passage through $C$ with $i$ and $j$ as 
above, we now get the factor
$x_{i-1}^{-1} x_{i-2}^{-1} \cdots x_1^{-1} x_d^{-1}\cdots x_{j}^{-1}/
(1+x_1^{-1}\cdots x_n^{-1})$.  This factor is actually equal to
$x_i\cdots x_{j-1}/(1+x_1\cdots x_n)$.  Thus the total contributions
of paths from $\P$ and $\P'$ to the boundary measurements in $N$ and
$N'$, respectively, are the same.  This implies that all boundary 
measurements in $N$ and $N'$ are the same.

Compare this argument with the proof of
Lemma~\ref{lem:inflaction_of_vertices}.
\end{proof}

We will need the following matrix identities.
All matrices below are square matrices of size $m$.
Let $E_{ij}$ denote the the matrix whose $(i,j)$th entry is
$1$ and all other entries are $0$.  Let $U:=\sum_{i\leq j} E_{ij}$ be the
upper-triangular and $L:=\sum_{i>j} E_{ij}$ be the strictly lower-triangular
matrices filled with $1$'s in all allowed places.  
Also let $E=\sum_i E_{i,\,i+1}$ be the superdiagonal matrix filled with $1$'s,
and let $X=(x_{ij})$ be the matrix filled with the formal variables $x_{ij}$.

\begin{lemma}
\label{lem:XULE}
For the $m\times m$-matrices $X$, $U$, $L$, $E$ as above, define the 
matrices $A=(a_{ij}):= U + UXU + UXUXU + \cdots$,
$B = (b_{ij}) :=L - LXL + LXLXL - \cdots$, 
$C =(c_{ij}) := (1-E-X) + (1-E-X) \,B \,(1-E-X)$,
$P =(p_{ij}) := 1 + (1-E-X) \,B$,
$Q =(q_{ij}) := 1 + B\,(1-E-X)$,
whose coefficients
are power series in the $x_{ij}$.
Then the following identities hold
\begin{enumerate}
\item $A^{-1} = 1- E- X$.
\item $C = a_{1m}^{-1}\, E_{m1}$. 
\item $p_{ij} = \delta_{im} \, a_{1m}^{-1}\, a_{1j}$.
\item $q_{ij} = \delta_{j1}\, a_{1m}^{-1}\, a_{im}$.
\item $b_{ij} = (a_{1j}\, a_{im} - a_{1m} \, a_{ij})\, a_{1m}^{-1}$,
for any $i,j\in [m]$.
\end{enumerate}
Here $1$ denotes the identity matrix and $\delta_{ij}$ is the Kronecker delta.
\end{lemma}

\begin{proof}
(1) \, We have $A= U + UXU + UXUXU+\cdots = U\,(1 - X\cdot U)^{-1}$,
thus $A^{-1} = (1-X\cdot U)\,U^{-1} = U^{-1} - X$.
Notice that $U^{-1} = 1 - E$.  Thus $A^{-1} = 1 - E - X$.

(3), (4), (5)  \, Let us show that these follow from (1) and (2).
We have $P = C\cdot A = a_{1m}^{-1} E_{m1} \cdot A$, which implies (3).  
Similarly, $Q = A\cdot C = a_{1m}^{-1} A \cdot E_{m1}$, which implies (4).
Also $B = A\cdot C\cdot A - A = a_{1m}^{-1}\,(A\cdot E_{m1} \cdot A) - A$.
Thus $b_{ij} =  a_{im}(a_{1m})^{-1} a_{1j} - a_{ij}$, which implies (5).


(2) \
The matrix $C$ expands as the alternating sum 
$$
C = 1 - E - X  + L - EL  - XL - EL - LX + LXL+ ELX + LXL + \cdots,
$$
of all words $w$ in the alphabet $\{L,X,E\}$ such that 
$w$ has the form ${X}{L}{X}{L}{X}\cdots$ or ${L}{X}{L}{X}{L}\cdots $, where
the initial and/or the final letter ``$X$'' can be replaced by 
the letter ``$E$''; and the sign of $w$ is equal to $(-1)^{xe}$, where
$xe$ is the total number of occurrences of letters ``$X$'' and ``$E$''
in $w$.

We need to show that the matrix $C$ has only one nonzero entry $c_{m1} = a_{1m}^{-1}$
in the lower left corner.

Suppose that $i < m$.  Let $C=C'+C''$, where 
$C' =(c_{ij}')  := - (EL) + (EL)X + (EL)E + (EL)XL - (EL)XLX - (EL)XLE + 
\cdots$ is the matrix 
 given by the alternating sum of all words $w$ starting with ``$EL$'';
and $C'' = (c_{ij}'') := 1 - E - X + L - XL - LX + XLX + XLE - LXL + \cdots $ 
is given by the alternating sum of the remaining words starting
with a letter ``$X$'' or ``$L$'' (plus $1-E$).
For a fixed index $i<m$, the contribution of 
a word $w=EL \tilde w$ to $c_{ij}'$ equals to the contribution
of $(E_{i,i+1})( E_{i+1,i} + \sum_{j<i} E_{i+1,j})\,\tilde w = 
\tilde w + L \tilde w$.
Thus the contribution of a word $w$ from $C'$ equals to the sum of 
contributions of two words $\tilde w$ and $L\tilde w$ from $C''$
(obtained from $w$ by erasing one or two initial letters).
Notice that that both words $\tilde w$ and $L\tilde w$ come with
signs opposite to the sign of $w$.  Also note that any word in $C''$
is of the form $\tilde w$ or $L\tilde w$.  Thus all terms from 
$c_{ij}'$ cancel all terms from $c_{ij}''$, implying that
$c_{ij} = c_{ij}'+c_{ij}'' = 0$.

In case when $j>1$ we can use the mirror image of the above argument 
to show that $c_{ij} = 0$.  In this case we need consider letters 
in the end of a word $w$.  Thus $c_{ij} = 0$ unless $(i,j) = (m,1)$.

In remains to prove that $c_{m1} = a_{1m}^{-1}$, or equivalently
$a_{1m}\cdot c_{m1}  = 1$.
We can express $a_{1m}$ as the sum
$$
a_{1m} = 1 + \sum x_{i_1,j_1} + \sum_{j_1\leq i_2} x_{i_1, j_1} x_{i_2, j_2} + 
\sum_{j_1\leq i_2,\, j_2\leq i_3} x_{i_1 j_1} x_{i_2 j_2} x_{i_3 j_3} + \cdots,
$$
over $i_1,j_1,i_2,j_2,\dots \in [m]$ such that
$j_1\leq i_2,\ j_2\leq i_3, \cdots$.
Note that any word $w$ that starts or ends with an ``$E$'' makes 
no contribution to $c_{m1}$.
Thus $c_{m1}$ is given by the alternating sum
$$
c_{m1} = 1 - \sum x_{k_1,l_1} + \sum_{l_1> k_2} x_{k_1, l_1} x_{k_2, l_2} 
- \sum_{l_1> k_2,\, l_2> k_3} x_{k_1 l_1} x_{k_2 l_2} x_{k_3 l_3} + 
\cdots,
$$
over $k_1,l_1,k_2,l_2,\dots\in[m]$ such that 
$l_1>k_2,\ l_2>k_3,\ \cdots$.

Let us use the involution principle to prove the equality
$a_{1m}\cdot c_{m1} = 1$.
The product $a_{1m} \cdot c_{m1}$ can be written as the sum of terms
$x_{i_1,j_i}\cdots x_{i_r j_r} (-1)^s x_{k_1,l_1}\cdots x_{k_s,l_s}$
over pairs of sequences 
$p=((i_1,j_1,i_2,j_2,\dots,i_r,j_r),(k_1,l_1,k_2,l_2,\dots,k_s,l_s))$
such that $j_1\leq i_2,\ j_2\leq i_3,\cdots$ and 
$l_1> k_2,\ l_2> k_3,\ \cdots$.
For $r+s\geq 1$, let us define the map $\iota$
from the set of such pairs of sequences to itself by
$$
\iota(p) := \left\{
\begin{array}{cl}
((i_1,j_1,\dots,i_r,j_r,k_1,l_1),(k_2,l_2,\dots,k_s,l_s))&
\text{if } j_r\leq k_1\text{ or } r=0, \\
((i_1,j_1,\dots,i_{r-1},j_{r-1}) (i_r,j_r,k_1,l_1,\dots,k_s,l_s))&
\text{if } j_r> k_1\text{ or } s=0. 
\end{array}
\right.
$$
Then $\iota$ is an involution, that is $(\iota)^2 = \mathit{id}$.
It preserves the monomial corresponding to $p$ and switches its sign.
Thus all terms in the product $a_{1m}\cdot c_{m1}$, 
except the constant term 1, cancel each other. 
This implies the needed identity.
\end{proof}

\begin{lemma}
\label{lem:switch_a_path}
 Let $N$ be a perfect network, and let $P_0$ be a directed 
path in $N$ from the boundary vertex $b_{i_0}$ to the boundary vertex $b_{j_0}$.
Let $N'$ be the network obtained from $N$ by switching the directions
of all edges in $P$ and inverting their weights.  Let $M_{ij}$ and 
$M_{ij}'$ be the boundary measurements of the networks $N$ and $N'$,
respectively.  Then $M_{i_0,j_0}\ne 0$, and the boundary measurements 
$M_{ij}'$ can be expressed through the boundary measurements $M_{ij}$, 
as follows:
\begin{enumerate}
\item If $(i,j) = (j_0,i_0)$, then $M_{j_0,i_0}' = M_{i_0,j_0}^{-1}$.
\item If $i=j_0$ and $j\ne i_0$, then $M_{j_0,j}' = M_{i_0,j}/M_{i_0,j_0}$.
\item If $i\ne j_0$ and $j=i_0$, then $M_{i,i_0}' = M_{i,j_0}/M_{i_0,j_0}$.
\item If $i\ne j_0$ and $j\ne i_0$, then $M_{ij}' = 
\Delta_{(I\setminus \{i_0, i\})\cup \{j_0,j\}}(A(N))/M_{i_0,j_0}$,
where $I$ is the source set of the network $N$.
\end{enumerate}
This implies that the boundary measurement map $\Mes$ maps 
the networks $N$ and $N'$ to the same point $\Mes(N)=\Mes(N')$ in the Grassmannian.
\end{lemma}

\begin{proof}  The measurement $M_{i_0,j_0}$ is nonzero because 
there is at least one path from $b_{i_0}$ to $b_{j_0}$ in the network $N$,
e.g., the path $P_0$.  Let $P_0 = (b_{i_0},v_1,v_2,\dots,v_m,b_{j_0})$.
We may assume that the path $P_0$ has no self-intersections, because
we can get rid of all self-intersections using Lemma~\ref{lem:switch_a_cycle}.
We may also assume that the weights of all edges in $P_0$ are $1$.
Indeed, we can apply the same gauge 
transformation \eqref{eq:gauge_transformation} to the networks $N$
and $N'$ that transforms the weights of all edges of $P_0$, expect
a single edge $e_0 = (b_{i_0},v_1)$, into $1$'s.
Since the weight of $e_0$ produces same factors in both sides
of all identities (1)--(4), we may assume that this weight in also $1$.

Let us use the notation $b_{\overline{1}},\dots,b_{\overline{m}}$ for 
the vertices $v_1,\dots,v_m$ of the path $P_0$.
For $i,j\in [n]\cup\{\overline{1},\dots,\overline{m}\}$, let
$$
\tilde M_{ij}:= \sum_{P: b_{i}\to b_{j}} (-1)^{\wind(P)} \prod_{e\in P} x_e
$$ 
be the generating function for all directed paths in $N$ from $b_i$ to $b_j$ 
that have no common points with the path $P_0$ (expect the first and the 
last point in $P$).
Note that all paths $P$ should lie in one of the two networks on which 
the path $P_0$ subdivides the network $N$.  
Thus $\tilde M_{ij}$ are the boundary measurements for these smaller networks.
We have $\tilde M_{ii} = 0$, because the network is perfect.
Let us also define the $m\times m$-matrix $X=(x_{ab})$ such that 
$x_{ab} = \tilde M_{\overline{a}, \overline{b}}$ if $a< b$,
$x_{ab} = - \tilde M_{\overline{a}, \overline{b}}$ if $a>b$,
and $x_{ab} = 0$ if $a=b$.

Any path $P'$ in $N$ from $b_{i_0}$ to $b_{j_0}$ first goes along the path
$P_0$; then $P'$ may depart from $P_0$ at some vertex $v_{k_1}$ and then arrive
to $P_0$ at $v_{l_1}$; then $P'$ may depart from $P_0$ again at $v_{k_2}$ and
arrive at $v_{l_2}$; etc.  Let $k_1,l_1,k_2,l_2,\dots,k_s,l_s$ be the indices
of these departure and arrival points.  
Then we have $l_1 < k_2,\ l_2 < k_3,\
\cdots$, because $P'$ coincides with $P_0$ on the segments from $v_{l_i}$ to
$v_{k_{i+1}}$.  The total contribution to $M_{i_0,j_0}$ of all paths $P'$ with
given departure and arrival points is $x_{k_1,l_1} \cdots x_{k_s, l_s}$.
Indeed, each path $P'$ breaks into segments between the departure and arrival
points, which give the factors $x_{k_i,l_i}$.  The extra factor $-1$
in the case when $k_i>l_i$ accounts for an extra cycle in $P'$ that we get 
in this case because the path $P'$ bumps into itself at the arrival
point $v_{l_i}$.  This shows that
$$
M_{i_0,j_0} = 1 + \sum x_{k_1,l_1} + 
\sum_{l_1\leq k_2} x_{k_1,l_1} \, x_{k_2,l_2} + 
\sum_{l_1\leq k_2,\,l_2\leq k_3} x_{k_1,l_1} \, x_{k_2,l_2} \, x_{k_3,l_3} 
+ \cdots,
$$
where the sum is over departure-arrival sequences
$k_1,l_1,k_2,l_2,\dots,k_s,l_s\in [m]$ such that $l_i\leq k_{i+1}$, for
$i\in [s-1]$.

The boundary measurement $M_{j_0,i_0}'$ in the network $N'$ is given by 
a similar expression where we need to sum over departure-arrival
sequences such that $l_i>k_{i+1}$ for $i\in[s-1]$ and we need to 
switch the signs of all $x_{ab}$.

Let $A$, $B$, $C$, $P$, $Q$ be the matrices 
as in Lemma~\ref{lem:XULE}.  Then $M_{i_0,j_0} = a_{1m}$
and $M_{j_0,i_0}' = c_{m1}$.
According to Lemma~\ref{lem:XULE}(2), we have $c_{m1} = a_{1m}^{-1}$.
Thus $M_{j_0,i_0}'= M_{i_0,j_0}^{-1}$. 

Similarly, we can express other boundary measurements of the 
networks $N$ and $N'$ in terms of these matrices.
For $i\ne i_0$ and $j\ne j_0$, we have
$$
\begin{array}{c}
\displaystyle
M_{i_0,j}  = \sum_{c=1}^m a_{1 c} \,\hat M_{\overline{c}, j},
\qquad
M_{i,j_0}  = \sum_{c=1}^m \hat M_{i,\overline{c}}\, a_{c m},\\[.1in]
\displaystyle
M_{j_0,j}' = \sum_{c=1}^m p_{mc}\, \hat M_{\overline{c},j},
\qquad 
M_{i,i_0}' = \sum_{c=1}^m \hat M_{i,\overline{c}}\,q_{c1}.
\end{array}
$$
According to parts (3) and (4) of Lemma~\ref{lem:XULE},
we have $p_{mc} = a_{1c}/a_{1m}$ and $q_{c1} = a_{cm}/a_{1m}$.
Thus $M_{j_0,j}' = M_{i_0,j}/M_{i_0,j_0}$ and 
$M_{i,i_0}' = M_{i,j_0}/M_{i_0,j_0}$. 

We also have 
$$
M_{i,j}  = \epsilon\cdot\hat M_{i,j} + \delta\cdot
\sum_{c,d\in[m]} \hat M_{i,\overline{c}}\,a_{cd} \,
\hat M_{\overline{d},j},
\quad
M_{i,j}' = \epsilon'\cdot \hat M_{i,j} + \delta'\cdot \sum_{c,d\in[m]} 
\hat M_{i,\overline{c}}\,b_{cd}\, \hat M_{\overline{d},j},
$$
where 
$(\epsilon,\delta,\epsilon',\delta') = (1,0,1,0)$ 
if the cords $[b_{i_0},b_{j_0}]$, $[b_i,b_j]$ form a crossing
(see Figure~\ref{fig:crossing_alignment_misalignment}),
$(\epsilon,\delta,\epsilon',\delta') = (1,1,-1,1)$ if the cords 
form an alignment, and
$(\epsilon,\delta,\epsilon',\delta') = (-1,1,1,1)$ 
if the cords form a misalignment.
By Lemma~\ref{lem:XULE}(5),  
$b_{cd} = (a_{1d} a_{cm} - a_{1m}a_{cd})\, a_{1m}^{-1}$.
According to Proposition~\ref{prop:Delta_pi_cn},
the minor 
$\Delta_{(I\setminus \{i_0, i\})\cup \{j_0,j\}}(A(N))$
is equal to $M_{i_0\,j}M_{i,j_0}- M_{i_0,j_0} M_{ij}$
the above cords form an crossing,
to $M_{i_0,j_0} M_{ij}- M_{i_0\,j}M_{i,j_0}$ if the cords
form an alignment, 
or to $M_{i_0,j_0} M_{ij}+ M_{i_0\,j}M_{i,j_0}$ if the cords
form a misalignment.
In all three cases, we get
$M_{ij}' = \Delta_{(I\setminus \{i_0, i\})\cup \{j_0,j\}}(A(N))/M_{i_0,j_0}$,
as needed, which proves (4).

Finally, note that the proved relations (1)--(4) mean that
$\frac{\Delta_{J}(A(N'))}{\Delta_{I'}(A(N'))} =
\frac{\Delta_J(A(N))}{\Delta_{I'}(A(N))}$ for all $k$-subsets $J\subset [n]$
such that $|J\setminus I'|= 1$, where $I'=(I\setminus\{i_0\})\cup\{j_0\}$ is
the source set of the network $N'$.  
Since the $k(n-k)$ quotients of the Pl\"ucker coordinates 
$\frac{\Delta_{J}}{\Delta_{I'}}$, for all such $J$'s, form a 
coordinate system on $Gr_{kn}\setminus\{\Delta_{I'}=0\}$, we deduce that
the matrices $A(N)$ and $A(N')$ represent the same 
point in the Grassmannian $Gr_{kn}$, as needed.
\end{proof}

\begin{proof}[Proof of Theorem~\ref{th:switch_of_orientation}]
Let $H$ be the subset of edges of $G$ whose orientations are switched.
For any internal vertex $v$, the fact that the color $\col(v)$ is preserved
implies that $H$ contains zero or exactly two edges adjacent to $v$.
Thus $H$ is the disjoint union of several cycles and/or paths connecting 
pairs of boundary vertices. 
According to Lemmas~\ref{lem:switch_a_cycle} 
and~\ref{lem:switch_a_path} we can switch the orientations of edges
in these cycles and paths one by one without changing the boundary 
measurements.
\end{proof}

\section{Plabic networks}
\label{sec:plabic}

In this section we will define new weights $y_f$ assigned to faces $f$
of a network, which are obtained from the edge weights $x_e$ by a simple
transformation.
Then we define plabic graphs and networks which are no longer
directed.

\medskip

For a planar graph $G$ (directed or undirected) drawn inside a disk, 
let $V:=V(G)$ be the set of
its internal vertices, $E:=E(G)$ be the set of its edges, and $F:=F(G)$ be the
set of its {\it faces,} that is the regions on which the edges subdivide
the disk.  Let us say that a connected component of $G$ is {\it isolated\/} if
it does not contain a boundary vertex; let $c$ be the number of such isolated
components.  The Euler formula says that $|V|-|E|+|F| = 1+c$.  If $c=0$, then
all faces are simply connected.  If $c\geq 1$, then there are some faces which
are not simply connected, because they contain isolated components inside of
them.  Clearly, if we can always remove isolated components from a directed 
network without affecting the boundary measurements.

\begin{lemma}
\label{lem:gauge_is_nontrivial}
Suppose that $G$ is a planar directed graph without isolated components.
Then the space of directed networks with given graph $G$, modulo 
the gauge transformations, is isomorphic to
$$
\R_{>0}^{E}/\{\text{gauge transformations}\} \simeq
\R_{>0}^{|E|-|V|} = \R_{>0}^{|F|-1}.
$$
\end{lemma}

\begin{proof}
In this case all gauge transformations \eqref{eq:gauge_transformation}
are nontrivial.  Indeed, if we have a gauge transformation
$x_{e}' = x_e \,t_u t_v^{-1}$ such that $x_{e}'=x_e$ for all edges $e$,
then we should have $t_u = 1$ for all vertices $u$ adjacent to the boundary 
vertices $b_i$, then we should have $t_{v} = 1$ for all vertices $v$ adjacent
to the vertices $u$, etc.  Thus $t_v=1$ for all internal vertices $v$.
This implies the first isomorphism.
The Euler formula $|V|-|E|+|F|=1$ implies the second equality.
\end{proof}

\begin{figure}[ht]
\input{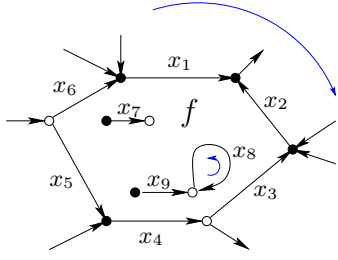}
\caption{A face with weight $y_f=(x_1\, x_2^{-1} x_3^{-1} x_4^{-1} x_5^{-1} 
x_6 )\,(x_7 \,x_7^{-1})\,(x_8^{-1} x_9^{-1}x_9)$}
\label{fig:face_weight}
\end{figure}

Let $f\in F(G)$ be a face in a planar directed network $N=(G,x)$.
The exterior boundary of $f$ consists of some edges $e_1,\dots,e_k$.
If $f$ is not simply connected then it has one or several holes
corresponding to isolated components inside $f$.  
Let $e_{k+1},e_{k+2},\dots,e_l$ be the edges in these holes.
(Note that the same edge might occur twice in this sequence.)
Let us assume that the exterior boundary of $f$ is oriented clockwise
and the boundaries of all holes are counterclockwise.
Let $I^{+}_f\subseteq [l]$ be the index set of edges $e_i$ whose orientations
in the graph $G$ agree with the orientation of the boundary of $f$,
and let $I^{-}_f = [l]\setminus I^{+}$ be the index set of edges
whose orientations disagree.
Then we define the {\it face weight\/} $y_f$ of the face $f$ in $N$ as
$$
y_f := \prod_{e_i\in I^+_f} x_{e_i} \cdot 
\prod_{e_j\in I^{-}_f} x_{e_j}^{-1},
$$
see Figure~\ref{fig:face_weight}.
Note that, if we switch directions of some edges inverting their
weights (as shown on Figure~\ref{fig:edge_switch}), then the face 
weights $y_f$ will not change.
There is one relation for the faces weights, namely,
$\prod_{f\in F} y_f = 1$.  Indeed, this product includes 
exactly one $x_e$ and exactly one $x_e^{-1}$ for each edge $e\in E$.

Let $\R_{>0}^{F-1}\simeq \R_{>0}^{|F|-1}$ be the set of 
$(y_f)_{f\in F}\in \R_{>0}^F$ such that $\prod y_f = 1$.

\begin{lemma}
\label{lem:E/g=F}
For a planar directed graph $G$,
the map $(x_{e})_{e\in E} \mapsto (y_f)_{f\in F}$ defined as above
gives the isomorphism
$$
\R_{>0}^{E}/\{\text{gauge transformations}\}\simeq \R_{>0}^{F-1}.
$$
\end{lemma}

\begin{proof}  
Let us first assume that $G$ has no isolated components.
We prove the claim by induction on $|F|$.
By Lemma~\ref{lem:gauge_is_nontrivial} we already know that both sides
have the same dimension.  Thus it is enough to show that the kernel 
of the map consists of a single point. 
Let $(x_e)_{e\in E}\in \R_{>0}^E$ be a point that maps to 
$(y_f)_{f\in F}$ with all $y_f = 1$.  We need to show that one
can transform all edge weights $x_e$ into $1$'s by the gauge transformations.
Let us pick a face $f_0\in F$ adjacent to the boundary of the disk.
The boundary of $f_0$ contains a segment $(b_{i},b_{i+1})$ of the boundary
of the disk, and a path $P=(b_{i},v_1,\dots,v_k,b_{i+1}]$.  Applying
the gauge transformations at the vertices $v_1,\dots,v_k$, we can transform 
all weights of the $k+1$ edges in the path $P$, except a single edge, 
into $1$'s.
But, since we have $y_{f_0} =1$ and $y_{f_0}$ is the product of weights 
of edges in $P$ (some of them might be inverted), we deduce that the weight
of the last edge in $P$ should also be $1$.   We can now apply 
the induction hypothesis to the smaller graph $G'$, which is obtained from 
$G$ by removing the face $f_0$ and replacing the segment $(b_i,b_{i+1})$
of the boundary of the disk with the path $P$.  By induction, all weights 
of edges in $G'$ can be transformed into $1$'s by the gauge transformations
at the internal vertices of the graph $G'$, that is the vertices
in $V\setminus \{v_1,\dots,v_k\}$.
Thus all weights $x_e$ of edges in $G$ are now transformed into $1$'s,
as needed.

In the general case, we prove the claim by induction on the number
$c$ of isolated components.  We have already established the base case $c=0$.
Suppose that $c\geq 1$.  Let $G'$ be the graph obtained from $G$ by 
removing an isolated component $G''$ located in some face $f$ of $G$. 
Let $F',E',F'',E''$ be the face and edge sets of these two
graphs, and also let $\tilde E''\subset E''$ be the set of internal
edges of $G''$.  The face weights for $G'$ and $G''$ are the same 
as for the graph $G$ with a single exception:  
the face weight $y_f$ for $G$ is obtained 
by multiplying the corresponding weight for $G'$ by the product of 
all face weights for $G''$.
By induction, we have
$\R_{>0}^{E'}/\{\text{gauge transformations}\}\simeq \R_{>0}^{F'-1}$.
If $G''$ has no faces, i.e., it is a tree, then 
$\R_{>0}^{E}/\{\text{gauge transformations}\}\simeq 
\R_{>0}^{E'}/\{\text{gauge transformations}\}\simeq \R_{>0}^{F'-1}
=\R_{>0}^{F'-1}$.
Otherwise, $G''$ breaks into several disjoint subgraphs $G_1,\dots,G_k$
(connected by paths),
each of which is formed by a cycle with 
a small graph inside, so that each face of $G''$ belong to one 
of these graphs $G_i$.
Applying the induction hypothesis to each graph $G_i$, we deduce that
$\R_{>0}^{\tilde E''}/\{\text{gauge transformations}\}\simeq 
\R_{>0}^{F''-k}$.  Since the weights of the boundary
edges in $G''$ can be arbitrary, we get $\R_{>0}^{E''}
/\{\text{gauge transformations}\}\simeq \R_{>0}^{F''}$.
Thus $\R_{>0}^{E}/\{\text{gauge transformations}\} \simeq
\R_{>0}^{F'-1}\times  \R_{>0}^{F''} \simeq \R_{>0}^{F-1}$, as needed.
\end{proof}

\begin{remark}  
As we have already mentioned, isolated components do not affect
the boundary measurements.  The reason that we are 
considering graphs that might have isolated
components will be clear below, when we define certain transformations
of graphs.  Even if an original graph does not have isolated components,
after performing several moves we might create such components.
\end{remark}


The boundary measurement map
$\Mes_G:\R_{>0}^{E}/\{\text{gauge transformations}\}\to Gr_{kn}$,
\eqref{eq:M_G_map}, now transforms into the map 
$\R_{>0}^{F-1}\to Gr_{kn}$.
Below we will use the face weights $y_f$ instead of the edge weights $x_e$.
With these weights we no longer need to care about gauge 
transformations, and we no longer need to invert the weights when 
we switch edge directions; see Figure~\ref{fig:edge_switch}.

We can define the boundary measurements of a directed network $N$ 
purely in terms of 
the face weights $y_f$ without using the edge weights $x_e$.
A directed path $P$ without self-intersections
that connects two boundary vertices $b_i$ and $b_j$ 
subdivides the disk into two parts: the part that is on the right
side of $P$ and the part on left of $P$ (as we go from $b_i$ to $b_j$).
We define $wt(P,y)$ as the product of the weights $y_f$ for the faces in 
the right part of $P$.  Also, for a clockwise (resp., counterclockwise) 
closed cycle $C$, define $wt(C,y)$ as the product of the $y_f$ for the faces 
$f$ inside (resp., outside) of $C$.
Now, for an arbitrary path $P$ from $b_i$ to $b_j$, we can remove
cycles $C_1,\dots,C_k$ and reduce it to a path $\tilde P$ without
self-intersections.  In this case, define $wt(P,y) = wt(\tilde P)\cdot
\prod C_i$.

\begin{lemma}  
For any directed path $P$, we have $wt(P,y) = \prod_{e\in P} x_e$.
\end{lemma}

\begin{proof}  It is enough to prove this equality of a path
$P$ without self-intersections and also prove that
$wt(C) = \prod_{e\in C} x_e$ for a closed cycle $C$.
In all cases, the product of the $y_f$ for all faces $f$ in the right side 
of $P$ (resp., inside/outside of $C$) includes exactly one weight $x_e$ 
and exactly one weight $(x_e)^{-1}$ for all edges $e$ 
in the corresponding areas expect for the edges $e$ that belong to $P$
(resp., to $C$).  The remaining terms give the needed product
of edge weights.
\end{proof}

Thus the (formal) boundary measurements of a network can be 
defined as 
$$
M_{ij} =\sum_{P:b_i\to b_j} (-1)^{\wind(P)}wt(P,y),
$$
cf., \eqref{eq:Mij}.

\begin{definition}  A {\it planar bicolored graph,} or simply
a {\it plabic graph\/}  is a planar (undirected) graph $G$, defined 
as in Definition~\ref{def:planar_networks} but without orientations
of edges, such that each boundary vertex $b_i$ is incident
to a single edge, together with a function $\col:V\to\{1,-1\}$
on the set $V$ of internal vertices.
As before, we will display vertices with $\col(v)=1$ in black color,
and vertices with $\col(v)=-1$ in white.

A {\it plabic network\/}
$N=(G,y)$ is plabic graph $G$ together with positive real weights $y_f>0$ 
assigned to faces $f$ of $G$ such that $\prod y_f = 1$. 

A {\it perfect orientation\/} of a plabic graph or network 
is a choice of orientation of its edges such that each internal
vertex $v$ with $\col(v)=1$ is incident to exactly one edge 
directed away from $v$; and each $v$ with $\col(v)=-1$ is incident 
to exactly one edge directed towards $v$.
A plabic graph or network is called {\it perfectly orientable\/} 
if it has a perfect orientation.

Let us say that a plabic graph or network has {\it type\/} $(k,n)$ 
if its has $n$ boundary vertices and 
$k+(n-k) = \sum_{v\in V} \col(v)\,(\deg(v)-2)$.
\end{definition}

\begin{remark}  One can think about plabic graphs as some kind
of ``Feynman diagrams,'' where the black and white vertices 
represent certain ``elementary particles'' of two types and
edges represent ``interactions'' between these particles.
\end{remark}

According to Lemma~\ref{lem:E/g=F}, plabic networks with a choice of a perfect 
orientation correspond to perfect networks modulo gauge transformation.
Theorem~\ref{th:switch_of_orientation} says that two perfect networks
$N$ and $N'$ that correspond to two orientations of the same  plabic 
network should map into the same point $\Mes(N)=\Mes(N')\in Gr_{kn}$.
Lemma~\ref{lem:knN} says a perfect orientation of a plabic graph of
type $(k,n)$ should have $k$ sources and $n-k$ sinks.
Thus the boundary measurement map $\Mes$ gives a well defined map
$$
\tilde \Mes: \{\text{perfectly orientable plabic networks of type $(k,n)$}\}
\to Gr_{kn}^\tnn.
$$
For a perfectly orientable plabic graph $G$, we have the
induced map on the set $\R_{>0}^{F(G)-1}$ of plabic networks
with the given graph $G$:
$$
\tilde \Mes_{G}:\R_{>0}^{F(G)-1} \to Gr_{kn}^\tnn.
$$

Note that not any plabic network is perfectly orientable.
For example, a plabic network that contains an isolated component
with a single vertex is not perfectly orientable.

For a plabic graph $G$ of type $(k,n)$ and 
a perfect orientation $\mathcal{O}$ of $G$, let $I_{\mathcal{O}}
\subset[n]$ be the $k$-element source set of this orientation.
Define the {\it matroid\/} of $G$ as the set of the $k$-subsets  
$I_{\mathcal{O}}$ for all perfect orientations:
$$
\M_G:=\{I_{\mathcal{O}}\mid \mathcal{O} 
\text{ is a perfect orientation of $G$}\}.
$$

\begin{proposition}
\label{prop:plabic_matroid}
For any perfectly orientable plabic graph $G$,
the collection $\M= \M_G$ is a totally nonnegative matroid.
The boundary measurement map $\tilde \Mes_G$ sends plabic networks
with the graph $G$ into the totally nonnegative Grassmann cell $S_{\M}^\tnn$
associated with this matroid:
$$
\tilde \Mes_G(\R_{>0}^{F(G)-1}) \subseteq S_{\M}^\tnn.
$$
\end{proposition}

In Section~\ref{sec:matroids_perms}
we will prove that this inclusion
is actually the equality;
see Corollary~\ref{cor:Im=Cell}.

\begin{remark}
This proposition gives way to combinatorially determine the totally 
nonnegative Grassmann cell corresponding to a planar network without 
any calculations of the boundary measurements.  Indeed, first transform 
a network into a perfect network (see Section~\ref{sec:perfection}),
then transform it into a plabic network (forget orientations
of edges but remember colors of vertices), and calculate
the matroid $\M_G$.
\end{remark}

\begin{proof}[Proof of Proposition~\ref{prop:plabic_matroid}]
\end{proof}

\begin{example}
Figure~\ref{fig:perfect_orientations} shows five perfect orientations
of a plabic graph $G$.  Recording their source sets, we obtain
the following matroid with five bases $\M_G = \{\{1,4\},\, \{1,2\},\, 
\{1,3\},\, \{2,4\},\, \{3,4\}\}$.
\begin{figure}[ht]
\input{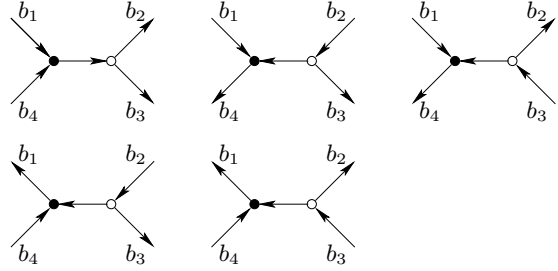}
\caption{Perfect orientations of a plabic graph}
\label{fig:perfect_orientations}
\end{figure}
\end{example}

Let us give two additional combinatorial descriptions of the matroid $\M_G$
in terms of paths and in terms of matchings.

Let us fix a perfect orientation $\mathcal{O}$ of $G$ and let 
$I=I_\mathcal{O}$.  Define the {\it path matroid\/} 
$\M_G^p$ as the set of 
$k$-element subsets $J\subset[n]$ such that the boundary vertices
$\{b_i\mid i \in I\setminus J\}$ can be connected with the 
boundary vertices $\{b_j\mid j\in J\setminus I\}$ by 
a family of pairwise noncrossing directed paths in the graph $G$ 
with edge orientation $\mathcal{O}$.

Let us say that a plabic graph $G$ is {\it bipartite\/}
if any edge in $G$ joins two vertices of different colors
(assuming that the colors of all boundary vertices are white).
Note that we can easily make any plabic graph bipartite by
inserting vertices of different color in the middle of unicolored
edges.
A {\it partial matching\/} in such graph $G$ is a subset $M$ 
of edges such that each internal vertex is incident to exactly one edge
in $M$.  (But the boundary vertices can be
incident to one or zero edges.)
Let $I_M\subseteq [n]$ be the set of indices
$i$ such that $b_i$ belongs an edge from $M$.
Define the {\it matching matroid\/} 
$$
\M_G^m:=\{I_M\mid M \text{ is a partial matching of } G\}.
$$ 

\begin{lemma}
\label{lem:path_matchings}
For any plabic graph $G$, we have $ \M_G^p = \M_G$.
Also, if $G$ is a bipartite plabic graph, then $\M_G^m = \M_G$.
\end{lemma}

\begin{proof}
Any perfect orientation $\mathcal{O}'$ is obtained from the fixed
perfect orientation $\mathcal{O}$ by switching edge directions
in a family of noncrossing directed paths between boundary vertices
or closed cycles.
This implies that $\M_G^p = \M_G$.
For a perfect orientation $\mathcal{O}$ in bipartite plabic graph,
let $M$ be the set of edges in $G$ directed from a black vertex to a white
vertex.  The map $\mathcal{O}\mapsto M$ is a bijection between
perfect orientations and partial matchings, which implies
that $\M_G^m=\M_G$.
\end{proof}

\section{Transformations of plabic networks}
\label{sec:transformations_plabic}

In this section we define several local transformation of 
plabic networks.  In all transformations below, we change 
a small fragment in a network and sometimes
change weights of adjacent faces.  The weights of 
remaining faces are not changed.
We will call the first three transformations (M1)--(M3)
the {\it moves,} and the next three transformations
(R1)--(R3) the {\it reductions.}
We think of the moves as invertible transformations of networks,
which we can perform in both directions.
On the other hand, we will perform reductions only in one
direction in order to simplify the structure of a network.
Essentially, the only nontrivial transformation of networks
is the {\it square move;} all other moves and reductions play
an auxiliary role.

\begin{figure}[ht]
\input{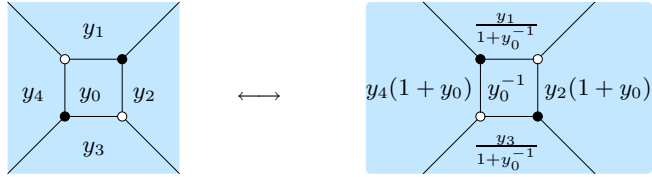}
\caption{Square move}
\label{fig:square_move}
\end{figure}

\smallskip
\noindent
(M1) {\sc square move.}  If a network has a square 
formed by four trivalent vertices whose colors alternate
as we go around the square, then we can switch colors of 
these four vertices 
and transform the weights of adjacent faces
as shown on Figure~\ref{fig:square_move}.
In other words, if $y_0$ is the weight 
of the face inside the square and $y_1,y_2,y_3,y_4$
are weights of the four adjacent faces, then we
transform these weights as follows
\begin{equation}
\label{eq:weight_sq_move}
y_0' = y_0^{-1},\
y_1'= \frac{y_1}{1+y_0^{-1}},\
y_2' = y_2\,(1+y_0),\quad
y_3' = \frac{y_3}{1+y_0^{-1}},\
y_4'= y_4\,(1+y_0).
\end{equation}
In the case when some of the four areas marked by $y_1,y_2,y_3,y_4$ in
Figure~\ref{fig:square_move} belong to the same face (connected outside
of the shown fragment of the network), say,
if $y_1 = y_2$ are in the same face, then its
weight changes to $y_1\frac{1}{1+y_0^{-1}} (1+y_0)  = y_1 y_0$. 

\begin{figure}[ht]
\input{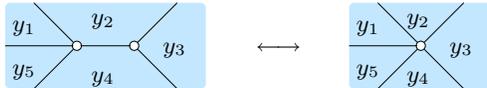}
\caption{Unicolored edge contraction}
\label{fig:edge_contraction}
\end{figure}

\noindent
(M2) {\sc unicolored edge contraction/uncontraction.}  If a network 
contains an edge with two vertices of the same color, then we can contract
this edge into a single vertex with the same color;
see Figure~\ref{fig:edge_contraction}.
The face weights $y_f$ are not changed. 
On the other hand, we can also uncontract a vertex into an edge
with vertices of the same color.

\begin{figure}[ht]
\input{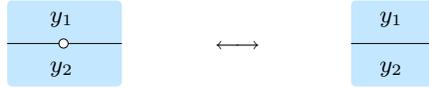}
\caption{Middle vertex insertion/removal}
\label{fig:vertex_removal}
\end{figure}

\noindent
(M3) {\sc middle vertex insertion/removal.} 
If a network contains a vertex of degree $2$,
then we can remove this vertex and glue the incident edges together;
see Figure~\ref{fig:vertex_removal}.
The face weights $y_f$ are not changed.
On the other hand, we can always insert a vertex (of any color)
in the middle of any edge.

\begin{figure}[ht]
\input{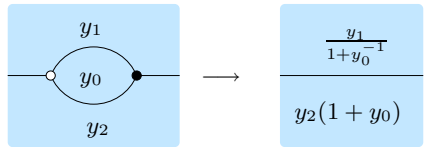}
\caption{Parallel edge reduction}
\label{fig:edge_reduction}
\end{figure}

\noindent
(R1) {\sc parallel edge reduction.}
If a network contains two trivalent vertices of different colors
connected by a pair of parallel edges,
then we can remove these vertices and edges,
glue the remaining pair of edges together, 
and transform the face weights as shown on Figure~\ref{fig:edge_reduction}. 
If $y_1=y_2$ correspond to the same face, then we change its weight
to $y_1\frac{1}{1+y_0^{-1}}(1+y_0) = y_1 y_0$.

\begin{figure}[ht]
\input{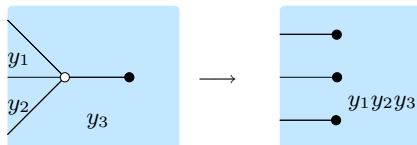}
\caption{Leaf reduction}
\label{fig:leaf_reduction}
\end{figure}

\noindent
(R2) {\sc leaf reduction.}  If a network contains a vertex (leaf) $u$ incident
to a single edge $e = (u,v)$ which in turn is incident to edges 
$e_1,\dots,e_k$, $k\geq 2$, such that $\col(u)=-\col(v)$,
then we can remove the vertex $u$ together with the edge $e$,
disconnect the edges $e_1,\dots,e_k$, and assign the color equal to $\col(u)$ 
to all newly created vertices of the edges $e_1,\dots,e_k$;
see Figure~\ref{fig:leaf_reduction}.
If this operation joins several faces $f_1,\dots,f_l$ into a single face $f$,
then its weight $y_f$ should be the product of weights of $f_1,\dots,f_l$.
(Note that is it possible all faces were already connected outside the shown
fragment so that this transformation does not reduce the number of faces.
In this case it creates new isolated components.)

\begin{figure}[ht]
\input{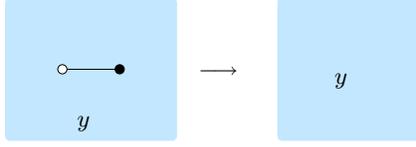}
\caption{Dipole reduction}
\label{fig:dipole_reduction}
\end{figure}

\noindent
(R3) {\sc dipole reduction.} 
If a network contains an isolated component $C$ that consists of a pair
of vertices of different colors connected by an edge, then we can remove $C$
from the network; see Figure~\ref{fig:dipole_reduction}.
The edge weights are not changed.

If forget about face weights in the above moves and reductions
we obtain corresponding transformations of plabic graphs.

Let us say that two plabic networks (or graphs) are {\it move-equivalent\/}
if they can be obtained from each other by moves (M1)--(M3).
Similarly, two plabic networks (graphs) are {\it move-reduction-equivalent\/}
if they can be transformed into the same network (graph) by moves (M1)--(M3) 
and reductions (R1)--(R3).

\begin{theorem}
\label{th:moves=boundary}
Let $N$ and $N'$ be two perfectly orientable plabic
networks of the same type $(k,n)$.  Then the boundary measurement
map $\tilde \Mes$ maps them into the same point $\tilde \Mes(N) = \tilde \Mes(N')$ in
the Grassmannian $Gr_{kn}$ if and only if the networks $N$ and $N'$
are move-reduction-equivalent.
\end{theorem}

We will prove this theorem in Section~\ref{sec:matroids_perms}.
In one direction this claim can be verified by direct calculation.

\begin{lemma}
Suppose that a plabic network $N'$ is obtained from $N$ by performing
a move {\rm (M1)--(M3)} or a reduction {\rm (R1)--(R3).}
Then $N'$ is perfectly orientable if and only if $N$ is perfectly orientable.
If this is the case, then $\tilde \Mes(N) = \tilde \Mes(N')$.
\end{lemma}

\begin{proof}  
It is quite easy to check in all six cases that a perfect orientation of $N$
gives a perfect orientation of $N'$, and vise versa.
The only nontrivial transformation of networks is the square move (M1).
Let us check that the boundary measurement map is invariant under this
transformation.  Let us pick two perfect orientations of the networks $N$ and $N'$, 
say, the orientations whose parts in the transformed fragment are shown on 
Figure~\ref{fig:square_move_arrow}
and which are are identical everywhere else.
The transformation of face weights in the square move corresponds to the 
following transformations of weights of the four edges that form the square:
$$
x_1' = \frac{x_3x_4}{x_2+x_1 x_3 x_4},\
x_2' =  x_2 + x_1 x_3 x_4,\
x_3' = \frac{x_2 x_3}{x_2 + x_1 x_3 x_4},\
x_4' = \frac{x_1 x_3}{x_2 + x_1 x_3 x_4},
$$
where we assume that the remaining edge weights are not changed.
\begin{figure}[ht]
\input{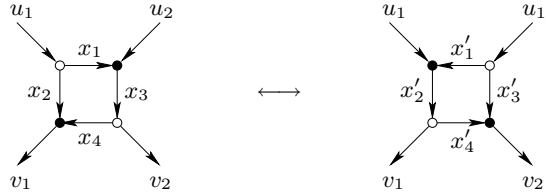}
\caption{Square move in directed networks}
\label{fig:square_move_arrow}
\end{figure}

Then for the both oriented network fragments shown on 
Figure~\ref{fig:square_move_arrow}, sums over paths
from $u_i$ to $v_j$ in $N$ and $N'$ are the same:  $x_2 + x_1 x_3 x_4 = x_2'$,
$x_1 x_3 = x_2' x_4'$, $x_3 x_4 = x_1' x_2'$, $x_3 = x_3' + x_1' x_2' x_4'$.
Thus all boundary measurements in the both networks should be the same,
implying $\tilde \Mes(N) = \tilde \Mes(N')$.

Similarly, for the the parallel edge reduction (R1), let us pick
a perfect orientation of edges as in Figure~\ref{fig:edge_reduction_arrow}.
The transformation of face weights in this move correspond to the 
transformation of edge weights given by $x_1' = x_1(x_2+x_3)x_4$.
Clearly, this transformation does not change the boundary measurements.

\begin{figure}[ht]
\input{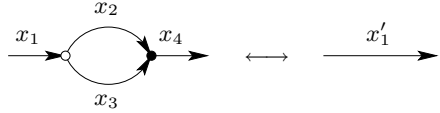}
\caption{Parallel edge reduction in directed networks}
\label{fig:edge_reduction_arrow}
\end{figure}

For the remaining moves and reductions is it clear
that $\tilde \Mes(N) = \tilde \Mes(N')$.
\end{proof}

\begin{remark}
Figure~\ref{fig:square_move_arrow} shows just one of several possible
square moves in {\it directed\/} networks.  If we pick another perfect
orientation of the edges (inverting the edge weights $x_i$ and $x_i'$ whenever
we switch orientations), then we get another legitimate square move.
\end{remark}

There are several special kinds of networks to which we
can easily transform any plabic network using moves (M2)--(M3)
and reductions (R1)--(R3).

\begin{figure}[ht]
\input{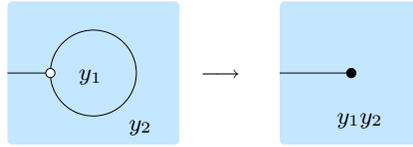}
\caption{Loop reduction}
\label{fig:loop_reduction}
\end{figure}

{\it Loop removal:}
We can remove all loops (i.e., edges whose both ends are the same vertices)
 from a network.  By uncontracting some edges we can make all loops
attached to trivalent vertices.
Then we apply the loop reduction shown on Figure~\ref{fig:loop_reduction}.
This reduction follows from parallel edge reduction (R1).
Indeed, insert an additional vertex of different color 
to the loop using (M3), then uncontract this vertex into an edge using (M2),
and apply parallel edge reduction (R1).
Let us call a network without loops {\it loopless.}

{\it Leaf removal:}
We can also easily get rid of all leaves in a plabic network,
except the leaves connected to boundary vertices.
(Let us call these special unremovable leaves the {\it boundary leaves.})
Indeed, if a leaf $u$ is attached to a vertex $v$ of the same color
then we can just contract the edge $(u,v)$ using move (M2).
If colors of $u$ and $v$ are different, then we can remove the 
edge $(u,v)$ using leaf reduction (R3) if $\deg(v)\geq 3$.
(If $\deg v=2$ then we can remove $v$ using (R2)
and if $\deg(v)=1$ then we can remove it by dipole reduction (R3).)
Then similarly treat all newly formed leaves, etc.
Let us call a network {\it leafless\/} if it has no leaves,
except the boundary leaves.

{\it Contraction:}
Any plabic network can be transformed
into a network that has no unicolored edges, no non-boundary leaves,
and no vertices of degree $2$.  
Indeed, first remove all non-boundary leaves, then contract all unicolored
edges, then remove all new vertices of degree $2$, then contract all new 
unicolored edges, etc.  Let us call such networks {\it contracted}.

{\it Making the graph trivalent:}
On the other hand, we can first get rid of non-boundary leaves and vertices of
degree $2$, and then uncontract all vertices of degree $>3$ by replacing them
with trivalent trees.  We obtain a network with all trivalent internal vertices
(expect boundary leaves).  Let us call such plabic networks {\it trivalent.}

\begin{corollary}  Any plabic network can be transformed (without changing the
boundary measurements) into a loopless contracted network.   On the other hand,
it can also be transformed into a loopless trivalent network.  
\end{corollary}

\begin{definition}
\label{def:reduced_plabic}
Let us say that a plabic network (or graph) is {\it reduced\/} if
it has no isolated connected components and 
there is no network/graph in its move-equivalence class 
to which we can apply a reduction (R1) or (R2).
A {\it leafless\/} reduced network/graph is a reduced network/graph
without non-boundary leaves.
\end{definition}

We will see that reduced networks are exactly the networks without
isolated component with the minimal possible value of
$E-V = F-c-1$ in its move-reduction-equivalence
class, where $V,E,F,c$ as in Section~\ref{sec:plabic}.

Note that if a network has no isolated components then all its move-equivalent networks
have no isolated components, so there is no chance to apply dipole reduction (R3).
In many cases it will be more convenient to use leafless reduced networks, to
which we can easily transform any reduced network by the leaf removal
procedure as described above.  An arbitrary reduced network can obtained from 
a leafless one
by uncontracting vertices into trees of the same color, then maybe inserting
vertices of different color in the middle of new edges, then maybe
uncontracting some of them into trees, etc.  So that we can grow a bicolored
tree of special kind at each vertex.  

Since we can never perform a leaf reduction
(R2) in a leafless graph, we obtain the following claim.

\begin{lemma}
\label{lem:leafless_reduced}
A leafless plabic graph without isolated components is reduced if and only
if it is impossible to transform it by the moves {\rm (M1)--(M3)} into a graph
where we can perform parallel edge reduction {\rm (R1).}
\end{lemma}

The next claim is the main result on reduced plabic graphs.

\begin{theorem}
\label{th:reduced_plabic_parametrization}
Let $G$ be a reduced plabic graph.
Then $G$ is perfectly orientable
and the map $\tilde \Mes_G:\R_{>0}^{F(G)-1}
\to S_\M^\tnn$ gives a subtraction-free rational parametrization of
the corresponding totally nonnegative Grassmann cell $S_\M^\tnn$.
In particular, the dimension of $S_\M^\tnn$ equals $|F(G)|-1$.

For any cell $S_\M^\tnn$ there is a reduced plabic graph $G$ such 
that $\tilde \Mes_G$ is a parametrization of $S_\M^\tnn$.

Any two different parametrizations $\tilde \Mes_G$ and $\tilde \Mes_{G'}$ of
the same cell $S_\M^\tnn$ can be obtained from each other
by the moves {\rm (M1)--(M3).}
\end{theorem}

We will prove this theorem in Section~\ref{sec:matroids_perms}.

\begin{remark}
\label{rem:non-reduced}
For a non-reduced plabic graph $G$ without isolated components,
the map $\tilde \Mes_G$ is either undefined (when $G$ is not
perfectly orientable) or this map is not injective.
Indeed, if we can do a reduction (possibly after performing 
some moves), then we can decrease the number of needed 
parameters.
\end{remark}


\section{Trips in plabic graphs}

In this section we give a criterion when a plabic graph is reduced
and describe move-equivalence classes of reduced graphs.
The results of this section can be related to Thurston's
work~\cite{Thurs} on triple diagrams (see Remark~\ref{rem:thurston_moves}).
\medskip

\begin{figure}[ht]
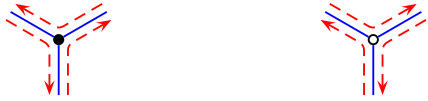

\psset{unit=.7pt}
\begin{center}
\pspicture(0,-40)(170,36)

\vblack(0,0){V1}
\cnode[linewidth=0](0,-30){0}{B1}
\cnode[linewidth=0](-25.98,15){0}{B2}
\cnode[linewidth=0](25.98,15){0}{B3}
\nccurve[angleA=90,angleB=-90]{-}{B1}{V1}
\nccurve[angleA=150,angleB=-30]{-}{V1}{B2}
\nccurve[angleA=30,angleB=-150]{-}{V1}{B3}
\psline[linearc=5,linestyle=dashed,linecolor=red]{->}(5,-30)(5,-2.89)(28.48,10.67)
\psline[linearc=5,linestyle=dashed,linecolor=red]{->}(23.48,19.33)(0,5.77)(-23.48,19.33)
\psline[linearc=5,linestyle=dashed,linecolor=red]{<-}(-5,-30)(-5,-2.89)(-28.48,10.67)

\vwhite(170,0){V1}
\cnode[linewidth=0](170,-30){0}{B1}
\cnode[linewidth=0](144.02,15){0}{B2}
\cnode[linewidth=0](195.98,15){0}{B3}
\nccurve[angleA=90,angleB=-90]{-}{B1}{V1}
\nccurve[angleA=150,angleB=-30]{-}{V1}{B2}
\nccurve[angleA=30,angleB=-150]{-}{V1}{B3}
\psline[linearc=5,linestyle=dashed,linecolor=red]{<-}(175,-30)(175,-2.89)(198.48,10.67)
\psline[linearc=5,linestyle=dashed,linecolor=red]{<-}(193.48,19.33)(170,5.77)(146.52,19.33)
\psline[linearc=5,linestyle=dashed,linecolor=red]{->}(165,-30)(165,-2.89)(141.52,10.67)
\endpspicture
\end{center}
\psset{unit=1pt}
\caption{Rules of the road for trips in plabic graphs}
\label{fig:trip_rules}
\end{figure}

For an (undirected) plabic graph $G$, a {\it trip\/} is a 
directed path $T$ in $G$ such that 
\begin{enumerate}
\item
$T$ either joins two boundary vertices ({\it one-way trip})
or it is a closed cycle that contains none of the boundary vertices ({\it round-trip}).
\item If $T$ arrives to an internal vertex $v$ with incident edges 
$e_1,\dots,e_d$ (in the clockwise order) though the edge $e_i$, 
then it should leave $v$ through the edge $e_{i-\col(v)}$. 
(Here indices $i$ in $e_i$ are taken modulo $d$.) 
In other words, $T$ obeys the following ``rules of the road'':
turn right at a black vertex, and turn left at 
a white vertex; see Figure~\ref{fig:trip_rules}.
\end{enumerate}

Note that these trips in {\it undirected\/} plabic graphs have nothing to do
with paths in directed networks that we used in the definition of the boundary
measurements.

If two trips pass along the same edge in the same directions then
they should be identical.  If a trip passes along the same edge twice in the
same direction then it is a round-trip.
Thus each edge of $G$ belongs to exactly two trips or to one trip with a 
self-intersection at this edge.

Each plabic graph $G$ with $n$ boundary vertices defines 
the {\it trip permutation\/} $\pi_G\in S_n$ such that $\pi_G(i)=j$ 
whenever the trip that starts
at the boundary vertex $b_i$ ends at the boundary vertex $b_j$.

The following claim is established by direct examination.

\begin{lemma}
Let a plabic graph $G'$ is obtained from $G$ by one of the moves {\rm (M1)--(M3).}
Then $\pi_G = \pi_{G'}$.
In other words, each one-way trip in $G$ is transformed into a one-way 
trip in $G'$ with the same end points.
Also each closed trip in $G$ is transformed into a closed trip in $G'$.
\end{lemma}

Note that reductions (R1)--(R2) (unlike the moves) change the trip permutation.

For an edge $e$ with vertices of different colors, we say that two trips
(resp., one trip) that pass(es) through the edge $e$ in two different
directions have/has an {\it essential intersection\/} (resp., {\it essential
self-intersection}) at this edge $e$.  All other (self)-intersections are
called {\it inessential.} Note that in an inessential (self)-intersection the
trips do not cross but rather touch each other.  We can can always remove an
inessential (self)-intersection by performing edge contraction/uncontraction moves (M2).

Let us say that two trips $T_1\ne T_2$ in $G$ have a {\it bad double
crossing\/} if they have two essential intersections at edges $e_1$ and $e_2$
such that both trips $T_1$ and $T_2$ are directed from $e_1$ to $e_2$.


The following theorem gives a criterion when a plabic graph is reduced.

\begin{theorem}
\label{th:reduced_criterion}
Let $G$ be a leafless plabic graph 
without isolated connected components.
Then $G$ is reduced if and only if the following conditions hold
\begin{enumerate}
\item $G$ has no round-trips.
\item $G$ has no trips with essential self-intersections.
\item $G$ has no pair of trips with a bad double crossing.
\item If $\pi_G(i)=i$ then $G$ has a boundary leaf attached
to the boundary vertex~$b_i$.
\end{enumerate}
\end{theorem}

Note that condition (1) implies that $G$ contains no
isolated connected components.

\begin{definition}  A {\it decorated permutation\/} 
$\pi^:=(\pi,\col)$ is a permutation
$\pi\in S_n$ together with a coloring function $\col$ 
from the set of fixed points $\{i\mid \pi(i) = i\}$ to $\{1,-1\}$.
That is a decorated permutation is a permutation with fixed
points colored in two colors.
\end{definition}

Suppose for a moment that we have already established 
Theorem~\ref{th:reduced_criterion}.  Then we can decorate
the trip permutation $\pi_G$ of a {\it reduced\/} plabic graph $G$ 
by coloring each fixed point $\pi_G(i) = i$ to the color $\col(i) := \col(v)$, 
where is $v$ is the boundary leaf attached to the boundary vertex $b_i$.
This gives the {\it decorated trip permutation\/} $\pi_G^:$ of $G$.

\begin{theorem}
\label{th:moves=permutations}
Let $G$ and $G'$ be two reduced plabic graphs with the same number
of boundary vertices.  Then the following claims are equivalent:
\begin{enumerate}
\item $G$ can be obtained from $G'$ by moves {\rm (M1)--(M3).}
\item These two graphs have the same decorated trip permutation $\pi_G^:=\pi_{G'}^:$.
\end{enumerate}
\end{theorem}

Since any reduced plabic graph can be transformed into a leafless
graph by moves (M2) and (M3), it is enough to prove 
Theorem~\ref{th:moves=permutations} for leafless reduced graphs.

We will also need the following auxiliary claim.

\begin{figure}[ht]
\input{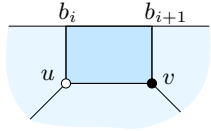}
\caption{Square boundary face}
\label{fig:square_boundary_face}
\end{figure}

\begin{lemma}
\label{lem:square_face}
Let $G$ be a reduced plabic graph such that $\pi_G$ has no fixed points.  
Let $i<j$ be two indices such that $\pi_G(i) =j$ or $\pi_G(j)=i$
and there is no pair $i',j'\in[i+1,j-1]$ such that $\pi_G(i')=j'$.
Then one can transform $G$ by moves
{\rm (M1)--(M3)} into a graph with a square face 
that is attached to the boundary interval $[b_i,b_{i+1}]$
and has two other internal vertices $u$ and $v$;
see Figure~\ref{fig:square_boundary_face}.  

Moreover, if $\pi_G(i)=i+1$ and $\pi_G(i+1) =i$, and if $G$
is leafless and has no vertices of degree $2$,
then the boundary vertices $b_i$ and $b_{i+1}$ are connected in $G$ by an edge.
\end{lemma}

Note that, for any pair $i<j$ such that $\pi(i)=j$ or $\pi(j)=i$, either this pair
itself satisfies the condition of Lemma~\ref{lem:square_face} or there is another pair
inside the interval $[i,j]$ that satisfies the condition of this lemma.

\begin{proof}[Proof of Theorems~\ref{th:reduced_criterion} 
and~\ref{th:moves=permutations}]
We will prove Theorem~\ref{th:reduced_criterion}, 
Theorem~\ref{th:moves=permutations}, and Lemma~\ref{lem:square_face}
all together by induction on the number faces in $G$.

Let us assume that $G$ is a leafless plabic graph with $1$ face and without
isolated connected components.  Then $G$ consists of trees attached to the boundary
vertices.  Since $G$ cannot have non-boundary leaves, that means the $G$
contains only of boundary leaves attached to all boundary vertices.  All such
graphs are reduced because it is impossible to further reduce them and clearly
they satisfy conditions in Theorem~\ref{th:reduced_criterion}.
Theorem~\ref{th:moves=permutations} and Lemma~\ref{lem:square_face} are also
trivial in this case.  This gives the base of induction.

Let us now assume that $G$ is a leafless plabic graph with $>1$ faces and
without isolated components.  By the induction hypothesis we have already
established Theorems~\ref{th:reduced_criterion},
and~\ref{th:moves=permutations}, and Lemma~\ref{lem:square_face} 
for all graphs with fewer number of faces than in $G$.

If a graph $G$ is not reduced then after performing some moves (M1)--(M3) we
should be able to reform a parallel edge reduction (R1);
see Lemma~\ref{lem:leafless_reduced}.
It is easy to see
that right before the reduction one of the conditions (1)--(4) in
Theorem~\ref{th:reduced_criterion} fails.  Indeed, suppose that we get a pair
of parallel edges between two vertices $u$ and $v$ of different color.  Then
contract all unicolored edges and consider several cases: if $\deg(u)=2$ or
$\deg(v)=2$ then (1) fails; if $\deg(u)=3$ or $\deg(v)=3$ then (2) or (4)
fails; if $\deg(u),\deg(v)>3$ then (3) fails.  Note that moves (M1)--(M3) can
never remove a failed condition (1)--(4).  Since moves (M1)--(M3) are
invertible, that means that in the original graph $G$ we also get a failed
condition (1)--(4).  This proves Theorem~\ref{th:reduced_criterion} in one
direction.

Let us prove Theorem~\ref{th:reduced_criterion} in the other direction.
Suppose that one the conditions (1)--(4) in Theorem~\ref{th:reduced_criterion}
fails.  Let us show that it is possible transform the graph $G$ by the moves
(M1)--(M3) to a graph where we can perform reduction (R1). 

In all cases a segment of a problematic trip $T_1$ (or a pair segments in a pair of
trips $T_1$ and $T_2$) surrounds an area $A$ that consists of some faces in $G$.  
Let us assume that $A$ is the area between closest essential (self)-intersections
so that there are no other essential (self)-intersections of the $T_i$'s inside $A$.
Let us show that we can always undo all inessential (self)-intersections
of the $T_i$'s using moves (M2).  Let us first contract all unicolored edges in $G$.
Let $v$ be a vertex, say, with $\col(v) = 1$ and with incident edges 
$e_1,\dots,e_d$ (in the clockwise order).  Suppose we have an inessential
(self)-intersection of the trip(s) $T_i$ at $v$.  That means that one
of these trips arrives to $v$ through the edge $e_j$ 
and then leaves through $e_{j-1}$ (according to the ``rules of the road'')
and then (the same or the other) trip arrives to $v$ through $e_l$ and leaves 
through $e_{l-1}$.  Note that all edges $e_{j}, e_{j-1}, e_{l}, e_{l-1}$ should
be different.  (Otherwise we get an essential intersection.)
We might get some other pairs of edges at the vertex $v$ corresponding to
other passages through this vertex.  Since all these pairs of edges $(e_j, e_{j-1})$,
$(e_l, e_{l-1})$, \dots are disjoint, we can always uncontract the vertex $v$ into a trivalent
tree such that all these segments of paths no longer intersect.
This argument shows that we may assume that the area $A$ is homeomorphic to a
disk.

Let us remove all vertices of degree 2 from $G$ and contract all unicolored
pairs of edges on the boundary of the disk $A$ so that colors of vertices now
alternate as we go around the disk $A$.  Remind that we assume that $G$
is leafless.

Suppose that the area $A$ contains only one face.  Let us consider several cases.

I.  Suppose that $A$ is surrounded by a round-trip, say, a clockwise round-trip. 
If there is a white vertex $v$ on the boundary of $A$, then all edges incident to $v$
should lie {\it inside\/} $A$.  Thus there is a tree attached to $v$, which is impossible
because we assume that $G$ is leafless.  So $v$ should have degree $2$, which is again
impossible because have removed all such vertices from $G$.  That mean that there 
is no white vertices on the boundary of $A$.  Thus there is only one black vertex on
the boundary of $A$.  So we get a loop; see Figure~\ref{fig:loop_reduction}.
This loop can be transformed into a pair of parallel edges.
Then we can apply reduction (R1), as needed.
Similarly, if $A$ is surrounded by a counterclockwise round-trip, 
there is no black vertices on its boundary and thus $A$ is again formed by a loop.

II.  Suppose that $A$ is surrounded by a segment of a trip with an 
essential self-intersection, say, a clockwise segment.  
Again in this case there is no white vertices on the boundary $A$,
except the white vertex that belongs to the essential self-intersection.  
So the boundary of $A$ has at most two vertices, one black and one white.
That means that we either get a loop or get a pair of parallel edges.
In both cases, we can perform reduction (R1).

III. The case when the area $A$ is surrounded by a segment of the trip that starts 
and ends at the same boundary vertex $b_i$, is exactly the same as the case II.

IV.  Suppose that $A$ is surrounded by a pair of segments $S_1$ and $S_2$
of two trips between a bad double crossing;
cf.\ Figure~\ref{fig:double_crossing_bi}.
Assume that the segment $S_1$
is directed clockwise, then $S_2$ is directed counterclockwise.
The same argument as above shows that there are no white vertices on $S_1$
and similarly there are no black vertices on $S_2$.
%
That means that $S_1$ consists of a single black vertex and $S_2$ consists of 
a single white vertex.  Again the boundary of $A$ is formed
by a pair of parallel edges, so that we can perform reduction (R1).


Let us now assume that $A$ has more that two faces.
Let $\tilde b_1,\dots,\tilde b_r$ be the vertices on the boundary of $A$
(in the clockwise order) that have incident edges {\it inside\/} the area $A$.
Using uncontractions (M2), we can transform the vertices $\tilde b_i$ into trivalent vertices,
i.e., for each $\tilde b_i$ there is now exactly one incident edge that lies inside of $A$.
Then we obtain a smaller a plabic graph $\tilde G$ (with at least one trip)
inside the area $A$. 
If $\tilde G$ is not reduced then so is $G$.  Assume that $\tilde G$ is reduced.
By the induction hypothesis we have already 
established all needed claims for the graph $\tilde G$.

Let $\pi_{\tilde G}\in S_r$ be the trip permutation of the graph $\tilde G$.
Note that $\tilde G$ has no boundary leaves, because we assumed that $G$ is
leafless.  Thus the trip permutation $\pi_{\tilde G}$ has no fixed points.

If $A$ is given by one segment of a trip, then all $\tilde b_i$'s have the same color
(all white if the segment is clockwise, and all black if the segment is 
counterclockwise).  For each pair $\tilde b_i$ and $\tilde b_{i+1}$,
the graph $G$  has one vertex of the opposite color between these two vertices.
By the induction hypothesis Lemma~\ref{lem:square_face} holds for the
graph $\tilde G$.  So we can transform the graph $\tilde G$ to
a graph $\tilde G'$ that has a square face $f$ attached to some
boundary segment $[\tilde b_l,\tilde b_{l+1}]$. 
Furthermore, in the case of a trip with an essential self-crossing 
we may assume that the boundary segment $[\tilde b_l,\tilde b_{l+1}]$ does not contain this self-crossing.
(Just label the vertices $\tilde b_i$ so that the self-crossing is between $\tilde b_r$ and $\tilde b_1$.)
In the graph $G$, the face $f$ includes $5$ vertices
(because there is one extra vertex of $G$ between $\tilde b_l$ and 
$\tilde b_{l+1}$).  Note that, as we go around $f$, the colors of the $5$ 
vertices cannot change more than $4$ times.  That means that we can always 
merge at least two vertices of $f$ together by move (M2).  
If the colors of vertices change less than $4$ times, that means that
we can transform a graph $G$ into a graph where we can perform
a reduction (R1).  If the colors change exactly $4$ times, then 
(after some uncontractions for vertices that are not trivalent) we can perform 
a square move for the face $f$.
In the obtained graph $G'$, the corresponding face now lies outside the area $A'$
formed by the corresponding trip.  Thus the number of faces in $A'$ is strictly less
than the number of faces in $A$.  Then we can repeatedly apply the same procedure 
until we get a graph where a problematic trip surrounds exactly one face.
This case was already considered above.

If the area $A$ is given by two segments $S_1$ (clockwise) and $S_2$ (counterclockwise)
of two trips with a bad double crossing, then the $\tilde b_i$ that belong
to $S_1$ are white and the $\tilde b_j$ that belong to $S_2$ are black;
see Figure~\ref{fig:double_crossing_bi}.
Suppose that there are two vertices $\tilde b_i$ and $\tilde b_j$ that both 
lie on the side $S_1$ or both lie on $S_2$ such that they are connected by a trip in $\tilde G$.
We may assume that they are closest such vertices, so that the condition
of Lemma~\ref{lem:square_face} holds.  Thus again we will can transform
$\tilde G$ to a graph that contain a square face attached to 
a segment in $S_1$ or in $S_2$.  Then we can decrease the number of faces inside $A$
as above.
Otherwise all trips in $\tilde G$ that start at $S_1$ should end at $S_2$, 
and vise versa.  (In particular, both $S_1$ and $S_2$ contain the same number
of vertices.)  Let $\tilde b_s$ be the last vertex in $S_2$ and be
the first vertex in $\tilde b_{s+1}$; see Figure~\ref{fig:double_crossing_bi}. 

If $l=\pi_{\tilde G}(s+1)\ne s$ or $l=\pi_{\tilde G}^{-1}(s+1)\ne s$,
then the pair $(l,s+1)$ satisfies condition of Lemma~\ref{lem:square_face}.
Thus again we can transform the adjacent face adjacent to $[\tilde b_l,\tilde b_{l+1}]$
into a square and then apply square move (M1) and
reduce the number of faces inside of $A$, as above. 

In the remaining case we have $\pi_{\tilde G}(s) = s+1$ are $\pi_{\tilde G}(s+1)=s$.
According to the second part of Lemma~\ref{lem:square_face},
in this case the vertices $\tilde v_{s}$ and $\tilde v_{s+1}$ are connected by 
an edge $e$.
Thus the graph $G$ contains a square face below the edge $e$,
so again we can perform a square move (M1) at this face
and reduce the number of faces inside $A$.
Thus in all cases we can repeatedly decrease the number of faces in $A$
until we get an area with one face.
This finishes the inductive step in the proof of 
Theorem~\ref{th:reduced_criterion}.

\begin{figure}[ht]
\input{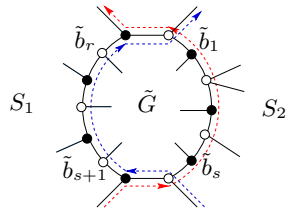}
\caption{A bad double crossing with a graph $\tilde G$ inside}
\label{fig:double_crossing_bi}
\end{figure}

Let us now prove Theorem~\ref{th:moves=permutations} and
Lemma~\ref{lem:square_face}.
In one direction Theorem~\ref{th:moves=permutations}
is straightforward because the moves (M1)--(M3) never change the
decorated permutation $\pi_G^:$.  
Let us assume that $G$ and $G'$ two reduced graphs such that
$\pi_G^: =\pi_{G'}^:$.  Also assume that $|F(G)|\geq |F(G')|$.
We will show that the graphs $G$ and $G'$ can be obtained from each other by the
moves (M1)--(M3).  By Theorem~\ref{th:reduced_criterion} 
(which is already proved for $G$ and $G'$) all
fixed points in $\pi_G$ should correspond to boundary leaves in the graphs $G$
and $G'$, which should be located in the same positions and should have same
colors (given by the decoration in the decorated permutation $\pi_G^:$).
Thus, without loss of generality, we may assume that $\pi_G$ has no fixed points.

Let us pick a pair $i<j$ satisfying the condition of Lemma~\ref{lem:square_face}.
Let us assume that $\pi_G(i) =j$.  (The other case is completely analogous.)
By Theorem~\ref{th:reduced_criterion} the trip $T$ in $G$ from $b_i$ to $b_j$ has no 
essential self-crossings, and we can eliminate inessential self-crossings 
by the moves (M2), as above. 
This trip $T$ subdivides the graph $G$ into two 
smaller graphs $G_1$ and $G_2$,  where $G_1$ is the graph containing
the boundary segment $(b_i,b_{i+1},\dots,b_j)$.
We may assume that the graph $G$ satisfies the property that $G_1$ has
the smallest possible number of faces for all graphs in the move-equivalence
class of $G$.  (Otherwise transform $G$ to such graph by the moves.)
Let us remove vertices of degree $2$ from $G$ and contract unicolored edges in $T$.

Let $\pi_{G_1}$ be the trip permutation for the graph $G_1$.
The trip permutation $\pi_{G_1}$ contains no fixed points.  
(They would correspond to boundary leaves in $G_1$, but we assume
that $G$ is leafless.)
If there is a trip in $G_1$ that starts and ends at $T$,
then using the same argument as above, we can transform $G_1$ into
a graph that has a square face attached to a segment of $T$, then apply
the square move at this face, and decrease the number of faces in $G_1$,
which contradicts to our assumption of minimality of $G_1$.
This means that all trips in $G_1$ that start at a (black) vertex on $T$ 
should end at one of the vertices $b_{i+1},\dots,b_{j-1}$, and vise versa.
Let $\tilde b_{i+1},\dots,\tilde b_{j-1}$ be the black vertices on $T$
as we go from $b_i$ to $b_j$.

If the trip of $G_1$ that starts at $b_{i+1}$ ends at a vertex $\tilde b_l \ne \tilde b_{i+1}$
or the trip that ends at $b_{i+1}$ starts at a vertex  $\tilde b_l\ne \tilde b_{i+1}$,
again by Lemma~\ref{lem:square_face} we can transform $G_1$ to a graph that
has a square face attached to $[\tilde b_{l-1},\tilde b_l]$, then perform 
a square move, and thus reduce the number of faces in $G_1$, which contradicts
to minimality of $G_1$.
Thus vertices $b_{i+1}$ and $\tilde b_{i+1}$ are connected by $2$ trips in both ways.
By Lemma~\ref{lem:square_face}, this implies that $b_{i+1}$ 
and $\tilde b_{i+1}$ should be connected by an edge.
Let us remove this edge and apply the same argument again to show that 
the vertices $b_{i+2}$ and $\tilde b_{i+2}$ are also connected an edge in $G$, etc.
There should be a white vertex in $T$ between two adjacent 
black vertices $\tilde b_{l}$, $\tilde b_{l+1}$ in $T$.
If needed, we can also insert white vertices between $b_i$ and  $\tilde b_{i+1}$
and between $\tilde b_{j-1}$ and $b_j$, and then make all white vertices in $T$
trivalent.  (So now $G$ might contain one or two leaves.)
Therefore we may assume that the trip $T$ in the graph $G$ goes along the 
boundary and has the form as shown in Figure~\ref{fig:boundary_trip}.
Note that if $T$ has this boundary form, then the trip permutation $\pi_{G_2}$
of the graph $G_2$ that lies outside of this boundary strip is uniquely
determined by the trip permutation $\pi_G$ of $G$.

\begin{figure}[ht]
\input{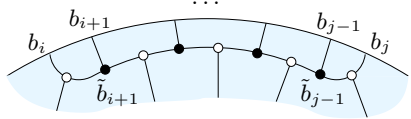}
\caption{Boundary trip from $b_i$ to $b_j$}
\label{fig:boundary_trip}
\end{figure}

We can also transform by the moves the second graph $G'$ (with $\pi_{G'} =
\pi_G$) to the form where the trip $T'$ from $b_i$ to $b_j$ has exactly the same
boundary form as in Figure~\ref{fig:boundary_trip}.  Thus its part $G_2'$ outside
the boundary strip has exactly the same trip permutation 
$\pi_{G_2'} = \pi_{G_2}$ as $G_2$.  By induction hypothesis, this implies
that $G_2$ and $G_2'$ are move-equivalent.  Then $G$ and $G'$ are move equivalent.
This finishes the inductive step for Theorem~\ref{th:moves=permutations}.

In order to prove Lemma~\ref{lem:square_face}, notice that in the above 
boundary trip from $b_i$ and $b_j$ shown on Figure~\ref{fig:boundary_trip}
the boundary segment $[b_i,b_{i+1}]$ is adjacent to a square face.
This proves the first claim of 
Lemma~\ref{lem:square_face}.

Now assume that $\pi_G(i)=i+1$ and $\pi_G(i+1) = i$.  Again we can transform $G$
so that the trip $T$ from $b_i$ to $b_{i+1}$ has the boundary form as above.  
Let us remove leaves and vertices of degree $2$.
The trip $T$ contains at most one vertex --- a white vertex $v$.
If $T$ contain a vertex $v$, then the other trip $T_2$ from $b_{i+1}$ to $b_i$ goes away from 
$T$ at this vertex and then returns back to $T$ at the same vertex $v$.
Let $A$ be the area surrounded by the part of $T_2$ between two visits of $v$.
Applying the same argument as above in proof of Theorem~\ref{th:reduced_criterion} (the case of
trip with a self-crossing), we can reduce the number of faces inside $A$ one by one
until there is only one face left.  Now all other (possible) vertices on the boundary of $A$
can be contracted into a single black vertex.   So the graph contains a pair of parallel
edges and we can perform reduction (R1), which is impossible because we assume that $G$
is reduced.  Thus there is no vertex $v$ on the trip $T$, that is $b_i$ and $b_{i+1}$
are connected by an edge $e$.  By performing moves (M1)--(M3) we can only 
insert middle vertices into $e$ and grow some trees at these vertices.
But since we assume that the original graph $G$ is leafless and has no vertices of degree $2$,
the edge $e$ should be present in $G$. 
This proves the second claim of Lemma~\ref{lem:square_face}.
This finishes the inductive proof of 
Theorem~\ref{th:reduced_criterion}, 
Theorem~\ref{th:moves=permutations}, and Lemma~\ref{lem:square_face}.
\end{proof}

Let us say that a {\it singleton\/} is an isolated connected component with a
single vertex and no edges.

\begin{lemma}
\label{lem:any_plabic_reduced}
Any plabic graph $G$ can be transformed 
by moves {\rm (M1)--(M3)} and reductions {\rm (R1)--(R3)}
into a reduced plabic graph possibly together with some singletons.
\end{lemma}

\begin{proof}
Let us first apply the leaf removal procedure to $G$.
If the obtained graph has a non-singleton isolated component,
then there is a round-trip $T$ in this component.
We can decrease the number of faces inside $T$ and
then apply a reduction as in the above proof. 
Repeatedly applying reductions we end up with a reduced
graph possibly together with some singletons.
\end{proof}

\section{Alternating strand diagrams}
\label{sec:alt_strand}

One can transform reduced plabic graphs into the following objects.

\begin{definition}
An {\it alternating strand diagram\/} consists of $n$ directed curves,
called the {\it strands,} which are drawn inside a disk and connect pairs
of the boundary vertices $b_1,\dots,b_n$ such that the following conditions hold:
\begin{enumerate}
\item For any boundary vertex $b_i$, there is exactly one strand that enters $b_i$
and exactly one strand that leaves $b_i$.
\item No three strands can intersect at the same point.
\item There is a finite number of pairwise intersection points of the strands.
All intersection points are transversal, i.e., the tangent vectors to the strands
at the intersection points are independent.
\item Let $S$ be a strand, and let $v_1=b_j,v_2,\dots,v_l = b_j$ be all points 
on $C$ where it intersects with other strands $S_1,\dots, S_l$ as it goes from $v_i$ to $v_j$.
(The same strand might occur several times in this sequence.)
Then the orientations of the strands $S_1,\dots,S_l$ at the points $v_1,\dots,v_l$ alternate.
In other words, if, say, $S_1$ is oriented at $v_1$ from left to write (with respect to $S$)
then $S_2$ is oriented at $v_2$ from right to left, $S_3$ is oriented from left to right, etc.
\item A strand has no self-intersections, except the case when the strand is a loop
(clockwise or counterclockwise) attached to a boundary vertex $b_i$.
\item If two strands  have two intersection points $u$ and $v$, then
one of these strands is oriented from $u$ to $v$ and the other is oriented from $v$ to $u$.
\end{enumerate}
Needless to say that alternating strand diagrams are considered modulo homotopy.
\end{definition}

Each alternating strand diagram $D$ has the {\it decorated strand
permutation\/} $\pi^:_D$ such that $\pi_D^:(i) = j$ whenever $D$ has a strand
from $b_i$ to $b_j$ and if there is a counterclockwise (resp., clockwise) loop
attached to $b_i$, then the fixed point $i$ is colored in black $\col(i) = 1$
(resp., white $\col(i)=-1$).

\begin{figure}[ht]
\input{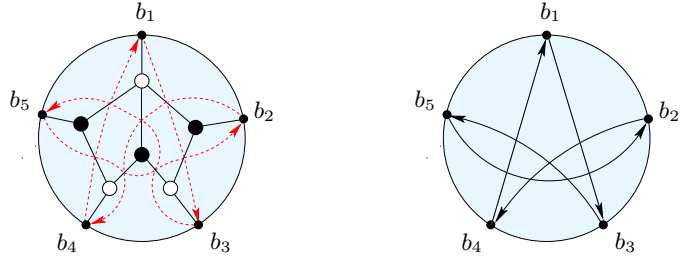}
\caption{A reduced plabic graph and the corresponding alternating strand diagram}
\label{fig:alternating_diagram}
\end{figure}

The right side of Figure~\ref{fig:alternating_diagram} displays an example of alternating strand
diagram $D$ with the strand permutation $\pi_D(i) = i+2 \pmod 5$.  Note that
the diagram where the points $b_i$ and $b_{\pi(i)}$ are connected straight
chords is not an alternating strand diagram.

{\it Faces\/} of an alternating strand diagram are the regions on which the strands
subdivide the disk.  There are 3 types of such faces: {\it clockwise\/} (whose boundary
is directed clockwise), {\it counterclockwise\/} (whose boundary is directed counterclockwise),
and {\it alternating\/} (where directions of the strands alternate when we go around
the face.  All {\it boundary faces,} i.e., faces adjacent to a segment of the boundary of
the disk, are alternating.  For a clockwise or counterclockwise face, all its adjacent faces are 
alternating.  On the other hands, an alternating face is adjacent to both clockwise 
and counterclockwise faces (in an alternating order).

A reduced plabic graph $G$ can be easily transformed into an alternating strand
diagram $D_G$ as follows:

\begin{enumerate}
\item Remove all non-boundary leaves and contract all unicolored edges in $G$.
\item Draw a dot in the middle of each edge connecting two internal vertices,
and also draw dots at the boundary vertices $b_1,\dots,b_n$.
\item  For any black internal vertex $v$, connect the dots $d_1,\dots,d_l$
(in clockwise order) on the incident edges by new edges $(d_2,d_1), \dots, (d_l,d_{l-1}), (d_1,d_l)$
oriented counterclockwise.
\item  Similarly, for any white internal vertex $v$, connect the adjacent dots $d_1,\dots,d_l$
by new edges $(d_1,d_2), \dots, (d_{l-1},d_l), (d_l,d_1)$
oriented clockwise.
\item  We obtain a new directed graph where all internal vertices (internal dots) have degree $4$
and all boundary vertices $b_i$ have degree $2$.  The strands are directed paths in this graph connecting
the boundary vertices $b_i$ and intersecting each other at the dots.
\end{enumerate}

Theorems~\ref{th:reduced_criterion} and~\ref{th:moves=permutations}
(upon some observation of Figure~\ref{fig:alternating_diagram})
imply the following result.

\begin{figure}[ht]
\input{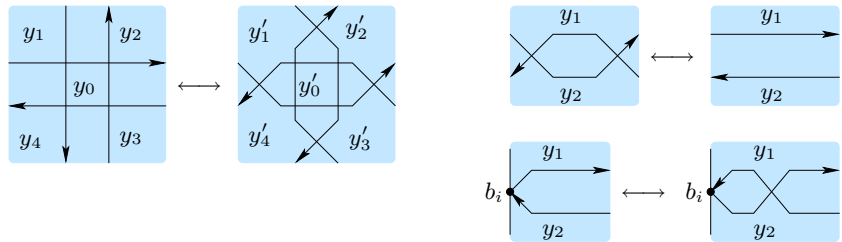}
\caption{Moves of alternating strand diagrams}
\label{fig:alternating_moves}
\end{figure}

\begin{corollary}
The diagram $D=D_G$ constructed from a reduced plabic graph $G$ as above is 
an alternating strand diagram.  The map $G\mapsto D_G$ is a bijection between 
reduced plabic graphs (without non-boundary leaves and unicolored edges)
and alternating strand diagrams.
Trips in $G$ are transformed into strands of $D$.
Thus diagram $D$ has the same decorated strand permutation 
$\pi_D^:=\pi_G^:$ as the decorated trip permutation of $G$.

Black (resp., white) vertices in $G$ correspond to counterclockwise (resp.,
clockwise) faces of $D$.  Faces of $G$ correspond to alternating faces of
$D$.

Two alternating strand diagrams have the same decorated strand permutation
if and only if they can be obtained from each other by the moves shown 
in Figure~\ref{fig:alternating_moves}.
\end{corollary}

In alternating strand diagrams weights are assigned only to
alternating faces.  These weights are transformed as shown on Figure~\ref{fig:alternating_moves},
where the subtraction-free transformation $(y_0,\dots,y_5)\mapsto (y_0',\dots,y_5')$ is 
given by~\eqref{eq:weight_sq_move}.

A special kind of alternating strand diagrams corresponds to Thurston's triple
diagrams~\cite{Thurs}.  Triple diagrams defined below are Thurston's minimal
triple diagrams.  Note that this definition is slightly different from the one
given in~\cite{Thurs}.

\begin{definition}
Consider a disk with $2n$ boundary vertices $b_1,b_1',b_2,b_2',\dots,b_n,b_n'$
(in the clockwise order) on its boundary.
A {\it triple diagram\/} is a diagram with $n$ directed strands drawn inside
the disk such that 
\begin{enumerate}
\item 
Each strand starts a boundary vertex $b_i$ and ends at another vertex $b_j'$.
For each boundary vertex $b_i$ (resp., $b_i'$) there is exactly one strand
starting (resp., ending) at this vertex.
\item Only triple intersections of strands are allowed inside the disk
such that the directions of the six rays at this point alternate
(as we go around the intersection point).
\item Strands have no self-intersections.
\item 
If two strands intersect each other at two points $u$ and $v$,
then one of these strands is directed from $u$ to $v$ and the other
strand is directed from $v$ to~$u$.
\end{enumerate}
The {\it strand permutation\/} $\pi_T\in S_n$ of a triple diagram $T$ is given 
by $\pi_T(i)=j$ whenever $T$ contain a strand from $b_i$ to $b_j'$.
\end{definition}

The map $T\mapsto D$ from triple diagrams to alternating strand diagrams is
quite simple:  Slightly deform each triple crossing of the strands in $T$ and
replace it by 3 simple crossings so that the boundary of newly created triangle
is oriented clockwise, then merge all pairs of boundary vertices $b_i$ and
$b_i'$.

We will call faces of a triple diagram the {\it chambers.}  A triple diagram
has two types of chambers --- the {\it dark chambers\/} whose boundary is
oriented counterclockwise, and the {\it light chambers\/} with clockwise
boundary.  Figure~\ref{fig:triple_diagram} below shows a triple diagram
with dark chambers colored in a darker shade.

\begin{figure}[ht]
\input{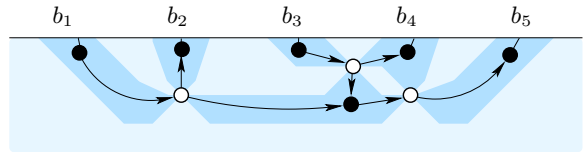}
\caption{The plabic graph $G$ corresponding to the triple diagram $T$
(with a perfect orientation of edges)}
\label{fig:triple_plabic}
\end{figure}

Let us transform a triple diagram $T$ into an alternating strand diagram $D$ and
then into a plabic graph $G$, as above. 
We have the following correspondences:
$$
\begin{array}{l}
\{\text{dark chambers of } T\} \leftrightarrow
\{\text{counterclockwise faces of } D\} \leftrightarrow
\{\text{black vertices of } G\} \\
\{\text{triple crossings in } T\} \leftrightarrow
\{\text{clockwise faces of } D\} \leftrightarrow
\{\text{white vertices of } G\} \\
\{\text{light chambers of } T\} \leftrightarrow
\{\text{alternating faces of } D\} \leftrightarrow
\{\text{faces of } G\}. 
\end{array}
$$
Figure~\ref{fig:triple_plabic} shows the plabic graph $G$ associated
with the triple diagram $T$ on Figure~\ref{fig:triple_diagram}.

The following claim is straightforward.

\begin{lemma}
The above map identifies triple diagrams with alternating strand diagrams 
$D$ such that all clockwise faces of $D$ are triangles and $D$ has no clockwise loops
attached to boundary vertices.
Equivalently, triple diagrams correspond to reduced plabic graphs 
(without non-boundary leaves and unicolored edges) 
such that all white vertices are trivalent and all boundary leaves are black.
\end{lemma}

\begin{remark}
\label{rem:thurston_moves}
Thurston~\cite{Thurs} proved that any two triple diagrams with the same strand
permutation can be obtained from each other by certain moves.  These moves can
be related to the moves of alternating strand diagrams on
Figure~\ref{fig:alternating_moves}.  
\end{remark}

Thurston \cite[Theorem 1]{Thurs} showed that, for each permutation $\pi \in
S_n$, there is a triple diagram $T$ with strand permutation $\pi_T = \pi$.  Let
us give another construction of a triple diagram $T$ with a given strand
permutation $\pi$, which is different from the construction in
\cite{Thurs}.  Note that both constructions are quite simple.

We will arrange the vertices $b_1,b_1',b_2,b_2',\dots,b_n,b_n'$ 
on the $x$-axis in the $xy$-coordinate plane.
We will draw a triple diagram in the half-space below the $x$-axis. 
Each strand will be a continuous curve $(x(t),y(t))$, $t\in[0,1]$
such that $x(t)$ is monotonically increasing or monotonically decreasing
function.  Let us call such special triple diagrams {\it monotone.}

For a strand $S$ from $b_i$ to $b_j'$,
we say that $S$ is {\it rightward\/} if $i\leq j$ 
(because it is directed from left to right),
and we say that $S$ is {\it leftward\/} if $i>j$.
Let $S_i$ denotes the strand starting at $b_i$ and 
$S_j'$ denotes the strand ending at $b_j'$ (so each strand
has two labels).

We will draw all strands in $T$ by adding little pieces to them
as we go from left to right.  Start with the vertices $b_1$ and $b_1'$.
If $\pi(1)=1$, then draw a (very short) strand from $b_1$ to $b_1'$.
Otherwise draw initial segments of 
the strands $S_1$ and $S_1'$ attached to the vertices $b_1$ and $b_1'$,
so that we get two strands with loose ends.  
Then proceed to the pair of vertices $b_2$ and $b_2'$.
If $\pi(2)=2$ then draw a short strand from $b_2$ to $b_2'$.
If $\pi(2)=1$ then attach the loose end of $S_1'$ to $b_2$.
If $\pi(1)=2$ then draw the initial segment of $S_2$, add a triple
crossing of the strands $S_1,S_1', S_2$ and attach the loose end of $S_1$ to $b_2'$.
Otherwise draw two initial segments of $S_2$ and $S_2'$, so that we get $4$ loose
ends of strands.  Then proceed to the pair of vertices $b_3$, $b_3'$, etc.
Note that a each moment we have some number of loose strands with alternating
directions (right, left, right, left, \dots) as we list them from the bottom.
So between an adjacent pair of rightward loose strands there is exactly one leftward loose strand,
and vise versa.
This means we can always switch two adjacent rightward (resp., leftward) loose strands by adding 
a triple crossing.   
Suppose now that we process the pair of vertices $b_i, b_i'$. 
If $\pi(i)=i$, we just draw a short strand $S_i=S_i'$.
If $\pi(i) = j<i$, then we extend the (already drawn) loose end of $S_j'$ all the way up by these 
``adjacent transpositions'' and connect it with $b_i$.  
Similarly, if $\pi^{-1}(i) = j<i$, then we extend the loose end of $S_j$ all the way up
and connect it with $b_i'$.  Otherwise we just add two new strands $S_i$ and $S_i'$
with loose ends.
When we finish processing all boundary vertices, all loose ends
should be attached to the corresponding vertices, 
and we obtain a triple diagram.
Notice that any pair of rightward strands 
(or a pair of leftwards strands) will not intersect
more than once.   (A leftward and a rightward strands can intersect many times but this
is not prohibited.)
Figure~\ref{fig:triple_diagram} shows a monotone triple diagram obtained by this procedure
for the permutation $\pi = 24513$.

\begin{figure}[ht]
\input{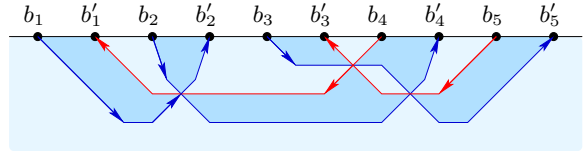}
\caption{A triple diagram $T$ with strand permutation $\pi_T =24513$}
\label{fig:triple_diagram}
\end{figure}

\begin{lemma}
\label{prop:exist_monotone}
For each permutation $\pi\in S_n$ there is a monotone triple 
diagram $T$ with strand permutation $\pi_T=\pi$.
\end{lemma}

\begin{corollary}
For any decorated permutation $\pi^:$ there is a reduced plabic graph $G$ with
trip permutation $\pi_G^: = \pi$.
\end{corollary}

\begin{remark}
Note that monotone triple diagrams are quite similar to 
{\it double wiring diagrams\/} of Fomin-Zelevinsky~\cite{FZ1}.
Such triple diagrams are obtained by superimposing two usual wiring diagrams
formed by all rightwards strands and by all leftward strands.
These two wiring diagrams are shown in blue and red colors on Figure~\ref{fig:triple_diagram}.
Actually, double wiring diagrams from~\cite{FZ1} are exactly 
the monotone triple diagrams in the case when $n=2m$ and $\pi:[m]\to [m+1,2m]$,
$\pi:[m+1,2m]\to [m]$.
\end{remark}

We will need the following property of monotone triple diagrams constructed
above.
For a permutation $\pi\in S_n$, we say that $i\in[n]$ is  
an {\it anti-exceedance\/} of $\pi$ if $\pi^{-1}(i)>i$.
Let $I(\pi)\subset[n]$ be the set of all anti-exceedances of $\pi$.

Recall Proposition~\ref{prop:plabic_matroid} that combinatorially 
describes the matroid corresponding to a perfectly orientable plabic graph $G$
as the matroid $\M_G$ of source sets of perfect orientations of $G$.

\begin{lemma}
\label{lem:monotone_lex_min}
Let $T$ be a monotone triple diagram with strand permutation $\pi$,
and let $G$ be the associated reduced plabic graph.
Then $G$ is perfectly orientable and $I=I(\pi)$ is the lexicographically 
minimal base of the matroid $\M_G$.

Any other base of $\M_G$ is obtained from $I$ by replacing 
some entries $i_1,\dots,i_s\in I$ with some $j_1,\dots,j_s$ such that
$j_1>i_1$, \dots, $j_s>i_s$.
\end{lemma}

\begin{proof}
Define the {\it nose\/} of a dark chamber $C$ in $T$ as its rightmost point,
i.e., the point with the maximal $x$-coordinate.
The monotonicity of a triple diagram implies that each dark
chamber has a unique nose.  Notice that $i\in I(\pi)$ 
if and only if $b_i'$ is not the nose of the dark chamber adjacent to 
the boundary segment $[b_i,b_i']$.
In the triple diagram shown of Figure~\ref{fig:triple_diagram},
the elements of $I(\pi)=\{1,3\}$ correspond to the targets of 
the leftward strands (shown in red).

Let us direct edges in of the plabic graph $G$, as follows. 
For each black vertex $u$ of $G$ (which corresponds to a dark chamber $C$ in $T$)
the only outgoing edge from $u$ goes to the white vertex at the nose
of $C$.  Each white vertex $v$ of $G$ (triple crossing in $T$)
is adjacent to three dark chambers $C_1,C_2,C_3$ such that
$v$ is the nose of $C_1$.  Then the $3$ edges incident to $v$ 
are directed away from $C_1$ and towards $C_2$ and $C_3$;
see Figure~\ref{fig:triple_plabic}.
That means that this orientation of edges is perfect,
which proves perfect orientability of $G$.

Notice that, for this orientation of edges, the boundary source set $I$
is exactly the index set $I(\pi)$ of the ``noseless'' boundary vertices $b_i'$.

Any other perfect orientation of $G$ is obtained from the constructed one by
switching edge directions in a family of noncrossing directed paths 
joining pairs of boundary vertices; see Lemma~\ref{lem:path_matchings}.
A directed path $P:b_i\to b_j$ in this digraph $G$
correspond to a sequence of dark chambers $C_1,\dots,C_l$ in $T$ such that
$C_1$ and $C_l$ are adjacent to the boundary segments $[b_i,b_i']$ and
$[b_j,b_j']$, and $C_{i+1}$ is adjacent to the nose of $C_i$, for $i\in[l-1]$.
For each $i$, the nose $C_{i+1}$ is located strictly to the right of the nose of
$C_i$.  That means that $j>i$ for any directed path $P:b_i\to b_j$.

In other words, when we switch to any other perfect orientation of $G$,
we replace some sources $i_1,\dots,i_s$ by other sources $j_1,\dots,j_s$
such that $j_1>i_1$, \dots, $j_s>i_s$. 
This implies that $I=I(G)$
is the lexicographically minimal source set of a perfect orientation of $G$,
i.e., the lexicographically minimal base of $\M_G$, as needed.
\end{proof}

\begin{corollary}
\label{cor:reduced_orientable}
Any reduced plabic graph is perfectly orientable.
\end{corollary}

\begin{proof}  
Let $G$ be a reduced plabic graph.  We may assume that
it has no boundary leaves, so that the trip permutation
$\pi_G$ has no fixed points.
Let $G'$ be a reduced plabic graph coming from a monotone triple diagram as
above such that $G'$ has the same trip permutation $\pi_{G'}=\pi_G$.
By Lemma~\ref{lem:monotone_lex_min}, $G'$ is reduced and by 
Theorem~\ref{th:moves=permutations} $G'$ is move-equivalent to $G$.
Since the moves never change perfect orientability, the graph $G$ is 
also perfectly orientable.
\end{proof}

\begin{remark}
Note that not all plabic graphs are perfectly orientable.
For example, if a a graph has a singleton (isolated component 
with a single vertex), then it is not perfectly orientable.
Essentially this is the only obstruction for perfect orientability.
Indeed, any graph can be transformed into a reduced graph 
possibly together with some singletons; see Lemma~\ref{lem:any_plabic_reduced}.
Such graph is perfectly orientable if and only if it has no
singletons.
\end{remark}


\section{Mutations of dual graphs}

\section{From matroids to decorated permutations}
\label{sec:matroids_perms}

\begin{definition}
A {\it Grassmann necklace\/} is a sequence $\I=(I_1,\dots,I_n)$ of subsets $I_r\subseteq[n]$ such
that, for $i\in[n]$, 
if $i\in I_i$ then $I_{i+1} = (I_i\setminus\{i\})\cup \{j\}$, for some $j\in[n]$;
and if $i\not\in I_r$ then $I_{i+1} = I_i$.
(Here indices $i$ are taken modulo $n$.)
In particular, we have $|I_1|=\cdots = |I_n|$.
\end{definition}

Such necklaces are in bijection with decorated permutations.  For a Grassmann
necklace $\I$, we construct the decorated permutation $\pi^:(\I)=(\pi,\col)$
such that
\begin{enumerate}
\item
if $I_{i+1}=(I_i \setminus\{i\})\cup\{j\}$, $j\ne i$, then $\pi(i)=j$;
\item
if $I_{i+1}=I_i$ and $i\not\in I_i$ then $\pi(i)=i$ is colored in black
$\col(i)=1$;
\item
if $I_{i+1}=I_i$ and $i\in I_i$ then $\pi(i)=i$ is colored in white
$\col(i)=-1$;
\end{enumerate}
where the indices $i$ are taken modulo $n$.
Notice that as we transform $I_1$ to $I_2$, then to $I_3$, and so on until get
get back to $I_1$, we can remove an element $i$ at most once (at the $i$-th
step);  and thus we can add an element $j$ at most once.  This shows that $\pi$
is indeed a permutation in $S_n$.  Note that black fixed points of $\pi^:$ are
exactly the elements $i\in[n]$ that belong to none of the subsets
$I_1,\dots,I_n$ and white fixed points are exactly the elements $j\in[n]$ that
belong to all subsets $I_1,\dots,I_n$.

Let us describe the inverse map from decorated permutations to Grassmann necklaces.
For a decorated permutation $\pi^:=(\pi,\col)$ of size $n$, we say that $i\in
[n]$ is an {\it anti-exceedance\/} of $\pi^:$ if $\pi^{-1}(i)>i$ or $\pi(i) =i$ and
$\col(i) = -1$. (That is we regard white fixed points as anti-exceedances.)
Let $I(\pi^:)\subset [n]$ be the set of all anti-exceedances in $\pi^:$.
For $r\in[n]$, let us also define the {\it shifted anti-exceedance set\/} $I_r(\pi^:)$
as the set of indices $i\in[n]$ such that $i<_r \pi^{-1}(i)$ or ($\pi(i) =i$ and
$\col(i) = -1$), where $<_r$ is the cyclical shift of the usual linear order on $[n]$:
$r<_r (r+1) <_r \cdots <_r n <_r < 1<_r \cdots <_r (r-1)$, i.e.,
$I_r(\pi^:)$ is the anti-exceedance set with respect to the linear order $<_r$.
In particular, $I_1(\pi^{:}) = I(\pi^:)$.
Let $\I(\pi^:) = (I_1,\dots,I_n)$, where $I_r = I_r(\pi^:)$, for $r\in[n]$.

\begin{lemma}
\label{lem:neck_perms}
The maps $\pi^:\mapsto \I(\pi^:)$ and $\I\mapsto \pi^:(\I)$ are inverse to each
other.  They give a bijection between decorated permutations $\pi^:$ of size
$n$ and Grassmann necklaces $\I$ of size $n$.  \end{lemma}

The proof of this lemma is quite straightforward
(an exercise for the reader).

We can graphically present decorated permutation $\pi^:$ by 
arranging the vertices $b_1,\dots,b_n$ clockwise on a circle,
drawing straight directed chords $[b_i,b_{\pi(i)}]$, whenever $\pi(i)\ne i$;
drawing a counterclockwise loop at $b_j$ for each black fixed point $\col(j)=1$;
and drawing a clockwise loop at $b_l$ for each white fixed point $\col(l)=-1$.

\begin{figure}[ht]
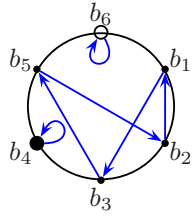

\psset{unit=.7pt}
\pspicture(-60,-40)(60,50)
\pscircle[linecolor=black](0,0){40}
\cnode*[linewidth=0,linecolor=black](34.64,20){2}{B1}
\cnode*[linewidth=0,linecolor=black](34.64,-20){2}{B2}
\cnode*[linewidth=0,linecolor=black](0,-40){2}{B3}
\cnode*[linewidth=0,linecolor=black](-34.64,-20){4}{B4}
\cnode*[linewidth=0,linecolor=black](-34.64,20){2}{B5}
\cnode[linewidth=1,fillcolor=white,linecolor=black](0,40){4}{B6}
\rput(43.30,25){$b_1$}
\rput(43.30,-25){$b_2$}
\rput(0,-50){$b_3$}
\rput(-43.30,-25){$b_4$}
\rput(-43.30,25){$b_5$}
\rput(0,50){$b_6$}

\ncline{->}{B1}{B3}
\ncline{->}{B2}{B1}
\ncline{->}{B3}{B5}
\ncline{->}{B5}{B2}
\nccurve[angleA=-60,angleB=-120,ncurv=10]{->}{B6}{B6}
\nccurve[angleA=0,angleB=60,ncurv=10]{->}{B4}{B4}

\endpspicture
\psset{unit=1pt}
\caption{A decorated permutation $\pi^:$}
\label{fig:dec_permut}
\end{figure}

For example, Figure~\ref{fig:dec_permut} shows
the decorated permutation $\pi^:$ with 
$\pi = 3\,1\,5\, \underline{4}\,2\,\underline{6}$ with
two fixed points $4$ and $6$ colored $\col(4)=1$ (black) and $\col(6)=-1$ 
(white).
This decorated permutation has the following shifted
anti-exceedance sets
$I_1 = \{1,2,6\}$, $I_2 = \{2, 3, 6\}$,
$I_3 = \{3, 6, 1\}$, $I_4 = \{5, 6, 1\}$,
$I_5 = \{5, 6, 1\}$, $I_6 = \{6, 1, 2\}$.


For a matroid $\M\subseteq\binom{[n]}{k}$ of rank $k$ on the set $[n]$, 
let $\I_\M=(I_1,\dots,I_n)$ 
be the sequence of subsets in $[n]$ such that, for $i\in[n]$, $I_i$ is the lexicographically minimal base of $\M$ 
with respect to the shifted linear order $<_i$ on $[n]$.

\begin{lemma}
\label{lem:IM_necklace}
For a matroid $\M$, the sequence $\I_\M$ is a Grassmann necklace.
\end{lemma}

\begin{proof}
Let $\I(\M)=(I_1,\dots,I_n)$.
By the cyclic symmetry, it is enough to check that 
$I_2 = (I_1\setminus\{1\})\cup\{j\}$ for some $j$, or $I_2 = I_1$.
The subset $I_1=\{i_1<\dots<i_k\}$ is the lex minimal base of $\M$ 
(with respect to the usual order on $[n]$).  
If $i_1\ne 1$ then  $I_1$ is also the lex minimal base 
with respect to the order $<_2$, and thus $I_2=I_1$.
Assume that $i_1 = 1$ and $I_2=\{j_1<\dots<j_k\}\ne I_1$.
Let $r$ be the index such that $j_s=i_{s+1}$ for all $s<r$
and $j_r\not\in I_1$.  Then $j_r\in [i_{r}+1,i_{r+1}-1]$
(or $j_r\in [i_{r}+1,n]$ if $r=k$) and 
$\M$ has a circuit (minimal dependence) involving $i_1, j_r$
and some of the $j_1,\dots,j_{r-1}$.
That implies that, for any $(k-r)$-element subset $S\subset [j_r+1,n]$,
$\{i_1,\dots,i_r\}\cup S$ is a base of $\M$ if and only if
$\{j_1,\dots,j_r\}\cup S$ is a base.
Because of the lex minimality of $I_1$ and $I_2$, we have $i_t=j_t$
for all $t>r$.
Thus $I_2 = (I_1\setminus\{1\})\cup\{j_r\}$, as needed.
\end{proof}


Recall that, for a plabic graph $G$, the image of the boundary measurement map
$\tilde \Mes_G$ belongs to a totally nonnegative Grassmann cell $S_\M^\tnn$
where $\M=\M_G$; see Proposition~\ref{prop:plabic_matroid}.

\begin{proposition}
\label{prop:matorid_to_perm}
Let $G$ be a reduced plabic graph, and let
$\M = \M_G$ be the associated matroid.  Then 
the Grassmann necklace $\I_\M$ of the matroid $\M$ equals the
necklace $\I(\pi_G^:)$ corresponding to the decorated trip permutation of $G$.
\end{proposition}

\begin{proof}
Black boundary leaves of $G$ correspond to isolated boundary sinks 
in directed network, which correspond to zeros of the matroid $\M$,
i.e., the elements $i\in[n]$ which do not appear in any base of $\M$.
These elements never appear in the necklace $\I_\M$.
They give black fixed points of the decorated permutation
$\pi^:(\I_\M)$, as needed.  
Similarly, white boundary leaves of $G$ give cozeros of $\M$,
i.e., the element $i\in[n]$ that appear in all bases of $\M$.
They give white fixed points in $\pi^:(\I_\M)$.

Thus we can remove all boundary leaves and assume that the 
reduced plabic graph $G$ has no boundary leaves.  Then the trip
permutation $\pi_G$ has no fixed points.

According to the cyclic symmetry (Remark~\ref{rem:cyclic_symmetry}) it is
enough to show that the anti-exceedance set $I(\pi_G)$ equals the 
lexicographically minimal base $I$ of $\M$. 
Let $G'$ be the plabic graph with the same trip permutation
obtained from a monotone triple diagram as shown 
in Section~\ref{sec:alt_strand}.
Then by Theorem~\ref{th:moves=permutations}
the graphs $G$ and $G'$ can be obtained from each other by
moves (M1)--(M3).   These moves never change the trip permutation $\pi_G$ and
never change the matroid $\M_G$ (because they do not change the image
$\Image(\tilde M_G)\subseteq S_\M^\tnn$).
Thus the needed claim follows from Lemma~\ref{lem:monotone_lex_min}.
\end{proof}

Now we can finally prove Theorems~\ref{th:moves=boundary}
and~\ref{th:reduced_plabic_parametrization}.

\begin{proof}[Proof of Theorems~\ref{th:moves=boundary}
and~\ref{th:reduced_plabic_parametrization}]
Let $N$ and $N'$ be two perfectly orientable plabic networks
with graphs $G$ and $G'$ such that $\tilde\Mes(N)=\tilde\Mes(N')$.
According to Lemma~\ref{lem:any_plabic_reduced}, we can 
transform these networks by the moves and reductions into networks
with reduced graphs and maybe some singleton components.
If there are singleton components then the graph(s) are
not perfectly orientable.
Thus we may assume that the plabic graphs $G$ and $G'$ are reduced.

The boundary measurement map sends $N$ and $N'$ into into the same cell
$S_\M^\tnn$.  Thus, by Proposition~\ref{prop:matorid_to_perm}, the graphs $G$
and $G'$ have the same necklaces $\I(\pi^:_G) = \I(\pi^:_G) = \I_\M$.  Thus by
Lemma~\ref{lem:neck_perms}, the decorated trip permutations $\pi^:_G$ and
$\pi^:_G$ are the same.  According to Theorem~\ref{th:moves=permutations}, the
graphs $G$ and $G'$ are move-equivalent.

We know that for any cell $S_\M^\tnn$ there is a plabic graph $G''$
such that $\tilde \Mes_{G''}$ is a subtraction-free parametrization of $S_\M^\tnn$.
Indeed, we can take the \Le-diagram associated with $S_\M^\tnn$
(see Theorem~\ref{th:g_D}) and transform it into a plabic graph.
The graph $G''$ must be reduced.  (Otherwise, we can kill
some parameter and the map $\tilde \Mes_{G''}$ would not be 
a parametrization; see Remark~\ref{rem:non-reduced}.)
Thus again $G''$ has the same decorated trip permutation $\pi^:_{G''} = \pi^:_G$
and is move-equivalent to $G$ and $G'$.

Note that every time when we perform moves (M1)--(M3), the face variables
$y_f$ are transformed by invertible subtraction-free rational maps.
Thus for the graph $G$ (and any other graph obtained from $G''$ by the moves), 
the map $\tilde \Mes_{G}$ is obtained from $\tilde \Mes_{G''}$ by a sequence
of these reparametrization maps, and thus $\tilde \Mes_{G}$ is also
a subtraction-free rational parametrization of $S_\M^\tnn$.

So any two networks with the same graph $G$ that maps into the same point in 
the Grassmannian must be equal to each other.
That means that if we transform the network $N'$ by the moves
into a network with the graph $G$ we will get the network $N$.
Thus the networks $N$ and $N'$ are move-equivalent.
\end{proof}

\begin{corollary}
\label{cor:Im=Cell}
For any perfectly orientable plabic graph $G$ (not necessarily reduced)
that corresponds to the cell $S_\M^\tnn$, we have
$$
\tilde \Mes_G(\R_{>0}^{F(G)-1}) = S_\M^\tnn,
$$
that is the image of $\tilde \Mes_G$ is the whole cell $S_\M^\tnn$.
\end{corollary}

\begin{proof} 
For a reduced graph $G$ this follow 
from Theorem~\ref{th:reduced_plabic_parametrization}.
Other graphs can be transformed into reduced ones
by the moves and reductions (Lemma~\ref{lem:any_plabic_reduced}),
but the moves and reductions do not change the image of the map
$\tilde \Mes_G$.
\end{proof}

\section{Circular Bruhat order}
\label{sec:circ_bruhat}

In this section we show that each nonnegative Grassmann cell
$S_\M^\tnn$ is an intersection of $n$ Schubert cells.
Then we combinatorially describe the partial
order on the cells $S_\M^\tnn$ by containment of their
closures.

Let us use notation of Section~\ref{sec:matroids_perms}.
Let us say that a decorated permutation $\pi^:$ has {\it type\/}
$(k,n)$ if $\pi^:$ has size $n$ and it has $k$
anti-exceedances.
Also say that a Grassmann necklace has {\it type\/} $(k,n)$
if it consists of $k$-element subsets in $[n]$.
Clearly, the types of corresponding decorated permutations
and Grassmann necklaces are the same.

\begin{theorem}
\label{th:SM_dec_perm}
The map $S_\M^\tnn \mapsto \pi^:(\I_\M)$ is a bijection between 
nonnegative Grassmann cells $S_\M^\tnn\subset Gr_{kn}^\tnn$
and decorated permutations of type $(k,n)$.
Equivalently, the map $S_\M^\tnn \mapsto \I_\M$ is a
bijection between nonnegative Grassmann cells in $Gr_{kn}^\tnn$
and Grassmann necklaces of type $(k,n)$.
\end{theorem}

\begin{proof}
By Theorems~\ref{th:reduced_plabic_parametrization}
and~\ref{th:moves=permutations},
all reduced plabic graphs with the same decorated trip permutation $\pi^:$
correspond to the same nonnegative Grassmann cell $S_\M^\tnn$.
By Propositions~\ref{prop:plabic_matroid} and~\ref{prop:matorid_to_perm}, 
we have $\pi^:=\pi^:(\I_\M)$.
Thus, two reduced plabic graphs with different trip permutations correspond 
to different Grassmann cell.
\end{proof}

Recall that $\Omega_\lambda^w: = w(\Omega_\lambda)$,
for $w\in S_n$, are the permuted Schubert cells;
see Section~\ref{sec:grassmannian}.  
Let us also use the subset notation for these Schubert cells
$$
\Omega^w_{I}:=\{V\in Gr_{kn}\mid I \text{ is the lex minimal base
of $\M_V$ with respect to $<_w$}\},
$$
where $<_w$ is the linear order on $[n]$ given $w(1)<w(2)<\cdots < w(n)$,
cf.~Section~\ref{ssec:Matroid_strata}.
The cells $\Omega^w_I$ are exactly the cells 
$\Omega_\lambda^w$ labelled by subsets rather than partitions.
Each matroid strata
$S_\M$ is an intersection of several permuted Schubert cells;
see Remark~\ref{rec:common_ref_factorial}.
For a nonnegative Grassmann cell $S_\M^\tnn$ only $n$ Schubert cells 
are needed.
Let $c=(1,\dots,n)\in S_n$ be the long cycle.

\begin{theorem}
\label{th:SM=circ_necklace}
Let $S_\M^\tnn\subset Gr_{kn}^\tnn$ be a nonnegative Grassmann cell,
and let $\I_\M = (I_1,\dots,I_n)$
be the Grassmann necklace corresponding to $\M$.
Then 
$$
S_\M^\tnn = \bigcap_{i=1}^n \Omega_{I_i}^{c^{i-1}}\cap Gr_{kn}^\tnn.
$$
Moreover, for any arbitrary collection of $k$-subsets 
$I_1,\dots,I_n\subset[n]$, the intersection in the right-hand-side
is nonempty if and only if $(I_1,\dots,I_n)$ is a Grassmann necklace.
\end{theorem}

\begin{proof}
Note that $<_{c^{i-1}}$ is exactly the shifted order $<_i$ on $[n]$.
By the definition of $\I_\M$, 
the cell $S_\M^\tnn$ belongs to the intersection of Schubert
cells in the right-hand-side.
Suppose that this intersection contains an element of another
cell $S_{\M'}^\tnn$.  Then $\I_{\M'} = \I_\M$,
so the cell $S_{\M'}^\tnn$ corresponds
to the same decorated trip permutation 
$\pi^:(\I_\M) = \pi^:(\I_{\M'})$, which is impossible by 
Theorem~\ref{th:SM_dec_perm}.
The second claim follows from Lemma~\ref{lem:IM_necklace}.
\end{proof}

Let $\overline{S_\M^\tnn} \subseteq Gr_{kn}^\tnn$ denotes the closure of the cell 
$S_\M^\tnn\subset Gr_{kn}^\tnn$. 
Define the partial order $\leq$ on nonnegative Grassmann cells
such that $S_\M^\tnn \leq S_{\M'}^\tnn$ if and only if 
$\overline{S_\M^\tnn}\subseteq \overline{S_{\M'}^\tnn}$.
Recall that the cell $S_\M^\tnn$ is given in the Pl\"ucker coordinates
by the conditions $\Delta_I>0$ for $I\in \M$, 
and $\Delta_J = 0$ for $J\not\in \M$; 
see Section~\ref{set:tot_nonneg_Grass}.
Thus, for two nonempty cells, we have 
$S_\M^\tnn \leq S_{\M'}^\tnn$ if and only if $\M\subseteq \M'$.

Recall that the partial order on the Schubert cells by the containment 
has a simple combinatorial description
$\overline{\Omega}_{\lambda}\subseteq \overline{\Omega}_{\mu}$
if and only if $\lambda\subseteq\mu$ (meaning the inclusion of Young
diagrams); see~\cite{Fult}.


Recall the standard bijection $\lambda\mapsto I(\lambda)$ 
between partitions $\lambda\subseteq(n-k)^k$ and  $k$-subsets $I\subset [n]$;
see Section~\ref{sec:grassmannian}.  Let $I\mapsto \lambda(I)$
be the inverse bijection.

\begin{lemma}
\label{lem:lambdaJI}
Let $S_\M^\tnn$ be a nonnegative Grassmann cell.
Let $I$ be the lexicographically minimal base of $\M$.
Then, for any other base $J\in \M$, we have $\lambda(J)\subseteq\lambda(I)$.
\end{lemma}

\begin{proof}
According to Theorem~\ref{th:SM_dec_perm},
Lemma~\ref{prop:exist_monotone}, and the second part of
Lemma~\ref{lem:monotone_lex_min}, 
any base $J$ of $\M$ is obtained from $I$ by switching
some entries $i_1,\dots,i_s\in I$ with some $j_1,\dots,j_s$
such that $j_1>i_1$, \dots, $j_s>i_s$. 
But this exactly means that for the corresponding partitions 
we have $\lambda(J)\subseteq \lambda(I)$.
\end{proof}

Let $\I_\M = (I_1,\dots,I_n)$ be the Grassmann necklace corresponding to a
matroid $\M\subseteq \binom{[n]}{k}$.
Let us also denote 
$\Lambda_\M = (\lambda_{(1)},\dots,\lambda_{(n)})$
the sequence of partitions $\lambda_{(i)} = \lambda(c^{-i+1}(I_i))$,
for $i\in[n]$.
In other words, the partitions $\lambda_{(i)}$ are chosen so that
$\Omega_{I_i}^{c^{i-1}} = \Omega_{\lambda_{(i)}}^{c^{i-1}}$.

\begin{proposition}
\label{prop:partial_order_lambda}
Let $S_\M^\tnn,S_{\M'}^\tnn$ be two cells in $Gr_{kn}^\tnn$,
$\Lambda_\M = (\lambda_{(1)},\dots,\lambda_{(n)})$
and $\Lambda_{\M'} = (\lambda_{(1)}',\dots,\lambda_{(n)}')$.
Then $S_\M^\tnn\leq S_{\M'}^\tnn$ if and only if 
$\lambda_{(i)}\subseteq \lambda_{(i)}'$,
for all $i\in[n]$.
\end{proposition}

\begin{proof}
If $\lambda_{(i)}\subseteq \lambda_{(i)}'$ for all $i\in[n]$,
then $S_\M^\tnn\leq S_{\M'}^\tnn$,
by Theorem~\ref{th:SM=circ_necklace}.
On the other hand, suppose that $S_\M^\tnn\leq S_{\M'}^\tnn$.
Then $\M\subseteq \M'$.  In particular, the lexicographically 
minimal base $I$ of of $\M$ is also a base of $\M'$.
By Lemma~\ref{lem:lambdaJI}, this implies
that $\lambda_{(1)}\subseteq \lambda_{(1)}'$.
Taking cyclic shifts and using the cyclic symmetry of 
the construction, we get $\lambda_{(i)}\subseteq\lambda_{(i)}'$,
for any $i\in[n]$.
\end{proof}

Let us now describe the partial order on the cells $S_\M^\tnn$ in terms of 
decorated permutations.  

\begin{definition}  
The {\it circular Bruhat order\/} $CB_{kn}$ is the partial order $\leq$ on
the set of decorated permutations of type $(k,n)$
such that, for two decorated permutations $\pi^:$ and $\sigma^:$
corresponding to the cells $S_{\M}^\tnn$ and $S_{\M'}^\tnn$ in $Gr_{kn}^\tnn$, 
we have $\pi^:\leq \sigma^:$ if and only if 
$S_\M^\tnn \leq S_{\M'}^\tnn$. 
\end{definition}

\begin{lemma}
\label{lem:CB_top_bottom}
The circular Bruhat order $CB_{kn}$ has a unique top element
given by the decorated permutation $\pi_{top}:i\mapsto i+k\pmod n$, for $i\in[n]$
(for $k=0$, all fixed points of $\pi_{top}$ are colored black,
and for $k=n$ all fixed points of $\pi_{top}$ are colored white).
The circular Bruhat order $CB_{kn}$ has $\binom{n}{k}$ minimal elements 
corresponding to the identity permutation permutation with some $k$ 
fixed points colored in white and remaining $(n-k)$ fixed points colored in black.
\end{lemma}

\begin{proof}
The top element of $CB_{kn}$ corresponds to the top cell
$S_\M^\tnn\subset Gr_{kn}^\tnn$, where $\M=\binom{[n]}{k}$ is 
the complete matroid of rank $k$ on $[n]$.  
By Lemmas~\ref{lem:neck_perms} and~\ref{lem:IM_necklace},
this matroid corresponds to the Grassmann necklace 
$\I_\M=(I_1,\dots,I_n)$ with $I_i=\{i,i+1,\dots,i+k\}$, for 
$i\in[n]$ (elements are taken modulo $n$);  and this 
necklace corresponds to decorated permutation $\pi_{top}:i\mapsto i+k$.

On the other hand, minimal elements of $CB_{kn}$ correspond to 0-dimensional
cells $S_\M^\tnn$.  These cells are fixed points of the torus action on 
the Grassmannian $Gr_{kn}$.  In other words, they correspond to matroids
with a single base $\M=\{I\}$.  Under the correspondence $\M\mapsto \pi^:$,
the $k$ elements of $I$ give $k$ white fixed points of $\pi^:$
and $n-k$ elements of $[n]\setminus I$ give $n-k$ black fixed points of $\pi^:$.
\end{proof}

For $a,b\in[n]$, define the {\it cyclic interval\/} $[a,b]^{cyc}$ as 
$\{a,a+1,\dots,b\}$
if $a\leq b$, and as $\{a,a+1,\dots,n,1,\dots,b-1\}$ if $a>b$.
In other words, a cyclic interval is a sequence of consecutive numbers
arranged on a circle in the clockwise order.
For a decorated permutation $\pi^:$ and a pair $a,b\in [n]$,
let us define the number $r_{ab}(\pi^:)$ as the number
of shifted anti-exceedances $i$ of $\pi^:$ with respect
to the shifted order $<_a$ such that $i\in [a,b]^{cyc}$.
In other words, if $\I(\pi^:)=(I_1,\dots,I_n)$,
then $r_{ab}(\pi^:) = |I_a\cap \{a,a+1,\dots,b\}|$.
In particular, for a decorated permutation of type $(k,n)$,
we have $r_{a,a-1}(\pi^:)=k$, for any $a\in[n]$.
(Here we take indices $a,b$ modulo $n$.)

\begin{corollary}
\label{cor:rab_pi}  
For two decorated permutations 
$\pi^:$ and $\sigma^:$ of the same type,
we have $\pi^:\leq \sigma^:$ if and only if
$r_{ab}(\pi^:)\leq r_{ab}(\sigma^:)$ for all $a,b\in[n]$.
\end{corollary}

\begin{proof}
Follows from Proposition~\ref{prop:partial_order_lambda} and the fact
that, for two partitions $\lambda,\mu\subseteq(n-k)^k$, we have
$\lambda\subseteq \mu$ if and only if 
$|I(\lambda)\cap[b]|\leq |I(\mu)\cap[b]|$ for any $b\in[n]$.
\end{proof}

Recall, that we presented each decorated permutation $\pi^:$ by
a chord diagram; see Figure~\ref{fig:dec_permut}.
Also recall {\it crossings, alignments,} and {\it misalignments\/}
from Section~\ref{sec:loop_erased}.
This notions can be adapted for
decorated permutations, as follows.

\begin{figure}[ht]
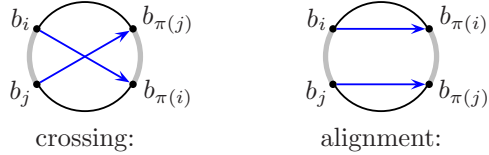

\psset{unit=.7pt}
\pspicture(-80,-40)(80,30)
\rput(0,-45){crossing:}
\pscircle[linecolor=black](0,0){30}
\psarc[linecolor=lightgray,linewidth=2pt](0,0){30}{150}{-150}
\psarc[linecolor=lightgray,linewidth=2pt](0,0){30}{-30}{30}
\cnode*[linewidth=0,linecolor=black](25.98,15){2}{B1}
\cnode*[linewidth=0,linecolor=black](25.98,-15){2}{B2}
\cnode*[linewidth=0,linecolor=black](-25.98,-15){2}{B4}
\cnode*[linewidth=0,linecolor=black](-25.98,15){2}{B5}
\ncline{->}{B5}{B2}
\ncline{->}{B4}{B1}
\rput(-34.64,20){$b_i$}
\rput(45,20){$b_{\pi(j)}$}
\rput(-34.64,-20){$b_j$}
\rput(45,-20){$b_{\pi(i)}$}
\endpspicture
\pspicture(-80,-40)(80,30)
\pscircle[linecolor=black](0,0){30}
\rput(0,-45){alignment:}
\psarc[linecolor=lightgray,linewidth=2pt](0,0){30}{150}{-150}
\psarc[linecolor=lightgray,linewidth=2pt](0,0){30}{-30}{30}
\cnode*[linewidth=0,linecolor=black](25.98,15){2}{B1}
\cnode*[linewidth=0,linecolor=black](25.98,-15){2}{B2}
\cnode*[linewidth=0,linecolor=black](-25.98,-15){2}{B4}
\cnode*[linewidth=0,linecolor=black](-25.98,15){2}{B5}
\ncline{->}{B5}{B1}
\ncline{->}{B4}{B2}
\rput(-34.64,20){$b_i$}
\rput(45,20){$b_{\pi(i)}$}
\rput(-34.64,-20){$b_j$}
\rput(45,-20){$b_{\pi(j)}$}
\endpspicture
\psset{unit=1pt}
\caption{A crossing and an alignment}
\label{fig:cross_align}
\end{figure}

Let $\pi^:$ be a decorated permutation, and let 
$(b_i,b_{\pi(i)})$ and $(b_j,b_{\pi(j)})$, $i\ne j$,
be a pair of chords (or loops).
We say that this pair is a {\it crossing\/} if
$\pi(j)\in [i,\pi(i)]^{cyc}$ and $j\in [\pi(i),i]^{cyc}$.
We say that this pair is an {\it alignment\/} if 
$\pi(i)\in [i, \pi(j)]^{cyc}$ and $j\in [\pi(j), i]^{cyc}$;
if $\pi(i)=i$ then $i$ must be colored $\col(i)=1$ 
(counterclockwise loop),
and if $\pi(j)=j$ then $j$ must be colored $\col(j)=-1$
(clockwise loop); see Figure~\ref{fig:cross_align}.
Note that in a crossing the vertex $b_i$ is allowed to 
coincide with $b_{\pi(j)}$, and the vertex $b_j$ is allowed
to coincide with $b_{\pi(i)}$.  But a loop can never participate
in a crossing.  On the other hand, two chords in an alignment
never have common vertices, but the vertex $b_i$ can coincide
with $b_{\pi(i)}$ and form a counterclockwise loop,
and similarly $(b_j,b_{\pi(j)}$ can form a clockwise loop.
In particular, any counterclockwise loop forms an alignment
with any clockwise loop.

Let us say that a crossing as above is a {\it simple crossing\/}
if there are no any other chords $(l,\pi(l))$ such that 
$l\in[j,i]^{cyc}$ and $\pi(l)\in [\pi(j),\pi(i)]^{cyc}$.
Similarly, a {\it simple alignment\/} is an alignment such that there
are not other chords $(l,\pi(l))$ such that 
$l\in[j,i]^{cyc}$ and $\pi(l)\in [\pi(i),\pi(j)]^{cyc}$.
In other words, simple crossings and alignments should have
no other chords that start on the left between $b_i$ and $b_j$
and end on the right between $b_{\pi(i)}$ and $b_{\pi(j)}$.
Notice if $\pi^:$ has a crossing/alignment then it should
also have a simple crossing/alignment.  (Just pick the one where
$b_i$ and $b_j$ are closest to each other.)

For two decorated permutations $\pi^:, \sigma:$ of the same type,
let $\pi^:\lessdot \sigma^:$ denotes the covering relation in
the circular Bruhat order.

\begin{theorem}
\label{th:circ_bruhat_cover}
In the circular Bruhat order $CB_{kn}$ we have
$\pi^:\lessdot \sigma^:$ if and only if the chord diagram $\pi^:$ is obtained
from the chord diagram of $\sigma^:$ by replacing a simple crossing with the 
corresponding simple alignment as shown on Figure~\ref{fig:cover_rel}.
\end{theorem}

Note again that we allow $i=\sigma(j)$ and/or $j=\sigma(i)$.
In this case $\pi^:$ should have a counterclockwise loop at $b_i=b_{\pi(i)}$
and/or a clockwise loop at $b_j=b_{\pi(j)}$.
In the case when $i=\sigma(j)$ and $j=\sigma(i)$, we can also switch $i$ and 
$j$ and get the decorated permutation clockwise loop at $b_i$ and counterclockwise
loop at $b_j$, which is also covered by $\sigma^:$.

\begin{figure}[ht]
\psset{unit=.7pt}
\pspicture(-50,-60)(170,40)
\pscircle[linecolor=black](0,0){30}
\psarc[linecolor=lightgray,linewidth=2pt](0,0){30}{150}{-150}
\psarc[linecolor=lightgray,linewidth=2pt](0,0){30}{-30}{30}
\rput(-34.64,20){$b_i$}
\rput(45,20){$b_{\pi(i)}$}
\rput(-34.64,-20){$b_j$}
\rput(45,-20){$b_{\pi(j)}$}
\cnode*[linewidth=0,linecolor=black](25.98,15){2}{B1}
\cnode*[linewidth=0,linecolor=black](25.98,-15){2}{B2}
\cnode*[linewidth=0,linecolor=black](-25.98,-15){2}{B4}
\cnode*[linewidth=0,linecolor=black](-25.98,15){2}{B5}
\ncline{->}{B5}{B1}
\ncline{->}{B4}{B2}
\rput(0,-50){$\pi^:$}
\pscircle[linecolor=black](120,0){30}
\psarc[linecolor=lightgray,linewidth=2pt](120,0){30}{150}{-150}
\psarc[linecolor=lightgray,linewidth=2pt](120,0){30}{-30}{30}
\rput(85.36,20){$b_i$}
\rput(165,20){$b_{\sigma(j)}$}
\rput(85.36,-20){$b_j$}
\rput(165,-20){$b_{\sigma(i)}$}
\cnode*[linewidth=0,linecolor=black](145.98,15){2}{B21}
\cnode*[linewidth=0,linecolor=black](145.98,-15){2}{B22}
\cnode*[linewidth=0,linecolor=black](94.02,-15){2}{B24}
\cnode*[linewidth=0,linecolor=black](94.02,15){2}{B25}
\ncline{->}{B25}{B22}
\ncline{->}{B24}{B21}
\rput(120,-50){$\sigma^:$}
\endpspicture
\psset{unit=1pt}
\caption{A covering relation $\pi^:\lessdot \sigma^:$ in the circular Bruhat order}
\label{fig:cover_rel}
\end{figure}

\begin{proof}
One directly verifies that for any pair $\pi^:$ and $\sigma^:$ 
related by ``undoing a crossing'' as above, we have $r_{ab}(\pi^:)
\leq r_{ab}(\sigma^:)$ for any $a,b\in[n]$.  Thus we have $\pi^:< \sigma^:$.
One the other hand, let $\I(\pi^:) = (I_1,\dots,I_n)$ 
and  $\I(\sigma^:) = (I_1',\dots,I_n')$. 
Then the relation $\pi^:\leq \sigma^:$ mean that $I_i$ is obtained
from $I_i'$ by moving some elements ``to the right'' (with respect
to the linear order $<_i$ on $[n]$).
We may assume that $I_1\ne I_1'$.  (Otherwise cyclic shift all elements
using the cyclic symmetry of the construction.)
Suppose that $I_1$ is obtained from $I_1'$ by switching the elements
$i_1,\dots,i_s\in I_1'$ with $j_1,\dots,j_s$, where $j_1>i_1$, \dots, $j_s>i_s$.
Let $a=\pi^{-1}(i_1)$ and $b=\pi^{-1}(j_1)$.
Then $a\leq i_1$ and $b\geq j_1$ because $j_1$ belongs to the anti-exceedance
set of $\pi^:$ and $i_1$ does not belong to this set.
Note that for $i\in [b+1,a]^{cyc}$, the set $I_i$ contains $j_1$ and does
not contain $i_1$.  
Let $(I_1'',\dots,I_n'')$ be the Grassmann necklace 
such that $I_i'' = (I_i \setminus \{j_1\})\cup \{i_1\}$,
for $i\in [b+1,a]^{cyc}$,  and $I''_i = I_i$ otherwise.
Let $\rho^:$ be the decorated permutation corresponding to $(I_1'',\dots,I_n'')$.
Then we have $\pi^:<\rho^:\leq \sigma^:$ and $\pi^:$ is obtained
from $\rho^:$ by undoing a crossing.
\end{proof}

In Section~\ref{sec:gluing} we will explicitly describe how the cells are glued
to each other, which will explain where these ``undoings of crossings'' are
coming from.

Let us define the {\it alignment number\/} $A(\pi^:)$ of a decorated
permutation as the total number of pairs of chords (or loops) 
in $\pi^:$ forming an alignment.
The following claim is obtained by an easy verification.

\begin{lemma}
Suppose that $\pi^:$ is obtained from $\sigma^:$ by 
undoing a simple crossing as shown on 
Figure~\ref{fig:cover_rel}.
Then $A(\pi^:)=A(\sigma^:)+1$.
\end{lemma}


Notice that the maximal element $\pi_{top}$ in $CB_{kn}$ has no alignments.
On the other hand, all minimal elements $\pi^:\in CB_{kn}$ have
$k(n-k)$ alignments; see Lemma~\ref{lem:CB_top_bottom}.
Also notice that the dimension of the top cell in $Gr_{kn}^\tnn$ is 
$k(n-k)$.
As we go down from the top element $\pi_{top}$ to a minimal element in 
$CB_{kn}$ by undoing simple crossings, the alignment number $A(\pi^:)$
increases by 1 at each step.  On the other hand, dimensions of the corresponding
cells should drop by at least 1 at each step.  Since the dimension of the 
top cell is the same as the number of steps, we obtain the following claim.

\begin{proposition}
\label{prop:CB_ranked}
For the decorated permutation $\pi^:\in CB_{kn}$ associated with a cell
$S_\M^\tnn$, we have $\dim S_\M^\tnn = k(n-k)- A(\pi^:)$.
The circular Bruhat order is a ranked poset with the rank
function $\mathrm{rank}(\pi^:) = k(n-k)-A(\pi^:)$.
\end{proposition}

Let us now explain the reason why we call the partial order $CB_{kn}$ 
the circular Bruhat order.
Actually, one can embed the usual (strong) Bruhat order on the symmetric group $S_k$ 
(and also the product of two copies of the Bruhat order) 
into $CB_{k,2k}$ as a certain interval.

Let $S_{\M}^\tnn$ be a cell in $Gr_{k,2k}^\tnn$ such that $[k],[k+1,2k]\in \M$.
Such cells are exactly the {\it double Bruhat cells\/} of Fomin-Zelevinsky, see
Remark~\ref{rem:FZ_double_Bruhat}.  Then the corresponding Bruhat necklace
$\I_\M$ has the entries $I_1=[k]$ and $I_{k+1} = [k+1,2k]$.  This means that
the corresponding permutation $\pi = \pi(\I_\M)$ satisfies the following
condition: $\pi:[k]\to [k+1,2k]$ and $\pi:[k+1,2k]\to [k]$.  (Such
permutations have no fixed points so there is no need to decorate them.)
In other words, $\pi$ can be subdivided into two permutations from $S_k$, as follows.
For two permutations $u,v\in S_k$, let $\pi=\pi(u,v)$ be the
permutation in $S_{2k}$, given by $\pi(i)=\overline{u(i)}$ and 
$\pi(\overline{i}) = v(i)$ for $i=1,\dots,k$, where $\overline{i} = 2k+1-i$.

\begin{proposition}  
For two permutations $\pi(u,v)$ and $\pi(u',v')$, where $u,v,u',v'\in S_k$
we have $\pi(u,v)\leq \pi(u',v')$ in the circular Bruhat order $CB_{k,2k}$
if and only if $u\leq u'$ and $v\leq v'$ in usual Bruhat order on $S_k$.
The interval $[\pi(1,1),\pi_{top}]$ in $CB_{k,2k}$ is isomorphic to the 
direct product of two copies of the usual Bruhat order on $S_k$.
\end{proposition}

\begin{proof}
For permutations of the form $\pi(u,v)$, the description of the circular Bruhat 
order from Corollary~\ref{cor:rab_pi} is equivalent to the well known
description of the usual Bruhat order on $S_k$:
$w\leq w'$ if and only if 
$|\{i\in[a]\mid w(i)\in[b]\}|\geq |\{i\in[a]\mid w'(i)\in[b]\}|$
for any $a,b\in[k]$.  The second claim follows from the fact
that any element element $\pi$ which is greater than $\pi(1,1)$ in $CB_{k,2k}$ 
should have the form $\pi=\pi(u,v)$.
\end{proof}

%

\section{Gluing of cells}
\label{sec:gluing}

In this section, we will explicitly describe how the nonnegative Grassmann cells 
$S_\M^\tnn$ are glued to each other using the network parametrization $\Mes_G$ 
of the cells.
\medskip


Let $G$ be a plabic (undirected) graph and let $\OO$ be a perfect orientation of
its edges.  Denote by $G'=(G,\OO)$ the corresponding directed graph.
According to Theorem~\ref{th:reduced_plabic_parametrization}
(or Theorem~\ref{th:Image_G_M_is_cell}),
for each cell $S_\M^\tnn\subset Gr_{kn}^\tnn$ there is a 
a directed graph $G'$ as above such that 
the boundary measurement map $\Mes_{G'}$
maps $\R_{>0}^{E(G')}/\{\text{gauge transformations}\}$ onto $S_\M^\tnn$.
By a slight abuse of notation, we will also denote by $\Mes_{G'}$ 
the map  $\R_{>0}^{E(G')}\to S_\M^\tnn$.

\begin{lemma} 
The map $\Mes_{G'}: \R_{>0}^{E(G')}\to Gr_{kn}^\tnn$ uniquely
extends to the continuous map 
$\overline{\Mes}_{G'}: \R_{\geq 0}^{E(G')}\to Gr_{kn}^\tnn$.
\end{lemma}

This claim is trivial for an acyclic graph $G'$. But we allow 
$G'$ to have cycles.

\begin{proof}  The uniqueness is clear because $\R_{>0}^{E(G')}$ is 
a dense subset in $\R_{\geq 0}^{E(G')}$.  We need to check that
$\Mes_{G'}$ does not have have poles as some of the edge variables 
$x_e$ approach $0$.
This follows from Proposition~\ref{prop:subtraction_free_minor}, 
which gives a rational subtraction-free expression for each maximal minor 
$\Delta_J(A)$ of the $k\times n$-matrix $A$ that represents 
$\Mes_{G'}(\{x_e\}_{e\in E(G')})$.
(This matrix has $\Delta_I(A)=1$ for the source set $I$ of $G'$.)
Note that the expression in the denominator contains the constant term $1$.
Thus we can specialize any subset of edge variables $x_e$ to $0$ without
getting a $0$ in the denominator.
\end{proof}

Clearly, the image of the map $\overline{\Mes}_{G'}$ belongs to 
the closure $\overline{S_\M^\tnn}$.  Moreover, this image
consists of the union of some cells $S_{\M''}^\tnn$.
Indeed, these are the cells that correspond 
(as in Theorem~\ref{th:Image_G_M_is_cell}) to all directed graphs 
$H'$ obtained from $G'$ by removing some edges.

Let us show that the opposite is true.
For any fixed (perfectly orientable) plabic graph $G$ and any point $p$ 
in the closure $\overline{S_\M^\tnn}$
of the corresponding cell, there exists a directed graph $G'$
obtained by a perfect orientation of edges of $G$ such that
$p\in \Image(\overline{\Mes}_{G'})$.

Let $p\in \overline{S_\M^\tnn}$.  Let us pick the graph $G'=(G,\OO)$
for some perfect orientation $\OO$ of $G$.
Then we can find nonnegative functions 
$x_e(t):]0,1]\to \R_{>0}$ on edges $e$ of $G'$ such that
$\lim_{t\to 0} \Mes_{G'}(\{x_e(t)\}) = p$
and each $x_e(t)$ is of the form
$x_e(t)= t^{m_e}f_e(t)$, where $m_e\in\R$ and 
$f_e(t)$ is real-valued analytic function such that $f_e(0)>0$.

\begin{lemma}
\label{lem:no_poles}  
For any collection of edge functions $x_e(t)$ as above 
one can apply gauge transformations at vertices and switch edge 
directions along some paths and/or closed cycles as in Section~\ref{sec:edge_switch}
to transform the $x_e(t)$ into functions $x_e'(t)$ that have no
poles at $t=0$, for all edges $e$. 
\end{lemma}

\begin{proof}
Let $M=\min(m_e)$.  We may assume that we have already transformed
the edge functions by gauge transformations and switching edge directions
so that $M$ is a big as possible and the number of edges $e$ such that
$m_e = M$ is as small as possible (for this $M$.)
If $M\geq 0$ then the $x_e(t)$ have no poles.
Assume that $M<0$.  
For any directed edge $e=(u,v)$ such that $m_e = M$, where
$v$ is an internal vertex, there should be an outgoing edge $e'$ from 
the vertex $v$ such that $m_{e'}=M$.
Otherwise, $m_{e'}> M$ for all edges $e'$ outgoing from $v$, and we could apply
a gauge transformation $t_v = t^{-\epsilon}$ at this vertex so that 
$m_e$ increases by $\epsilon$ and the $m_{e'}$  decrease by $-\epsilon$.  
For a sufficiently small $\epsilon>0$ this would make the number of edges
with $m_e=M$ smaller, which is impossible by our assumption.
Similarly, there should be an edge $e''$ incoming to the vertex $u$
such that $m_{e''} = M$.  This means that we can always find a directed
path $P$ in $G'$ joining two boundary vertices (or a closed cycle $C$)
such that for all edges $e$ in this path/cycle we have $m_e = M$.
Let us switch directions of edges in $P$ (or $C$) and invert the edge
functions $x_e(t)$ for these edges; see Section~\ref{sec:edge_switch}.
This switch would transform the $m_e=M$ into $-M$ for all edges in
the path/cycle, so again this would make the number of edges with $m_e = M$
smaller, which is impossible by our assumption.
Thus we should have $M\geq 0$, as needed.
\end{proof}

According to Lemma~\ref{lem:no_poles}, the graph $G'$ and the edge functions
$x_{e}(t)$ can be chosen so that the $x_e(t)$ have no poles at $t=0$.
That means that the point $p=\lim_{t\to 0}\Mes_{G'}(\{x_e(t)\})
=\overline{\Mes}_{G'}\{x_e(0)\}$ belongs to the image of 
$\overline{\Mes}_{G'}$.

We have proved the following result. 

\begin{theorem}
\label{th:closure_SM_O}
For a cell $S_\M^\tnn$, pick any (perfectly orientable)
plabic graph $G$ such that $\tilde \Mes_G$ maps onto $S_\M^\tnn$.
Then the closure of this cell equals
$$
\overline{S_\M^\tnn} = \bigcup_{G' = (G,\OO)} 
\overline{\Mes}_{G'}(\R_{\geq 0}^{E(G')})
\subseteq Gr_{kn}^\tnn,
$$
where $\OO$ ranges over all perfect orientations of $G$.
\end{theorem}

Assume that $G$ is a contracted plabic graph.
(That is $G$ is without unicolored edges, non-boundary leaves,
and  vertices of degree $2$).
For each $G'=(G,\OO)$, we have
$$
\overline{\Mes}_{G'}(\R_{\geq 0}^{E(G')})
=\bigcup_{H'} \Mes_{H'}(\R_{>0}^{E(G')}),
$$ 
where the union is over directed graphs $H'$ obtained from $G'$ by removing some 
edges (but keeping all vertices).
If we remove a directed edge $e=(u,v)$ where $u$ is a black vertex and $v$
is white, then $u$ will be an internal sink and $v$
will be an internal source in the obtained graph.
Instead of removing such edge $e$,
we could remove all incoming edges to the vertex $u$ (and get
a graph with the same image of $\Mes_{H'}$).
Thus in the above union it is enough to take only graphs $H'$ obtained
from $G'$ by removing some white-to-black directed edges $e=(u,v)$, 
i.e., $u$ is white and $v$ is black.
Notice that such graphs $H'$ will be perfectly oriented graphs with the same
coloring of vertices as in the graph $G$.

For a contracted plabic graph $G$, let us say that a plabic graph $H$ is a {\it
subgraph\/} of $G$, and write $H\subseteq G$, if $H$ is obtained from $G$ by
removing some edges
while keeping all vertices in $G$.  
When we remove a nonleaf edge $e=(b_i,v)$ attached to a boundary vertex,  
we need to create a new boundary leaf at $b_i$, whose color is opposite to $v$.  
When we remove an edge $e=(b_i,b_j)$ between two boundary vertices,
we need to create two boundary leaves at $b_i$ and $b_j$ of different colors.
(So there are two different ways to ``remove'' such edge.)
We are not allowed to remove boundary leaves.  
Note that each perfect orientation of a plabic subgraph $H\subseteq G$ 
uniquely extends to a perfect orientation $G$.  (Just direct all removed edges
from white vertices to black vertices.)

Let us call a plabic graph {\it almost-reduced\/} if it consists of
a reduced graph possibly together with some dipoles (isolated components
with a single edge and a pair of vertices of different colors).

Theorem~\ref{th:closure_SM_O} and the above discussion implies that the closure 
$\overline{S_\M^\tnn}$
is the union of $\Image(\tilde\Mes_H)$ over all 
perfectly orientable plabic subgraphs $H\subseteq G$.
Let us show that it is enough to take only almost-reduced plabic subgraphs $H$.

\begin{lemma}  Let $G$ be a perfectly orientable contracted plabic graph.
Then $G$ has an almost-reduced plabic subgraph $H\subseteq G$ such that
boundary measurement maps $\tilde \Mes_G$ and $\tilde \Mes_{H}$ 
have the same images in $Gr_{kn}^\tnn$.
\end{lemma}

Note that if we remove the edge from a dipole, then we would create
a graph with a singleton (isolated component with a single vertex),
which is not perfectly orientable.   This is why we need to 
consider almost-reduced subgraphs (and not just reduced ones).

\begin{proof}
Suppose that $G$ is not almost-reduced.  
Let us pick its perfect orientation $G'$ and transform this (directed) graph 
by the moves (M1)--(M3) into a (directed) graph $\tilde G$, where
we can apply a reduction (R1) or (R2); see Section~\ref{sec:transformations_plabic}.
Then we can find an edge $\tilde e=(u,v)$ in $\tilde G$ directed from a 
white vertex $u$ to a black vertex $v$ such that by removing $\tilde e$ 
from $\tilde G$ we will not change 
the image of $\Mes_{\tilde G}$.  (We need to check $3$ possible directions
of edges in (R1); for (R2) there is only one possible direction of edges;
see Figures~\ref{fig:edge_reduction} and~\ref{fig:leaf_reduction}.)
Let us transform $\tilde G$ back to $G'$ by  moves (M1)--(M3) and  keep 
track of the ``marked for removal'' edge $\tilde e$. 
In all cases when we perform a move $G''\mapsto G'''$
and $G''$ has a white-to-black marked edge $e''$, one can find a white-to-black
edge $e'''$ in $G'''$ such that the maps $\Mes_{G''\setminus\{e''\}}$ and 
$\Mes_{G'''\setminus\{e'''\}}$ have the same image.
The only case that needs a special attention is square move (M1).
Note that transformations of edge weights on Figure~\ref{fig:square_move_arrow} 
involve only weights of the four edges in the square.  Thus, if $e''$ is not one 
of these four edges, then we can just keep this edge $e'''=e''$.  
If $e''$ is one of the four edges in the square, the we can pick $e'''$
to be the edge opposite to $e''$ in the square.

This shows that in the (directed) graph $G'$ we can always find a white-to-black
edge $e'$ whose removal does not change the image of the map $\Mes_G$.
Thus in the undirected graph $G$ we can find an edge $e$ whose removal
preserves perfect orientability and does not change the image of $\tilde\Mes_G$.
If the graph $G\setminus\{e\}$ is still not almost-reduced then we can repeatedly
remove its edges using this procedure until we get an almost-reduced graph with the 
same image of the boundary measurement map.
\end{proof}


We proved the following result.

\begin{theorem}
\label{th:reduced_closure}
For a cell $S_\M^\tnn$, and any reduced contracted plabic graph $G$ associated with
this cell, we have
$$
\overline{S_\M^\tnn} = \bigcup_{H\subseteq G} \tilde 
\Mes_H(\R_{>0}^{F(H)-1}),
$$
where $H$ ranges over all almost-reduced plabic subgraphs $H\subseteq G$.
\end{theorem}

Note that each term in the right-hand side of this expression is a certain cell
$S_{\M'}^\tnn$ inside the closure $\overline{S_\M^\tnn}$.  Each map
$\tilde\Mes_H$ gives a parametrization of this cell. 

\begin{remark}
Let us describe bases of the matroid $\M$ using 
Theorem~\ref{th:reduced_closure}.
Any almost-reduced
subgraph $H\subseteq G$ should contain at least one edge incident to 
each internal vertex $v$ of $G$.  (If we remove all edges incident to $v$ 
then we would get a graph, which is not perfectly orientable.)
Thus {\it minimal\/} almost-reduced subgraphs are the subgraphs
containing exactly one edge at each internal vertex.  These are 
partial matchings of $G$, discussed in the end of Section~\ref{sec:plabic}.
Each minimal subgraph $H\subseteq G$ gives a cell $S_{\M'}^\tnn$ that
consists of a single point and $\M'$ contains a single base 
$I=\{i\mid b_i \text{ connected to a white boundary leaf}\}$.
For any base of $\M$, there is a minimal subgraph of this form.
We obtain the description $\M$ equivalent to the matching
matroid $\M_G^m$ from Lemma~\ref{lem:path_matchings}.
\end{remark}

\begin{figure}[ht]
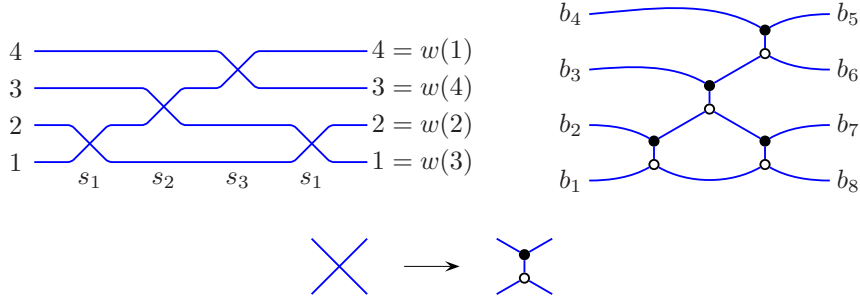

\psset{unit=0.7pt}
\pspicture(-10,-20)(220,70)
\psline[linearc=5]{-}(0,00)(20,00)(40,20)(60,20)(80,40)(100,40)(120,60)(140,60)(160,60)(180,60)
\psline[linearc=5]{-}(0,20)(20,20)(40,00)(60,00)(80,00)(100,00)(120,00)(140,00)(160,20)(180,20)
\psline[linearc=5]{-}(0,40)(20,40)(40,40)(60,40)(80,20)(100,20)(120,20)(140,20)(160,00)(180,00)
\psline[linearc=5]{-}(0,60)(20,60)(40,60)(60,60)(80,60)(100,60)(120,40)(140,40)(160,40)(180,40)
\rput(-10,0){$1$}
\rput(-10,20){$2$}
\rput(-10,40){$3$}
\rput(-10,60){$4$}
\rput(210,0){$1=w(3)$}
\rput(210,20){$2=w(2)$}
\rput(210,40){$3=w(4)$}
\rput(210,60){$4=w(1)$}
\rput(30,-10){$s_1$}
\rput(70,-10){$s_2$}
\rput(110,-10){$s_3$}
\rput(150,-10){$s_1$}
\endpspicture
\pspicture(-100,-10)(120,100)
\rput(-30,0){$b_1$}
\rput(-30,30){$b_2$}
\rput(-30,60){$b_3$}
\rput(-30,90){$b_4$}
\rput(120,90){$b_5$}
\rput(120,60){$b_6$}
\rput(120,30){$b_7$}
\rput(120,00){$b_8$}
\cnode[linewidth=0](-20,0){0}{B1}
\cnode[linewidth=0](-20,30){0}{B2}
\cnode[linewidth=0](-20,60){0}{B3}
\cnode[linewidth=0](-20,90){0}{B4}
\cnode[linewidth=0](110,90){0}{B5}
\cnode[linewidth=0](110,60){0}{B6}
\cnode[linewidth=0](110,30){0}{B7}
\cnode[linewidth=0](110,0){0}{B8}
\vwhite(15,8.66){V1}
\vblack(15,21.34){V2}
\vwhite(45,38.66){V3}
\vblack(45,51.34){V4}
\vwhite(75,8.66){V5}
\vblack(75,21.34){V6}
\vwhite(75,68.66){V7}
\vblack(75,81.34){V8}
\nccurve[angleA=00,angleB=-150]{-}{B1}{V1}
\nccurve[angleA=00,angleB=150]{-}{B2}{V2}
\nccurve[angleA=5,angleB=150]{-}{B3}{V4}
\nccurve[angleA=5,angleB=150]{-}{B4}{V8}
\nccurve[angleA=90,angleB=-90]{-}{V1}{V2}
\nccurve[angleA=-30,angleB=-150]{-}{V1}{V5}
\nccurve[angleA=30,angleB=-150]{-}{V2}{V3}
\nccurve[angleA=90,angleB=-90]{-}{V3}{V4}
\nccurve[angleA=-30,angleB=150]{-}{V3}{V6}
\nccurve[angleA=30,angleB=-150]{-}{V4}{V7}
\nccurve[angleA=90,angleB=-90]{-}{V7}{V8}
\nccurve[angleA=90,angleB=-90]{-}{V5}{V6}
\nccurve[angleA=30,angleB=180]{-}{V8}{B5}
\nccurve[angleA=-30,angleB=180]{-}{V7}{B6}
\nccurve[angleA=30,angleB=180]{-}{V6}{B7}
\nccurve[angleA=-30,angleB=180]{-}{V5}{B8}
\endpspicture

\pspicture(0,-5)(130,50)
\psline[linearc=5]{-}(0,0)(30,30)
\psline[linearc=5]{-}(0,30)(30,0)
\psline[linecolor=black,linewidth=.3pt]{->}(50,15)(80,15)

\vwhite(115,8.66){V1}
\vblack(115,21.34){V2}
\cnode[linewidth=0](100,0){0}{B1}
\cnode[linewidth=0](130,0){0}{B2}
\cnode[linewidth=0](100,30){0}{B3}
\cnode[linewidth=0](130,30){0}{B4}
\nccurve[angleA=30,angleB=-150]{-}{B1}{V1}
\nccurve[angleA=150,angleB=-30]{-}{B2}{V1}
\nccurve[angleA=90,angleB=-90]{-}{V1}{V2}
\nccurve[angleA=-30,angleB=150]{-}{B3}{V2}
\nccurve[angleA=-150,angleB=30]{-}{B4}{V2}
\endpspicture
\caption{Transforming a wiring diagram of $w=s_1 s_2 s_3 s_1$ into 
a plabic graph $G$}
\label{fig:wiring_plabic}
\end{figure}

\begin{remark}
There is an analogy between the usual Bruhat order on a Weyl
group $W$ and the circular Bruhat order $CB_{kn}$, as follows.  
The cells $S_\M^\tnn$ are analogues of Weyl group elements $w\in W$;
plabic graphs $G$ are analogues of reduced decompositions
$w=s_{i_1}\cdots s_{i_l}$; plabic subgraphs $H\subseteq G$
are analogues of reduced subwords in a reduced decomposition.
Then Theorem~\ref{th:reduced_closure} corresponds to the following
well-known description of the Bruhat order on $W$.
For $w=s_{i_1}\cdots s_{i_l}$, all elements $u$ such that  $u\leq w$
(in the Bruhat order) have reduced decompositions 
obtained by taking subwords in $s_{i_1}\cdots s_{i_l}$.

This analogy can be made more precise.  Each type $A$ reduced decomposition 
$w=s_{i_1}\dots s_{i_l}\in S_k$ is graphically represented by a wiring
diagram.  We can transform the wiring diagram into a reduced plabic graph $G$
with $n=2k$ boundary vertices by replacing each crossing with a pair of
trivalent vertices as shown on Figure~\ref{fig:wiring_plabic}.  The trip
permutation of $G$ is $\pi:i\to 2k+1-w(i)$, for $i\in[k]$ and $\pi:i\to
2k+1-i$, for $i\in[k+1,2k]$.  Subgraphs $H\subseteq G$ obtained by removing
some vertical edges of $G$ correspond to subwords in the reduced decomposition.
\end{remark}



\begin{remark}
Note that two different subgraphs of $H\subset G$ might correspond to 
the same cell in the closure $\overline{S_\M^\tnn}$. 
(Similarly, to the situation when different two subwords in a reduced 
decomposition give the same Weyl group element.)
That means that possibly there are some nontrivial identifications of 
components in the right-hand side of expression in Theorem~\ref{th:reduced_closure}.
This why the description of the geometrical structure of the closure
$\overline{S_\M^\tnn}$ is a nontrivial problem.
Conjecturally, the closure $\overline{S_\M^\tnn}$ 
is homeomorphic to an open ball.
\end{remark}

Let us now describe the covering relation $\lessdot$ 
in the circular Bruhat order using Theorem~\ref{th:reduced_closure}.
Suppose that the cell $S_\M^\tnn$ covers $S_{\M'}^\tnn$.
Let us pick a reduced contracted plabic graph $G$ such that  
$S_\M^\tnn = \Image(\tilde\Mes_G)$ and $G$ has no leaves and vertices of degree $2$.
By Theorem~\ref{th:reduced_closure}, there is an almost-reduced
subgraph $H\subset G$ such that $S_{\M'}^\tnn = \Image(\tilde\Mes_H)$. 
Note that if $H$ is obtained from $G$ by removing two or more edges,
then it number of faces drops by at least $2$.
By Proposition~\ref{prop:CB_ranked}, the circular Bruhat order $CB_{kn}$ 
is a ranked poset
with the rank function equal $\dim S_\M^\tnn$.  Thus the codimension
of $S_{\M'}^\tnn$ in $\overline{S_{\M}^\tnn}$ should be $1$, that is 
$H$ is obtained from $G$ by removing a {\it single\/} edge.
In this case we cannot create a dipole, so $H$ should be reduced.

Let us call an edge $e$ in a reduced contracted plabic graph $G$ {\it removable\/}
if $G\setminus\{e\}\subset G$ is a reduced plabic graph.
By Theorem~\ref{th:reduced_closure}, removable edges in $G$ are in one-to-one
correspondence with the cells covered by $S_\M^\tnn = \Image(\tilde \Mes_G)$.

Let $G$ be a reduced contracted plabic graph with the decorated trip permutation
$\pi^:=\pi^:(G)$.  For an edge $e$, let $T_1:b_{i}\to b_{\pi(i)}$
and $T_2:b_j\to b_{\pi(j)}$ be the two trips in $G$ that contain $e$ (and
pass this edge in two different directions).

\begin{lemma}  
The edge $e$ is removable if and only if the pair $(i,\pi(i))$
and $(j,\pi(j))$ is a simple crossing in the decorated trip permutation $\pi^:$.
In this case the decorated trip permutation of $G\setminus\{e\}$ is obtained
from $\pi^:$ by replacing this simple crossing with the corresponding
alignment; see Figure~\ref{fig:cover_rel}.  
\end{lemma}

Since we have already described covering relations in
Theorem~\ref{th:circ_bruhat_cover}, this lemma follows.  We can also easily
deduce it from the reducedness criterion in Theorem~\ref{th:reduced_criterion}. 

\begin{proof}
The trips of  $G\setminus\{e\}$ are exactly the same as the trips
of $G$, except that we need to switch tails of $T_1$ and $T_2$ at
the their (essential) intersection point at $e$.
According to Theorem~\ref{th:reduced_criterion}, 
the graph $G\setminus \{e\}$ is reduced if and only if the trips
$T_1$ and $T_2$ have only one essential intersection at $e$;
and there is no any other trip $T_3$ that intersects
the part of $T_1$ before $e$ and then the part of $T_2$ after $e$;
and vise versa.  The means that the pair $(i,\pi(i))$ and $(j,\pi(j))$
is a simple crossing in $\pi^:$.
On the other hand, for a simple crossing in $\pi^:$,
the corresponding pair of trips should intersect only once.
Otherwise, if they intersect $\geq 3$ times, then there is another trip
$T_3$ that passes through, say, the second intersection point of $T_1$ and $T_2$.
This trip cannot intersect the tails of $T_1$ and $T_2$ (after all their
intersection points).  Thus $T_3$ should end at a boundary point between 
$b_{\pi(j)}$ and $b_{\pi(i)}$.  
(Here we assume that the vertices are arranged on the circle as
in the crossing on Figure~\ref{fig:cross_align} and the word ``between''
means ``between in the clockwise order''.)
Similarly, $T_3$ should start at a
boundary point between $b_j$ and $b_i$. 
So the crossing in $\pi^:$ would not be simple.  Also if there is a trip $T_3'$ 
that first intersects with the initial part of $T_1$ before $e$ and then intersects
with the part of $T_2$ after $e$, then again the trip $T_3'$ would
form an obstruction for a simple crossing.
\end{proof}

\begin{corollary}  For a cell $S_\M^\tnn$ and a reduced contracted plabic graph $G$ 
such that $S_\M^\tnn=\Image(\tilde\Mes_G)$, the cells 
$S_{\M'}^\tnn$ that are covered by $S_\M^\tnn$ are in one-to-one correspondence
with removable edges $e$ of $G$.  They have the
form $S_{\M'}^\tnn = \Image(\tilde\Mes_{G\setminus\{e\}})$.
\end{corollary}

Let us show how to glue such adjacent cells $S_\M^\tnn$ and $S_{\M'}^\tnn$ together.
For a reduced plabic graph $G$ and a removable edge $e\in G$, pick a
perfect orientation $\OO$ of $G\setminus\{e\}$.  Let $G'$ be the graph obtained
by directing edges of $G$ so that the edge $e$ is directed from a white vertex 
to a black vertex and other edges are directed as in $\OO$.

Let $E=E(G')$ and $V=V(G')$, be the edge and vertex sets of $G'$.
Let $\R_{> 0}^{E-e}\times \R_{\geq 0}$ be the space of edge 
weights $\{x_f\}_{f\in E(G')}$ such that $x_f>0$, for $f\ne e$, and $x_e\geq 0$.
Then
$$
(\R_{> 0}^{E-e}\times \R_{\geq 0})/\{\text{gauge transformations}\}
\simeq \R_{>0}^{|E|-|V|-1}\times \R_{\geq 0}.
$$

\begin{corollary}
The map  $\overline{\Mes}_{G'}:\R_{\geq 0}^{E}\to \overline{S_\M^\tnn}$ 
induces the bijective map
$$
\R_{>0}^{|E|-|V|-1}\times \R_{\geq 0} \to \overline{S_\M^\tnn}.
$$
The restriction of this map to the subset given $x_e>0$ is
a parametrization of $S_\M^\tnn$ and the restriction of this map to 
the subset $x_e=0$ is a parametrization of $S_{\M'}^\tnn$.
\end{corollary}

\section{\protect\Le-diagrams and Bruhat intervals}
\label{sec:Bruhat_intervals}

Rietsch constructed cellular decomposition of 
the totally nonnegative part of $G/P$; 
see \cite{Riet1, Riet2, M-R}.  
In this section we show that Rietsch's cells are
in one-to-one correspondence with the cells $S_\M^\tnn$.

Let us fix the pair $(k,n)$ as before.
Recall that a permutation $w\in S_n$ is called {\it
Grassmannian\/} if $w$ is the minimal length representative
of a left coset $(S_k\times S_{n-k})\backslash S_n$.
In other words, $w$ is a Grassmannian permutation if
$w(1)<w(2)<\cdots <w(k)$ and $w(k+1)<w(k+1)<\cdots <w(n)$.

In case of the Grassmannian $Gr_{kn}$, Rietsch's cells are labeled by pairs of
permutations $(u,w)$ in $W=S_n$ such that $u\leq w$ in the Bruhat order on $W$
and $w$ is a Grassmannian permutation.  Let us construct a bijection between
such pairs and $\Le$-diagrams.

\begin{figure}[ht]
\input{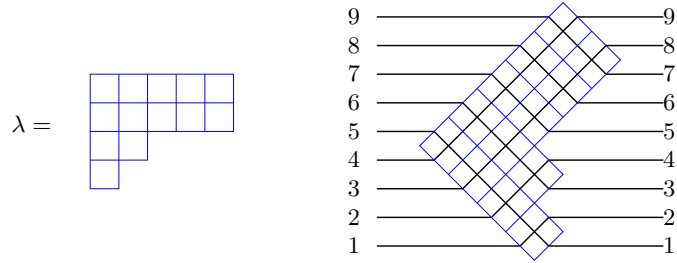}
\caption{A Young diagram $\lambda=(5,5,2,1)$ and the corresponding Grassmannian 
permutation $w_\lambda = 2\, 4\, 8\, 9\, 1\, 3\, 5\, 6\, 7$
for $(k,n) = (4,9)$}
\label{fig:grass_perm}
\end{figure}

Grassmannian permutations are in bijection with partitions $\lambda\subseteq
(n-k)^k$.  This bijection can be described as follows. 
Rotate the Young diagram of shape $\lambda$ by
$45^\circ$ counterclockwise and draw the wiring diagram with $k$ wires going
along the rows of $\lambda$ and $(n-k)$ wires going along the columns of
$\lambda$. 
Label both ends of the wires by numbers $1,\dots,n$ starting from 
the bottom.
This wiring diagram represents a permutation $w_\lambda$ such that
the wires connect the indices $i$ on the left with the $w_\lambda(i)$ on the right;
see Figure~\ref{fig:grass_perm}.  
Using our earlier notation, we can describe this permutation as
$w_\lambda = (\tilde i_k,\tilde i_{k-1},\dots,\tilde i_1,
\tilde j_{n-k},\tilde j_{n-k-1},\dots,\tilde j_{1})$, 
where $I(\lambda)=\{i_1<\cdots<i_k\}$,
$[n]\setminus I(\lambda) = \{j_1<\cdots< j_{n-k}\}$, and
$\tilde i := n+1-i$;
see Section~\ref{sec:grassmannian}.
It is clear that the length $\ell(w_\lambda)$ (the number of inversions)
equals $|\lambda|$.

Let $D$ be a \Le-diagram of shape $\lambda$. 
Again rotate it by $45^\circ$ counterclockwise.\footnote{After
rotation, it should be called $>$-diagram.}
For each box of $\lambda$ filled with a $1$ in $D$ (shown by a dot),
we replace the corresponding crossing in the wiring diagram of $w_\lambda$
by an uncrossing, as shown on Figure~\ref{fig:grass_Le_diag}.
The obtained ``pipe dream'' is a wiring diagram of a certain permutation in $S_n$,
which we denote by $u_D\in S_n$.

\begin{figure}[ht]
\input{grass_Le_diag.pstex_t}
\caption{A \Le-diagram $D$ and the corresponding permutation 
$u_D=1\,4\,2\,7\,3\,5\,9\,6\,8$}
\label{fig:grass_Le_diag}
\end{figure}

\begin{theorem}
\label{th:Le-bruhat}      
The map $D\mapsto u_D$ is a bijection between \Le-diagrams of shape
$\lambda$ and permutations $u\in S_n$ such that $u\leq w_\lambda$ 
in the Bruhat order on $S_n$.  The number of $1$'s in $D$ equals
$\ell(w_\lambda)-\ell(u_D)$.
\end{theorem}

\begin{lemma} 
\label{lem:distinguish_subword}
{\rm Marsh-Rietsch~\cite[Lemma~3.5]{M-R}}
Let $w=s_{i_1}\cdots s_{i_l}$ be a reduced decomposition of a
Weyl group element $w\in W$.
Then, for any $u$ such that $u\leq w$ in the Bruhat order,
there is unique subset $J=\{j_1<\cdots <j_s\}\subseteq [l]$ 
such that $u= s_{i_{j_1}}\cdots s_{i_{j_s}}$ and
$(s_{i_1}' s_{i_2}'\cdots s_{i_{a}}') \,s_{i_{a+1}} \geq
s_{i_1}' s_{i_2}'\cdots s_{i_{a}}'$, for all $a\in[l-1]$,
where $s_{i_j}'= s_{i_j}$ for $j\in J$ and $s_{i_j}'=1$ otherwise.
\end{lemma}

Clearly, $u= s_{i_{j_1}}\cdots s_{i_{j_s}}$ should be a reduced decomposition.

It is well-known that for any $u\leq w$ there is a subword in a reduced
decomposition $w=s_{i_1}\cdots s_{i_l}$ which gives a reduced decomposition of
$u$; see~\cite{Hump}.  However, there are usually many such subwords.  The
above lemma describes one distinguished subword for each $u\leq w$.
Remark that the subset $J$ described in the Lemma~\ref{lem:distinguish_subword}
is exactly the {\it lexicographically minimal\/} subset in $[l]$ that produces 
a reduced decomposition of $u$.

In type $A$ case, the condition of Lemma~\ref{lem:distinguish_subword}
have a simple combinatorial description.  A reduced decomposition of $w$
corresponds to a wiring diagram $\mathcal{W}$ where the wires are not 
allowed to cross more than once.  Subwords in the reduced decomposition 
correspond to diagrams $\mathcal{W}'$ obtained from $\mathcal{W}$ by 
replacing some crossings of with uncrossings (like we did above for $w_\lambda$).
The condition of Lemma~\ref{lem:distinguish_subword} says that once 
two wires in $\mathcal{W}'$ intersect with each other they can never intersect
or even {\it touch} each other again.   Here ``touch'' means that the two 
wires participate in the same uncrossing of $\mathcal{W}'$.

For example, this condition fails for the wiring diagram shown on
Figure~\ref{fig:grass_Le_diag}.  Indeed, the two wires whose left ends are
labeled by $4$ and $6$ intersect each other and then arrive to the same
uncrossing.  However, if we take the mirror image of the condition (with 
respect to the vertical axis) then we will get exactly what we need.

Let $\MR(\lambda)$ be the set of pipe dreams obtained from the wiring
diagram of $w_\lambda$ by replacing some crossings with uncrossings
so that the following conditions hold:
\begin{enumerate}
\item Two wires can intersect at most once. 
\item If two wires do intersect at point $P$, then they cannot
participate in the same uncrossing to the left of $P$.
\end{enumerate}
The set $\MR(\lambda)$ corresponds to the subwords in a reduced decomposition
of $w_\lambda$ that satisfy the (mirror image) of Marsh-Rietsch's condition.

\begin{lemma}
The set $\MR(\lambda)$ is exactly the set of pipe dreams coming from 
\Le-diagrams, as described above.
\end{lemma}

\begin{proof}
Recall that \Le-diagrams $D$ can be described as follows.  For any 
box $x$ filled with a $0$ in $D$, either all boxes above $x$ or
all boxes below $x$ are filled with $0$'s.
Let us translate this \Le-condition in the language of pipe dreams.
It says that, for any crossing of two wires at point $P$, at least
one of these wires does not participate in any uncrossing before $P$.
So that one of these wires goes directly from the boundary of 
$\lambda$ (from the North-West or from the South-West) to the 
intersection point $P$ without 
making any turns.
Clearly, such pipe dreams belong to the set $\MR(\lambda)$.
Let us now show the opposite inclusion.
Suppose that the \Le-condition fails at some point $P$.
That means that both wires intersecting at $P$ 
diverge from the straight course before $P$.
Let $A$ and $B$ be the points where the wires make the last turns before 
arriving to $P$; see Figure~\ref{fig:ABCP}.
\begin{figure}[ht]
\input{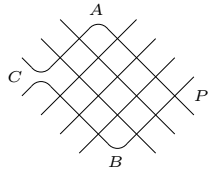}
\caption{A failed \Le-condition implies a failed MR-condition}
\label{fig:ABCP}
\end{figure}
Let us consider the rectangle $R$ with the vertices $A,B,P$
and another vertex $C$.  
Let us assume that our failed \Le-condition was chosen so that
the size of the rectangle $R$ (say, its perimeter) 
is a small as possible.
Then any other wire that intersects with 
the side $AP$ or $BP$ of $R$ cannot make any turn in the rectangle $R$.
Indeed, if it makes a turn then it gives another failed \Le-condition
with a smaller rectangle.  Thus all other wires go straight through
the rectangle $R$ as shown on Figure~\ref{fig:ABCP}.
Thus means that the two wires that intersect at $P$ should either
intersect or touch each other at $C$.  This means that the condition 
describing the set $\MR(\lambda)$ also fails.
\end{proof}

This proves Theorem~\ref{th:Le-bruhat} and shows that our cells are in
bijection with Rietsch's cells.  This implies Theorem~\ref{th:Rietsch=equal}
saying that these two cell decompositions coincide.

\section{From \protect\Le-diagrams to decorated permutations (and back)}

The cells $S_\M^\tnn$ are in bijection with \Le-diagrams
(Theorem~\ref{th:g_D}) and also in bijection with 
decorated permutations (Theorem~\ref{th:SM_dec_perm}).
In this section we discuss the induced bijection between 
\Le-diagrams and decorated permutations.
\medskip

\begin{figure}[ht]
\input{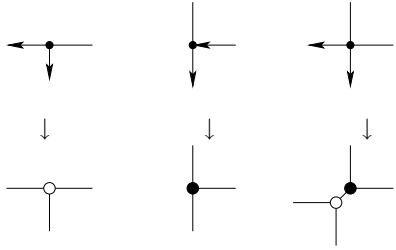}
\caption{Transforming $\Gamma$-graphs into plabic graphs}
\label{fig:gamma_plabic}
\end{figure}

Let us pick a \Le-diagram $D$ of shape $\lambda\subseteq (n-k)^k$, transform it
into $\Gamma$-graph $G^\Gamma_D$, and then transform this graph into a plabic
graph $G_D^{plabic}$ (see  Section~\ref{sec:perfection}).  In other words, we
need replace $4$-valent vertices of $G^\Gamma_D$ by pairs of trivalent vertices
and color the vertices as shown on Figure~\ref{fig:gamma_plabic}.
Theorems~\ref{th:g_D} and~\ref{th:SM_dec_perm}, and
Proposition~\ref{prop:matorid_to_perm} imply that the corresponding 
decorated permutation is the decorated trip permutation of the obtained graph.

\begin{corollary}
The map $D\mapsto \pi^:(G_D^{plabic})$ is a bijection between 
\Le-diagrams of shape $\lambda \subseteq (n-k)^k$
and decorated permutations with the anti-exceedance set $I(\lambda)$.
\end{corollary}

\begin{figure}[ht]
\input{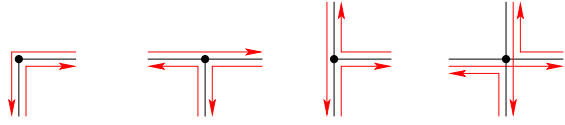}
\caption{Rules of the road for $\Gamma$-graphs}
\label{fig:gamma_rules}
\end{figure}

The rules of the road for plabic networks (see Figure~\ref{fig:trip_rules})
translate into the rules of the road for $\Gamma$-graphs shown on
Figure~\ref{fig:gamma_rules}.  So the decorated 
permutation $\pi^:$ corresponding to a \Le-diagram $D$
can be described as follows.  The empty rows (columns) or $D$
correspond to white (black) fixed points of $\pi^:$.
For other entries we need to follow trips in the graph $\Gamma$-graph 
$G_D^\Gamma$.  Let us trace {\it backwards\/} the
trip that ends at a boundary vertex $b_i$ located on a vertical segment of the
boundary.
We need to go from $b_i$ all the way to the left until we hit a $\Gamma$-turn
or a $\vdash$-fork (the first and the third segments shown on
Figure~\ref{fig:gamma_rules}), then turn down and go until the first junction,
then turn right and go until the first junction, then turn down, then right,
etc.  We need to keep going in this zig-zag fashion until we hit the boundary
at a boundary vertex $b_j$.  Then we should have $\pi(j) = i$ in the corresponding
permutation.
The rule for trips that end on a horizontal segment of the boundary 
is the symmetric to the above rule (with respect to the axis $x+y=0)$.

Steingrimsson-Williams~\cite{S-W} investigated nice properties of this
bijection $D\mapsto \pi^:(G_D^{plabic})$ with respect to various statistics on
permutations and \Le-diagrams.

Let us give a simpler description of the bijection between \Le-diagrams 
and decorated permutations using the bijection $D\mapsto u_D$
constructed in Section~\ref{sec:Bruhat_intervals}.

Let $w_\lambda = (w_1,\dots, w_n) = (\tilde i_k,\tilde i_{k-1},\dots,\tilde i_1,
\tilde j_{n-k},\tilde j_{n-k-1},\dots,\tilde j_{1})$ be the Grassmannian
permutation associated with $\lambda\subseteq (n-k)^k$, 
where $I(\lambda) = \{i_1<\cdots<i_k\}$,
$ [n]\setminus I(\lambda) = \{j_1<\cdots< j_{n-k}\} $, and $\tilde i := n+1-i$.

\begin{lemma}
\label{lem:uSn}
For a permutation $u\in S_n$, we have $u\leq w_\lambda$
in the Bruhat order if and only if 
$u_1\leq w_1,\dots, u_k\leq w_k$ and $u_{k+1}\geq w_{k+1},\dots,u_n\geq w_n$.
\end{lemma}

\begin{proof}
According to the well-known description of the Bruhat order on $S_n$,
we have $u\leq w_\lambda$ if and only if $\{i\in[a]\mid u(i)\in [b]\}
\leq \{i\in[a]\mid w_\lambda(i)\in [b]\}$, for any $a,b\in[n]$.
This translates into the needed inequalities.
\end{proof}

Let us now describe the map $u\mapsto \pi^:=\pi^:(u)$ from permutations
$u\leq w_\lambda$ to decorated permutations $\pi^:$ with the anti-exceedance 
set $I(\lambda)$.  It is given by $\pi^{-1}(i_1) = \tilde u_k$, \dots,
$\pi^{-1}(i_k) = \tilde u_1$ (with fixed points colored in white),
and $\pi^{-1}(j_1) = \tilde u_n$, \dots, $\pi^{-1}(j_{n-k}) = \tilde u_{k+1}$
(with fixed points colored in black).

\begin{theorem}
The map $u\mapsto \pi^:(u)$ is a bijection between 
the Bruhat interval $\{u\mid u\leq w_\lambda\}$ and decorated permutations with 
anti-exceedance set $I(\lambda)$.
The map $D\mapsto \pi^:(u_D)$ is the bijection between \Le-diagrams
and decorated permutations, which coincides with the above
map $D\mapsto \pi^:(G_D^{plabic})$.
\end{theorem}

\begin{proof}
The first claim follows directly from Lemma~\ref{lem:uSn}.
The second claim is obtained by comparing the constructions
of $u_D$ and $\pi^:(G_D^{plabic})$.
\end{proof}

\section{Cluster parametrization and chamber anzatz}

\section{Berenstein-Zelevinsky's string polytopes}

\section{Enumeration of nonnegative Grassmann cells}

Let $N_{kn}$ be the number of totally nonnegative cells in
the Grassmannian $Gr_{kn}$.
Recall that the {\it Eulerian number\/} $A_{kn}$ is the number of 
permutations $w\in S_n$ 
such that $w$ has $k-1$ descents:
$\#\{1\in[n-1] \mid w(i) > w(i+1)\} = k-1$.

\begin{proposition}
The numbers $N_{kn}$ are related to the Eulerian numbers by
$N_{kn} = \sum_{r=0}^n {n\choose r} A_{k,n-r}$.
Their generating function is
$$
1+\sum_{n\geq 1,\, 0\leq k\leq n} x^k {y^n\over n!}\, N_{kn}  =
e^{xy} {x-1\over x-e^{y(x-1)}}.
$$
\end{proposition}

\begin{proof}  
According to Theorem~\ref{th:SM_dec_perm},
$N_{kn}$ is the number of decorated permutations $\pi^:$ of size
$n$ with $k$ anti-exceedances.
If we remove black fixed points from $\pi^:$ we get a usual
permutation with $k$ anti-exceedances.  It is well-know that the Eulerian 
number count such permutations.
The second claim is obtained from the known generating
function for the Eulerian numbers.
\end{proof}

The numbers $N_{kn}$ appear in Sloan's On-Line Encyclopedia of Integer 
Sequences \cite{Slo} with ID number A046802.

\begin{figure}[ht]
\begin{tabular}{|l|cccccccc|}
\hline
$n\backslash k$  & \  0 & 1  & 2   & 3   & 4   & 5  &  6 & $\cdots$ \\
\hline
0 & \ 1 &    &     &     &     &    &    &  \\
1 & \ 1 & 1  &     &     &     &    &    &  \\
2 & \ 1 & 3  & 1   &     &     &    &    &  \\
3 & \ 1 & 7  & 7   & 1   &     &    &    &  \\
4 & \ 1 & 15 & 33  & 15  & 1   &    &    &  \\
5 & \ 1 & 31 & 131 & 131 & 31  & 1  &    &  \\
6 & \ 1 & 63 & 473 & 883 & 473 & 63 & 1  & \\
\vdots & \vdots&\vdots &\vdots  &\vdots  &\vdots  &\vdots  & \vdots & $\ddots$ 
\\
\hline
\end{tabular}
\caption{The numbers $N_{kn}$ of nonnegative cells in $Gr_{kn}^\tnn$}
\end{figure}

Let $N_n = \sum_{k=0}^n N_{kn}$ be the total number of decorated permutations
of size $n$.

\begin{proposition}
The numbers $N_n$ satisfy the recurrence relation
$N_n = n\cdot N_{n-1} +1$, $N_0=1$.
The exponential generating function for these numbers is
$\sum_{n\geq 0} N_n \frac{x^n}{n!} = e^x/(1-x)$.  
\end{proposition}

The sequence $N_n$ appears in~\cite{Slo} with ID number A000522.

Let $N_{kn}(q) =  \sum q^{\dim S_\M^\tnn}$ be the generating function for 
the nonnegative cells $S_\M^\tnn\subset Gr_{kn}^\tnn$ counted according to 
their dimension.  By Theorem~\ref{th:g_D}, we have
$$
N_{kn}(q) = \sum_{D} q^{|D|},
$$
where the sum is over \Le-diagrams whose shape
fits inside the rectangle $(n-k)^k$, and $|D|$ denotes the 
number of 1's in the diagram.

Williams gave a formula for the polynomials $N_{kn}(q)$ by 
counting \Le-diagrams.

\begin{theorem} {\rm Williams~\cite[Theorem~4.1]{W1}} \ 
We have
$$
N_{kn}(q) = \sum_{i=1}^{k-1} \binom{n}{i} q^{-(k-i)^2}\,
\left( [i-k]_q^i \, [k-i+1]_q^{n-i} - [i-k+1]^i\, [k-i]^{n-i}_q\right),
$$
where $[i]_q:= \frac{1-q^i}{1-q}$.
\end{theorem}

Steingrimsson-Williams \cite{S-W} studied various statistics on
\Le-diagrams.  Corteel-Williams \cite{CW} investigated the relationship
between \Le-diagrams and the asymmetric exclusion process.

\section{\protect\miscell}

Recall that a {\it rook placement\/} on some board
is a way to place rooks so that no rook attacks another rook.

Define a {\it teuton\/} as a chesspiece that can move 
only downwards and to the right. 
In order to attack a piece, two teutons need to simultaneously charge 
from two different 
directions in a wedge-shaped formation.\footnote{According to Tacitus,
      the ancient Teutons arranged their forces in the form of a wedge 
      ({\it Germania,} 6).
      Later this wedge-shape phalanx was known as {\it sv\'infylking.}}
A {\it teuton placement\/} 
is a placement of several teutons on a board such
that no pair of teutons can attack another teuton.
In other words, a teuton placement is a subset $S$ (possibly empty) of boxes 
on some board (say, a Young diagram) such that for any three boxes $a,b,c$
in a \Le-shaped pattern (as on Figure~\ref{fig:Le_property}),
if $a,c\in S$ then $b\not\in S$.

For $\lambda\subseteq(n-k)^k$, let $\{i_1<\cdots<i_k\}=I(\lambda)$,
$\{j_1<\dots<j_{n-k}\}=[n]\setminus I(\lambda)$.
Define the skew shape 
$\kappa_\lambda:=(n^{n-k},\tilde i_1,\dots,\tilde i_{k})/
(\tilde j_1-1,\dots,\tilde j_{n-k}-1)$, where $\tilde i:=n+1-i$.

For a permutation $w\in S_n$, define the hyperplane arrangement $A_w$
in $\R^n$ that consists of the hyperplanes $x_i - x_j =0$, for all
inversions $i<j$, $w(i)>w(j)$.

\begin{theorem}   The following numbers are equals:
\begin{enumerate} 
\item 
The number of nonnegative cells $S_\M^\tnn$ inside the Schubert cell 
$\Omega_\lambda$.
\item The number of \Le-diagrams of shape $\lambda$.
\item The number of decorated permutations with anti-exceedance set
$I(\lambda)$.
\item The number of teuton placements on the Young diagram of shape $\lambda$
\item The number of rook placements with $n$ rooks on the skew Young
diagram of shape $\kappa_\lambda$.
\item The number of permutations in the Bruhat interval 
$\{u\in S_n\mid u \leq  w_\lambda\}$.
\item The number of regions of the hyperplane arrangement $A_{w_\lambda}$.
\end{enumerate}
\end{theorem}

\begin{remark}  It seems that, for many permutations $w\in S_n$,
the number of regions in $A_{w}$ equals to the number of elements 
in the Bruhat interval $\{u\mid u\leq w\}$.
For example, this is true for the longest permutation $w_\circ\in S_n$.
However this is not true for the four permutations in $S_6$ with reduced 
decompositions $s_2 s_4 [s_1][s_5] s_3 s_2 s_4$, where the terms is the 
brackets might be omitted.  We suspect that the equality holds if and only if 
$w$ avoids the four patterns given by these permutations.
\end{remark}


\begin{theorem}  The number of \Le-diagrams $D$ of the triangular
shape $\lambda=(n,n-1,\dots,1)$ such that $D$ contains no 1's in
the $n$ corner boxes equals $n!$.
The bijection from this set of diagrams to permutations is constructed
as follows.   For $i\in[n]$, let $a_{i1}> a_{i2} > \dots > a_{i,k_i}$ 
be the positions of 1's in the $i$-th column of $D$.
The the permutation $w\in S_n$ corresponding to $D$ is given by 
the following product of cycles:
$$
w = (n,a_{11}, a_{12},\dots,a_{1,k_1})(n-1,a_{21},a_{22},\dots)
(n-2,a_{31},a_{32},\dots) \cdots (1).
$$
\end{theorem}

\end{document}